\documentclass[final,a4paper,11pt]{article}
\usepackage{mathptmx}
\usepackage{amssymb,amsmath, amsfonts,amsthm}
\usepackage{a4wide}
\usepackage{color}
\usepackage[colorlinks,citecolor=blue]{hyperref}
\usepackage[capitalise,nameinlink]{cleveref}
\usepackage[noadjust]{cite}
\usepackage{relsize}

\usepackage[shortlabels]{enumitem}
\setlist[enumerate]{label=(\alph*)}

\numberwithin{equation}{section}

\newcommand{\email}[1]{{\tt #1}}
\newcommand{\R}{\mathbb{R}}
\newcommand{\N}{\mathbb{N}}

\newcommand{\norm}[1]{\|#1\|}
\newcommand{\nnorm}[1]{{\vert\kern-0.25ex\vert\kern-0.25ex\vert #1 
    \vert\kern-0.25ex\vert\kern-0.25ex\vert}}
\newcommand{\nskalp}[1]{{\langle\kern-0.25ex\langle\kern-0.25ex\langle #1 
    \rangle\kern-0.25ex\rangle\kern-0.25ex\rangle}}

\newcommand{\dist}{\operatorname{dist}}
\newcommand{\bdry}{\operatorname{bdry}}
\newcommand{\diag}{\operatorname{diag}}

\newcommand{\vecmax}{\operatorname{vecmax}}

\newcommand{\mv}{\,\mid\,}

\newcommand{\Np}{\hat{\cal N}^p}

\newcommand{\skalp}[1]{\langle #1\rangle}

\newcommand{\xb}{\bar x}
\newcommand{\yb}{\bar y}
\newcommand{\zb}{\bar z}

\newcommand{\oo}{o}
\newcommand{\OO}{{O}}
\newcommand{\argmin}{{\rm arg\,min\,}}

\newcommand{\inn}{{\rm int\,}}

\newcommand{\gph}{\mathrm{gph}\,}
\newcommand{\dom}{\mathrm{dom}\,}
\newcommand{\epi}{\mathrm{epi}\,}
\newcommand{\dr}{{\rm d\,}}
\newcommand{\tto}{\rightrightarrows}

\newcommand{\mscLink}[1]{\href{http://www.ams.org/mathscinet/msc/msc2020.html?t=#1}{#1}}

\renewcommand\email[1]{\href{mailto:#1}{\texttt{#1}}}

\newcommand{\orcid}[1]{ORCID: \href{https://orcid.org/#1}{#1}}

\newtheorem{theorem}{Theorem}[section]
\newtheorem{proposition}[theorem]{Proposition}
\newtheorem{remark}[theorem]{Remark}
\newtheorem{lemma}[theorem]{Lemma}
\newtheorem{corollary}[theorem]{Corollary}
\newtheorem{definition}[theorem]{Definition}
\newtheorem{example}[theorem]{Example}

\title{
	Why second-order sufficient conditions are, in a way, easy\\
	-- or --\\
	revisiting calculus for second subderivatives\thanks{%
		This paper is dedicated to Roger Wets on the occasion of his 85th birthday.
	}%
	}

\author{
	Mat\'u\v{s} Benko\thanks{%
		Applied Mathematics and Optimization, 
		University of Vienna, 1090 Vienna, Austria,
		\email{matus.benko@univie.ac.at}.
		}%
	\and
	Patrick Mehlitz\thanks{%
		Institute of Mathematics, 
		BTU Cottbus-Senftenberg, 03046 Cottbus, Germany, 
		\email{mehlitz@b-tu.de},
		\orcid{0000-0002-9355-850X}.
		}%
	}

\date{}

\begin{document}

\maketitle

 \noindent
 {\bf Abstract.}\ 
 	In this paper, we readdress the classical topic of second-order sufficient
 	optimality conditions for optimization problems with nonsmooth structure.
 	Based on the so-called second subderivative of the objective function
 	and of the indicator function associated with the feasible set,
 	one easily obtains second-order sufficient optimality conditions
 	of abstract form.
 	In order to exploit further structure of the problem, e.g., composite
 	terms in the objective function or feasible sets given as (images of)
 	pre-images of closed sets under smooth transformations, to make these
 	conditions fully explicit, we study calculus rules for the second 
 	subderivative under mild conditions.
 	To be precise, we investigate a chain rule and a marginal function rule,
 	which then also give a pre-image and image rule, respectively.
 	As it turns out, the chain rule and the pre-image rule
 	yield lower estimates, desirable in order to obtain sufficient optimality conditions,
	for free.
 	Similar estimates for the marginal function and the image rule are valid 
 	under a comparatively mild inner calmness* assumption.
 	Our findings are illustrated by several examples including problems
 	from composite, disjunctive, and nonlinear second-order cone programming.



\section{Introduction}

The derivation of second-order necessary and sufficient optimality conditions has been
a key topic in mathematical optimization throughout the last decades. 
Clearly, using second-order information associated with the appearing
data provides sharper conditions than those ones obtained from pure first-order information.
Moreover, there is a huge interest in second-order sufficient
optimality conditions since these do not only guarantee stability of the
underlying local minimizer in certain sense, but can be applied in numerical
optimization in order to show local fast convergence of diverse solution methods.

Second-order optimality conditions for smooth standard nonlinear optimization problems (NLPs)
date back, e.g., to \cite{BenTal1980,McCormick1967} and essentially postulate positive
(semi-) definiteness of the Hessian of the associated Lagrangian function over a suitable
critical cone. We are interested in a much broader setting where the constraint region is given by so-called
geometric constraints, i.e., feasible sets of the form 
$F^{-1}(C):=\{x\in\R^n\,|\,F(x)\in C\}$ where 
$F\colon\R^n\to\R^m$ is a twice continuously differentiable mapping and $C\subset\R^m$
is a closed set.
Second-order optimality conditions typically combine the second-order derivative of the objective function
and the curvature of the feasible set $F^{-1}(C)$.
In case of standard NLPs, where the set $C$ is convex and \emph{polyhedral}, the curvature of the feasible set
comes solely from the constraint mapping $F$ and is incorporated in the Hessian of the Lagrangian.
As shown e.g.\ in \cite{Gfrerer2014,Mehlitz2019b}, the situation is analogous for so-called disjunctive programs,
where $C$ is still polyhedral, i.e., the union of finitely many convex polyhedral sets.
In general, however, the curvature of the feasible set comes from both, mapping $F$ and set $C$.
Exemplary, in nonlinear second-order cone programming, $C$ corresponds to the well-known second-order cone, i.e., a curved, closed, convex set,
and second-order optimality conditions additionally comprise a so-called curvature term associated with $C$, see \cite{BonnansRamirez2005}.
Moreover, it has been suggested in \cite{BeGfrYeZhouZhang1} that various rather challenging problem classes
like bilevel optimization problems, see e.g.\ \cite{Dempe2002}, or optimization
problems with (quasi-) variational inequality constraints, 
see e.g.\ \cite{FacchneiPang2003,LuoPangRalph1996,OutrataKocvaraZowe1998},
possess certain common structure where $C$ itself is quite complicated
and so the computation of its curvature is far from trivial.

Various variational tools have been used to capture the curvature of the feasible set.
Among the standard ones are the support function, applied to a suitable (approximation of some) second-order tangent set to the feasible set,
see \cite{BonSh00,BonnansCominettiShapiro1999},
and the second subderivative, applied to the indicator function associated with the feasible set,
see the classical references \cite{Rockafellar1989,RoWe98}, \cite[Section~7]{MohammadiMordukhovichSarabi2021} and \cite{ThinhChuongAnh2021} for some recent developments,
and \cite{ChristofWachsmuth2018,WachsmuthWachsmuth2022} for generalizations of this approach to optimization in abstract spaces.
We believe that these two tools are more suitable for \emph{sufficient} conditions
since the related second-order sufficient conditions in terms of $F$ and $C$ can be derived without any additional requirements, such as constraint qualifications.
This fact was shown in \cite[Section~3]{BeGfrYeZhouZhang1}, and we further clarify it in this paper from the perspective of calculus rules.
On the other hand, to derive the respective necessary conditions, one typically needs rather severe assumptions,
see e.g.\ \cite[Theorem~4]{GfrererYeZhou2022}. 
A new construction, the so-called lower generalized support function, was recently introduced in \cite{GfrererYeZhou2022} as another tool to measure the curvature of the feasible set.
The resulting second-order necessary conditions in terms of $F$ and $C$ are valid under very mild assumptions (so-called directional metric subregularity), see \cite[Theorem~2]{GfrererYeZhou2022}. Again, this approach seems to be less suitable for sufficient conditions.
We point out that the replacement of the standard support function (particularly useful in the convex setting) by the lower generalized support function
to some extent resembles the move from the convex subdifferential to the limiting one, which is mainly
suitable for (first-order) \emph{necessary} optimality conditions.
All these three tools and approaches are currently being closely investigated in a larger context of second-order variational analysis in \cite{BeGfrYeZhouZhang1} and a forthcoming continuation.

Dealing with second-order conditions, it is a common aim to derive so-called ``no-gap'' conditions,
meaning that the only difference between necessary and sufficient optimality conditions appears in the involved relation sign.
The above considerations suggest, however, that this might be only possible in a restrictive setting, regardless of what tool is used to handle the curvature of the feasible set,
and that it may not be of main interest.
Perhaps, it is better to focus on the purpose of necessary and sufficient conditions separately,
potentially even using different variational tools.
This seems to be in line with the comments of Poliquin and Rockafellar from their paper \cite{PoliquinRockafellar1998} introducing the concept of tilt stability:
``{\em The role of optimality conditions is seen
rather in the justification of numerical algorithms, in particular their stopping criteria,
convergence properties and robustness. From that angle, the goal of theory could be
different. Instead of focusing on the threshold between necessity and sufficiency, one
might more profitably try to characterize the stronger manifestations of optimality
that support computational work.}''
Let us also mention two related recent works of Rockafellar which use a conceptually very interesting approach
and deal with sufficient conditions with focus on stability properties of minimizers
as well as related aspects of numerical optimization, 
see \cite{Rockafellar2022a,Rockafellar2022b}.

In this paper, we focus solely on sufficient optimality conditions via second subderivatives, 
and we proceed by exploring basic calculus rules
that enable us to estimate the curvature of the feasible set by computable expressions.
Let us point out that for sufficient conditions, only lower estimates are relevant, and so we focus on those ones without trying to get the best possible upper estimates 
 which would be needed to infer necessary conditions.
The most important calculus principle is the pre-image rule which provides an estimate for the curvature of the feasible set $F^{-1}(C)$
in terms of second derivatives of $F$ and the curvature term associated with $C$.
For NLPs and disjunctive programs, $C$ has no curvature and so the sufficient conditions can be derived just from the pre-image rule.
For more challenging programs, one may need to further apply calculus rules to estimate the curvature of $C$ until eventually ending up with a set possessing no curvature,
e.g.\ a polyhedral set, and fully explicit optimality conditions.
This is precisely the case for the aforementioned setting from \cite{BeGfrYeZhouZhang1} which covers bilevel problems and problems with (quasi-) variational inequality constraints.
Therein, the set $C$ is the image of the pre-image of the graph of a normal cone mapping associated with a convex polyhedral set, and the latter is a polyhedral set as well.

The main message of the paper is that calculus for second subderivatives is very easy and so are the resulting sufficient second-order optimality conditions,
which can be readily applied to a large variety of optimization problems, including very challenging ones.
More precesily, the lower estimates relevant for sufficient conditions are valid under very mild assumptions.
Indeed, the pre-image rule yields a general lower estimate in the absence of constraint qualifications.
This is in fact a very simple observation, yet it seems that, with the exception of \cite[Section~3]{BeGfrYeZhouZhang1}, both standard works \cite{Rockafellar1989,RoWe98} 
as well as the more recent contributions 
\cite[Section~7]{MohammadiMordukhovichSarabi2021} and \cite{ThinhChuongAnh2021} 
postulated a superfluous constraint qualification for that purpose.
Let us mention that the sufficient conditions from \cite{MohammadiMordukhovichSarabi2021,ThinhChuongAnh2021}
can be derived by just applying the pre-image rule to estimate the second subderivative of the indicator function associated with the feasible set $F^{-1}(C)$,
while \cite[Theorem~3.3]{BeGfrYeZhouZhang1} provides a stronger result:
it needs no constraint qualification and a milder second-order condition, yet it yields a more stable minimizer.
Interestingly, the calculus rules are so versatile, that we are able to fully recover this result by estimating the second subderivative of the function
$f(x):=\max\{f_0(x)-f_0(\bar{x}), \dist(F(x),C)\}$, $x\in\R^n$, 
where $f_0$ is the objective function, $\xb$ is a fixed reference point, 
and $\dist(\cdot,C)$ stands for the distance function associated with $C$.
On the other hand, the image rule yields only an upper estimate for free, while the lower estimate is valid in the presence of a suitable qualification
condition. Here, we exploit the recently introduced inner calmness* property from
\cite{Be19} for this purpose. The latter is not very restrictive and can be efficiently verified.
Particularly, an easy consequence of the famous Walkup--Wets result on Lipschitzness of convex polyhedral mappings, see \cite{WaWe69},
is that polyhedral mappings enjoy inner calmness*, see \cite[Theorem~3.4]{Be19}.
General sufficient conditions for inner calmness* can be found in \cite{BenkoMehlitz2022a}.
Moreover, in applications to the challenging problem classes where the image rule is needed,
the inner calmness* assumption has to hold for an associated so-called multiplier
mapping which was shown to be inner calm* under reasonable conditions in \cite[Theorem~3.9]{Be19}.
Finally, in \cref{Prop:Ic*Proj} we prove that the projection mapping associated with any closed set is inner calm* at every point of the set.

The remainder of this paper is organized as follows.
In \cref{sec:preliminaries}, we provide the necessary preliminaries
from variational analysis and generalized differentiation which are
used in this paper. Particularly, we recall some essential theoretical
foundations of second subderivatives in \cref{sec:second_subderivative}
and, specifically, their role in the context of second-order optimality
conditions in \cref{sec:second_order_conditions}.
In \cref{sec:proximal_stuff}, we first address the calculus of second subderivatives 
associated with indicator functions in \cref{sec:indicator_function}, and as already
pointed out in \cite{BeGfrYeZhouZhang1}, this naturally leads to the introduction
of a so-called directional proximal normal cone. For later use, we also introduce a
directional proximal subdifferential in order to capture the
finiteness of second subderivatives in \cref{sec:directional_proximal_subdifferential}.
The essential \Cref{sec:calculus} is dedicated to the derivation of calculus rules for
second subderivatives. We derive a quite general composition rule 
as well as a marginal function rule which can be used to infer a pre-image and an
image rule for indicator functions, respectively.
Based on these findings, we are in position to easily
derive second-order sufficient optimality conditions in constrained optimization
over geometric constraints in \cref{sec:SOSC_geometric} 
and composite optimization \cref{sec:SOSC_composite}. Particular applications of these
results in disjunctive optimization, second-order cone programming, and
optimization with structured geometric constraints are presented in
\cref{sec:disjunctive_programs,sec:SOSC_second_order_cone,sec:structured_geometric_constraints},
respectively. The paper closes with some concluding remarks in \cref{sec:conclusions}.

\section{Preliminaries}\label{sec:preliminaries}

\subsection{Notation}

In this paper, we mainly use standard notation as utilized in \cite{RoWe98}.
For brevity of notation, we do not properly distinguish between a sequence and its terms, writing simply $z_k$ instead of, say, $\{z_k\}$
or $\{z_k\}_{k\in\N}$.

Throughout the paper, we equip $\R^n$ with the Euclidean norm $\norm{\cdot}$ and the Euclidean inner product $\skalp{\cdot,\cdot}$.
The open and closed $\varepsilon$-ball around some point $\bar z\in\R^n$ are denoted by $\mathbb U_\varepsilon(\bar z)$ and
$\mathbb B_\varepsilon(\bar z)$, respectively.
For given $w\in\R^n$, $\{w\}^\perp:=\{z^*\in\R^n\,|\,\skalp{z^*,w}=0\}$ is the annihilator of $w$.
By $\dist(\bar z,\Omega):=\inf_{z\in\Omega}\norm{z-\bar z}$ we denote the distance of $\bar z\in\R^n$ to a nonempty set $\Omega\subset\R^n$.
Additionally, for each nonempty index set $I\subset\{1,\ldots,n\}$, the vector $\bar z_I$
results from $\bar z$ by deleting all components whose indices do not belong to $I$.
We use $\mathtt e_1,\ldots,\mathtt e_n\in\R^n$ to denote the canonical unit vectors of $\R^n$.
For a twice continuously differentiable function $f_0\colon\R^n\to\R$ and some point $\bar z\in\R^n$, $\nabla f_0(\bar z)$ and
$\nabla^2f_0(\bar z)$ denote the gradient and the Hessian of $f_0$ at $\bar z$. For each $w\in\R^n$, we exploit
$\nabla^2f_0(\bar z)(w,w):=w^\top\nabla^2f_0(\bar z)w$.
For a mapping $F\colon\R^n\to\R^m$ and a vector $y^*\in\R^m$, $\skalp{y^*,F}\colon\R^n\to\R$ given by
$\skalp{y^*,F}(z):=\skalp{y^*,F(z)}$, $z\in\R^n$, is the associated scalarization mapping.
In case where $F$ is twice continuously differentiable and $\bar z$ is fixed, $\nabla F(\bar z)$ is the Jacobian of $F$ at $\bar z$.
Furthermore, we use
\[
	\skalp{y^*,\nabla^2F(\bar z)(w,w)}
	:=
	\sum_{i=1}^my_i^*\nabla^2F_i(\bar z)(w,w)
\]
for each $w\in\R^n$ for brevity of notation.
A single-valued mapping $G\colon\R^n\to\R^m$ is referred to as calm at $\bar z\in\R^n$
in direction $w\in\R^n$ if there is a constant $\kappa>0$ such that, for each sequences 
$t_k\downarrow 0$ and $w_k\to w$,
$\norm{G(\bar z+t_kw_k)-G(\bar z)}\leq\kappa t_k\norm{w_k}$
holds for sufficiently large $k\in\N$.
For $w:=0$, this notion recovers the classical property of $G$ to be calm at $\bar z$,
and this is weaker than local Lipschitz continuity of $G$ at $\bar z$.

Fix a closed set $\Omega\subset\R^n$ and some point $\bar z\in\Omega$. 
The sets
\begin{align*}
	\widehat N_\Omega(\bar z)
	&:=
	\{z^*\in\R^n\,|\,\skalp{z^*,z-\bar z}\leq\oo(\norm{z-\bar z})\;\forall z\in\Omega\},\\
	\widehat N_\Omega^p(\bar z)
	&:=
	\{z^*\in\R^n\,|\,\skalp{z^*,z-\bar z}\leq\OO(\norm{z-\bar z}^2)\,\forall z\in\Omega\}
\end{align*}
are referred to as the regular and proximal normal cone to $\Omega$ at $\bar z$, respectively, and these
are closed, convex cones by definition. In case $\tilde z\notin\Omega$, we set
$\widehat N_\Omega(\tilde z):=\widehat N_\Omega^p(\tilde z):=\varnothing$.
Based on the regular normal cone, we can define the so-called limiting normal cone to $\Omega$
at $\bar z$ by means of
\[
	N_\Omega(\bar z)
	:=
	\{z^*\in\R^n\,|\,\exists z_k\to\bar z,\,\exists z_k^*\to z^*,\,\forall k\in\N\colon\,z_k^*\in\widehat N_\Omega(z_k)\}.
\]
The latter is a closed cone which does not need to be convex. 
Again, we set $N_\Omega(\tilde z):=\varnothing$ for each $\tilde z\notin\Omega$.
Clearly, we have $\widehat N^p_\Omega(\bar z)\subset\widehat N_\Omega(\bar z)\subset N_\Omega(\bar z)$, and all these
cones coincide with the standard normal cone of convex analysis whenever $\Omega$ is convex.
For some direction $w\in\R^n$, let us recall that
\[
	N_\Omega(\bar z;w)
	:=
	\{z^*\in\R^n\,|\,\exists t_k\downarrow 0,\,\exists w_k\to w,\,\exists z_k^*\to z^*,\,\forall k\in\N\colon\,z_k^*\in\widehat N_\Omega(\bar z+t_kw_k)\}
\]
is referred to as the directional limiting normal cone to $\Omega$ at $\bar z$ in direction $w$.
	This notion has been introduced in \cite{Gfr13a,GinchevMordukhovich2011} in the Banach space setting,
	and a simplified finite-dimensional counterpart of the definition appears in \cite{Gfrerer2014}.
	The calculus of the directional normal cone has been explored in the paper \cite{BeGfrOut19}.
Whenever $w\in T_\Omega(\bar z)$ holds, where
\[
	T_\Omega(\bar z)
	:=
	\{w\in\R^n\,|\,\exists t_k\downarrow 0,\,\exists w_k\to w,\,\forall k\in\N\colon\,\bar z+t_kw_k\in\Omega\}
\]
denotes the tangent cone to $\Omega$ at $\bar z$, 
then $N_\Omega(\bar z;w)$ is a closed cone which satisfies $N_\Omega(\bar z;w)\subset N_\Omega(\bar z)$.
For arbitrary $w\notin T_\Omega(\bar z)$, $N_\Omega(\bar z;w)$ is empty.
We also set $N_\Omega(\tilde z;w):=\varnothing $ for each $\tilde z\notin\Omega$.

Let $\overline\R:=\R\cup\{-\infty,\infty\}$ denote the extended real line.
Recall that a function $h\colon\R^n\to\overline\R$ is called proper whenever $h(z)>-\infty$ hold for all $z\in\R^n$
while there is at least one $\bar z\in\R^n$ such that $h(\bar z)<\infty$.
Let us note that whenever $h$ is lower semicontinuous at $\bar z\in\R^n$ such that $|h(\bar z)|<\infty$,
then $h(z)>-\infty$ holds for all $z\in\R^n$ in a neighborhood of $\bar z$, i.e., $h$ is locally proper at $\bar z$.
The sets $\dom h:=\{z\in\R^n\,|\,h(z)<\infty\}$ and $\epi h:=\{(z,\alpha)\in\R^n\times\R\,|\,h(z)\leq\alpha\}$ are
called the domain and epigraph of $h$, respectively.
Observe that whenever $h$ is proper, we have $\dom h=\{z\in\R^n\,|\,|h(z)|<\infty\}$.
Furthermore, whenever $h$ is lower semicontinuous, then $\epi h$ is closed.

For a lower semicontinuous function $h\colon\R^n\to\overline\R$ and some point $\bar z\in\R^n$ where $|h(\bar z)|<\infty$,
the mapping $\dr h(\bar z)\colon\R^n\to\overline\R$ defined by
\[
	\dr h(\bar z)(w)
	:=
	\liminf\limits_{t\downarrow 0,\,w'\to w}\frac{h(\bar z+tw')-h(\bar z)}{t}
	\qquad\forall w\in\R^n
\]
is called the subderivative of $h$ at $\bar z$. 
Clearly, $\dr h(\bar z)$ is a lower semicontinuous mapping.
Observe that we have $\epi \dr h(\bar z)=T_{\epi h}(\bar z,h(\bar z))$.
Whenever $h$ is continuously differentiable at $\bar z$, we find $\dr h(\bar z)(w)=\skalp{\nabla h(\bar z),w}$ for all $w\in\R^n$.
Let us mention that $\dr h(\bar z)(w)$ is also referred to as lower Hadamard directional derivative of $h$ at $\bar z$ 
in direction $w$ in the literature.
	Note that the definition of the subderivative is also possible for arbitrary functions
	which are not necessarily lower semicontinuous.
	However, keeping in mind that our main purpose behind the consideration of variational objects
	is related to applications in mathematical optimization, it is reasonable to focus on lower semicontinuous
	functions. Some of the results in this paper can, via some nearby adjustments, be extended to functions 
	which are not lower semicontinuous. 
	Anyhow, for simplicity and brevity of presentation, we will use lower semicontinuity as a standing
	assumption in the remainder of the paper.

Recall that the sets
\begin{align*}
	\widehat\partial h(\bar z)
	&:=
	\{z^*\in\R^n\,|\,(z^*,-1)\in\widehat N_{\epi h}(\bar z,h(\bar z))\},
	\\
	\widehat\partial^p h(\bar z)
	&:=
	\{z^*\in\R^n\,|\,(z^*,-1)\in\widehat N_{\epi h}^p(\bar z,h(\bar z))\},
	\\
	\partial h(\bar z)
	&:=
	\{z^*\in\R^n\,|\,(z^*,-1)\in N_{\epi h}(\bar z,h(\bar z))\}
\end{align*}
are called the regular, proximal, and limiting subdifferential of $h$ at $\bar z$, respectively.
Clearly, we have 
$\widehat\partial^p h(\bar z)\subset\widehat\partial h(\bar z)\subset\partial h(\bar z)$, 
and whenever $h$ is convex, all these subdifferentials coincide
with the subdifferential in the sense of convex analysis.

\subsection{Inner semicompactness and inner calmness* of set-valued mappings}\label{sec:svm}

For a set-valued mapping $\Gamma\colon\R^n\tto\R^m$, the sets $\dom\Gamma:=\{x\in\R^n\,|\,\Gamma(x)\neq\varnothing\}$ and
$\gph\Gamma:=\{(x,y)\in\R^n\times\R^m\,|\,y\in\Gamma(x)\}$ are referred to as domain and graph of $\Gamma$, respectively.
Furthermore, $\Gamma^{-1}\colon\R^m\tto\R^n$ given by $\Gamma^{-1}(y):=\{x\in\R^n\,|\,y\in\Gamma(x)\}$, $y\in\R^m$, is the
inverse of $\Gamma$.

Let $(\bar x,\bar y)\in\gph\Gamma$ be fixed. 
We refer to $D\Gamma(\bar x,\bar y)\colon\R^n\tto\R^m$ given via $\gph D\Gamma(\bar x,\bar y):=T_{\gph\Gamma}(\bar x,\bar y)$
as the graphical derivative of $\Gamma$ at $(\bar x,\bar y)$.
If $\Gamma$ is single-valued at $\bar x$, we exploit the notation $D\Gamma(\bar x):=D\Gamma(\bar x,\Gamma(\bar x))$
for brevity. Note that whenever $G\colon\R^n\to\R^m$ is a single-valued and continuously differentiable mapping,
then we have $DG(\bar x)(u)=\{\nabla G(\bar x)u\}$ for each $u\in\R^n$.

Fix $\bar x\in\dom\Gamma$.
For a set $\Omega\subset\R^n$, we say that $\Gamma$ is inner semicompact at $\bar x$ w.r.t.\ $\Omega$
whenever for each sequence $x_k\to\bar x$ such that $x_k\in\Omega$ for all $k\in\N$, there exist
$y\in\R^m$ and a sequence $y_\ell\to y$ such that $y_\ell\in\Gamma(x_{k_\ell})$ holds for each $\ell\in\N$ where $x_{k_\ell}$ is a subsequence of $x_k$. 
If $\Omega:=\R^n$ can be chosen, $\Gamma$ is called inner semicompact at $\bar x$
for simplicity.

The next lemma shows that inner semicompactness at a single point already ensures inner
semicompactness in a neighborhood of this point.

\begin{lemma}\label{lem:inner_semicomactness_locally_preserved}
	Let $\Gamma\colon\R^n\tto\R^m$ be a set-valued mapping and let
	$\bar x\in\dom \Gamma$ be a point where $\Gamma$ is inner 
	semicompact w.r.t.\ $\dom\Gamma$.
	Then there exists $\varepsilon>0$ such that $\Gamma$ is
        inner semicompact w.r.t.\ $\dom\Gamma$ at all points
        from $\mathbb B_\varepsilon(\bar x)\cap\dom\Gamma$.
        Moreover, $\mathbb B_\varepsilon(\bar x)\cap\dom\Gamma$ is closed
        provided $\gph \Gamma$ is closed.
\end{lemma}
\begin{proof}
	Suppose that the assertion is not true.
	Then we find a sequence $x_k\to\bar x$ such that
	$x_k\in\dom\Gamma$ and $\Gamma$ is not inner semicompact
	at $x_k$ w.r.t.\ $\dom\Gamma$ for each $k\in\N$.
	Thus, for each $k\in\N$, there is a sequence $x_k^\ell\to x_k$
	such that $x_k^\ell\in\dom\Gamma$ for each $\ell\in\N$ and
	$\dist(0,\Gamma(x_k^\ell))\to\infty$ as $\ell\to\infty$.
	For each $k\in\N$, fix $\ell(k)\in\N$ such that
	$\norm{x_k^{\ell(k)}-x_k}\leq 1/k$ and $\dist(0,\Gamma(x_k^{\ell(k)}))\geq k$.
	Due to $x_k^{\ell(k)}\to\bar x$ and $x_k^{\ell(k)}\in\dom\Gamma$ for
	each $k\in\N$, $\Gamma$ is not inner semicompact at $\bar x$ w.r.t.\
	$\dom\Gamma$ which is a contradiction.
	
	The statement about closedness follows from \cite[Lemma~2.1]{Be19}
	by closedness of $\gph\Gamma$.	
\end{proof}

Let us also recall the inner calmness* property which has been coined in \cite{Be19} and further studied in \cite{BenkoMehlitz2022a}.
It can be interpreted as a quantitative version of the inner semicompactness from above.
For $\bar x\in\dom\Gamma$ and $\Omega\subset\R^n$, $\Gamma$ is called inner calm* at $\bar x$ 
w.r.t.\ $\Omega$ whenever there exists $\kappa>0$ such that for each sequence $x_k\to\bar x$
satisfying $x_k\in\Omega$ for all $k\in\N$, there exist $y\in\R^m$ and a sequence $y_\ell\to y$
with $y_\ell\in\Gamma(x_{k_\ell})$ and $\norm{y_\ell-y}\leq\kappa\norm{x_{k_\ell}-\bar x}$
for each $\ell\in\N$ where $x_{k_\ell}$ is a subsequence of $x_k$.
If $\Omega:=\R^n$ can be chosen, $\Gamma$ is called inner calm*
at $\bar x$ for brevity.
If, for given $u\in\R^n$, the above property holds just for all sequences $\bar x+t_ku_k=:x_k\in\Omega$
where $t_k\downarrow 0$ and $u_k\to u$, then $\Gamma$ is referred to as inner calm* at
$\bar x$ in direction $u$ w.r.t.\ $\Omega$.
Similarly as above, inner calmness* at $\bar x$ in direction $u$ is defined.

Inner calmness* is not very restrictive. Let us mention some situations when it is satisfied,
see also \cref{Prop:Ic*Proj} below.
\begin{itemize}
	\item[(i)] In \cite[Theorem~3.4]{Be19}, it was shown that polyhedral set-valued mappings are inner calm* at every point of the domain w.r.t.\ the domain.
	This result provides a certain lower/inner counterpart to Robinson's result on upper/outer Lipschitzness of polyhedral mappings from \cite{Rob81},
	and it is an easy consequence of the famous Walkup--Wets result on Lipschitzness w.r.t.\ the domain of convex polyhedral mappings, see \cite{WaWe69}.
	\item[(ii)] In \cite[Theorem~3.9]{Be19}, inner semicompactness and inner calmness* of a certain multiplier mapping associated with standard geometric constraints
	was established under suitable constraint qualifications.
	Let us mention that this multiplier mapping is very relevant for the analysis of the normal cone mapping associated with this constraint system.
	Moreover, (essentially) the same mapping also appears in applications of our second-order conditions to the most challenging optimization problems,
	see \cref{rem:multiplier_mappings}.
	\item[(iii)] Finally, due to \cite[Lemma~4.3]{BenkoMehlitz2022a}, whenever $\Gamma$ is isolatedly calm at all points
	$(\bar x,y)$ such that $y\in\Gamma(\bar x)$ and inner semicompact at $\bar x$ w.r.t.\ $\dom\Gamma$,
	then $\Gamma$ is inner calm* at $\bar x$ w.r.t.\ $\dom\Gamma$.
	Let us recall that $\Gamma$ is called isolatedly calm at some point $(\bar x,\bar y)\in\gph\Gamma$
if there exist $\kappa>0$, $\varepsilon>0$, and $\delta>0$ such that
\[
	\Gamma(x)\cap\mathbb U_\varepsilon(\bar y)
	\subset
	\bar y+\kappa\norm{x-\bar x}\mathbb B_1(0)
	\qquad
	\forall x\in\mathbb U_\delta(\bar x).
\]
It is well known that $\Gamma$ is isolatedly calm at $(\bar x,\bar y)$ if and only if
$D\Gamma(\bar x,\bar y)(0)=\{0\}$ holds, and the latter has been named Levy--Rockafellar criterion
in the literature.
\end{itemize}

The main purpose of introducing inner calmness* was its role in certain calculus rules, particularly, the so-called image rule for tangents, see \cite[Section~4]{Be19}.
In this paper, we further pursue these developments as inner calmness* is the essential assumption we rely on to derive calculus rules for second subderivatives
(it is essential for one (of two) pattern of related results while the other one requires only the most basic assumptions).

More precisely, in \cref{sec:calculus}, we will study a marginal function rule for the second subderivative, and we now state some preparatory results.
Therefore, we choose a proper, lower semicontinuous function $\varphi\colon\R^n\times\R^m\to\overline\R$
and consider the associated marginal function $h\colon\R^n\to\overline\R$ given by
\begin{equation}\label{eq:marginal_function}
	h(x)
	:=
	\inf\limits_{y\in\R^m}\varphi(x,y)\qquad \forall x\in\R^n.
\end{equation}
Closely related to $h$ are the set-valued mappings $\Upsilon\colon\R^n\times\R\tto\R^m$
and $\Psi\colon\R^n\tto\R^m$ given by
\begin{equation}\label{eq:M_and_S}
	\begin{aligned}
    	\Upsilon(x,\alpha)
    	&:=
    	\{y\in\R^m\,|\,(x,y,\alpha)\in\epi\varphi\}
    	\quad
    	&&\forall x\in\R^n,\,\alpha\in\R,
    	\\
    	\Psi(x)
    	&:=
    	\argmin_{y\in\R^m}\varphi(x,y)
    	\quad
    	&&\forall x\in\R^n.
    	\end{aligned}
\end{equation}
We further specify that $\Psi(x):=\varnothing$ is used for each $x\in\R^n$ such that
$\varphi(x,y)=\infty$ holds for all $y\in\R^m$.
Let us note that $\Psi$ is often called solution mapping in the literature and a tool of major
interest in parametric optimization. To the best of our knowledge, $\Upsilon$ is rarely used in
this regard, see e.g.\ \cite[Section~5.2]{BenkoMehlitz2022a}. 
However, we observe that $\Psi(x)=\Upsilon(x,h(x))$ is valid for all $x\in\R^n$ with $|h(x)|<\infty$,
i.e., $\Upsilon$ and $\Psi$ are related via the value function $h$.
Additionally, one can easily check that $\dom \Upsilon\subset\epi h$ is valid.
Finally, by lower semicontinuity of $\varphi$, $\gph \Upsilon$ is naturally closed.

Note that, given a closed set $\Omega \subset \R^n$, the distance function fits into
\eqref{eq:marginal_function} by observing
\begin{equation}\label{eq:distance_function_as_marginal_function}
	\dist(x,\Omega)=\inf_{y\in\R^n}\big(\norm{y - x} + \delta_\Omega(y)\big)
	\quad
	\forall x\in\R^n
\end{equation}
where $\delta_\Omega\colon\R^n\to\overline{\R}$
is the so-called indicator function of $\Omega$ which vanishes on $\Omega$ and is set to $\infty$ on $\R^n\setminus \Omega$.
Clearly, the argmin mapping $\Psi$ from \eqref{eq:M_and_S} 
corresponds to the projection mapping $P_\Omega\colon\R^n\tto\R^n$ given by
$P_\Omega(x) := \argmin_{y\in\Omega}\norm{y - x}$, $x\in\R^n$, in this case.
Interestingly, $P_\Omega$ is inner calm* at every point belonging to $\Omega$.
\begin{proposition}\label{Prop:Ic*Proj}
    Given a closed set $\Omega \subset \R^n$, the projection mapping $P_\Omega$ is inner calm*
    at every $\xb \in \Omega$.
    Moreover, for sequences $\xb \neq x_k \to \xb$ satisfying $\dist(x_k,\Omega)/\norm{x_k - \bar x} \to 0$,
    we even get for all $y_k \in P_\Omega(x_k)$,
    \begin{equation*}
      \frac{y_k - \xb}{\norm{x_k - \xb}} - \frac{x_k - \xb}{\norm{x_k - \xb}} \to 0,
    \end{equation*}
    i.e., along a subsequence, $y_k$ converges to $\xb$ from the same direction as $x_k$.
\end{proposition}
\begin{proof}
    Consider a sequence $x_k \to \bar x = P_\Omega(\bar x)$.
    By \cite[Example 1.20]{RoWe98}, $P_\Omega(x_k) \neq \varnothing$, and for every $y_k \in P_\Omega(x_k)$,
    we easily get
    \[\norm{y_k - \xb} \leq \norm{y_k - x_k} + \norm{x_k - \bar x} \leq 2 \norm{x_k - \bar x},\]
    showing inner calmness* of $P_\Omega$ at $\xb$. Moreover, the second claim follows from
    $(y_k - \xb) - (x_k - \xb) = y_k - x_k$, which means $\norm{(y_k - \xb) - (x_k - \xb)} = \dist(x_k,\Omega)$.
\end{proof}

To see that the statement of \cref{Prop:Ic*Proj} does not remain valid if $\xb \notin \Omega$,
one can consider $\xb$ to be the center of the unit sphere $\Omega\subset\R^2$. 
Then, the sequence $x_k$ can be chosen to approach $\xb$ rapidly along slowly changing rays. 
For instance $x_k := 1/k^2(\cos(1/k),\sin (1/k)) \to 0$, while
we have $y_k = (\cos(1/k),\sin (1/k)) \to (1,0)=:\yb$, $\norm{x_k - \xb} = 1/k^2$, 
and $\norm{y_k - \yb} \approx 1/k$.
Since this is not essential in the further parts of this paper, we omit the details.

We proceed with some basic results for the marginal function and the mappings from \eqref{eq:M_and_S}.

\begin{lemma}\label{lem:isc_of_M}
	Let $\varphi\colon\R^n\times\R^m\to\overline\R$ be a proper, lower semicontinuous function
	and consider the associated mappings defined in \eqref{eq:marginal_function} and
	\eqref{eq:M_and_S}.
	Fix $\bar x\in\dom \Psi$ and let $\Upsilon$ be inner semicompact at $(\bar x,h(\bar x))$ w.r.t.\ $\dom \Upsilon$.
	Then $h$ is lower semicontinuous at $\bar x$ and, locally around $(\bar x,h(\bar x))$, the
	sets $\epi h$ and $\dom \Upsilon$ coincide and are closed.
\end{lemma}
\begin{proof}
    From \cref{lem:inner_semicomactness_locally_preserved} we know that $\dom \Upsilon$ is closed locally around $(\bar x,h(\bar x))$ 
    and $\dom \Upsilon\subset\epi h$ holds always true.
    Thus, all the claims follow once we show that near $(\bar x,h(\bar x))$ also the opposite inclusion holds.
    This is, however, an easy consequence of the local closedness of $\dom \Upsilon$.
    Indeed, consider a closed neighborhood $U\subset\R^n\times\R$ of $(\bar x,h(\bar x))$
    such that $\dom \Upsilon \cap U$ is closed and let $(x,\alpha) \in \epi h \cap U$. Clearly, $(x,\alpha) \in \dom \Upsilon$ if $\alpha > h(x)$.
    If $\alpha = h(x)$, the definition of $h$ guarantees the existence of a sequence $y_k$ such that
    $\alpha = \lim_{k \to \infty} \varphi(x,y_k)$.
    Particularly, $(x,\varphi(x,y_k)) \in \dom \Upsilon$ and since $(x,\varphi(x,y_k)) \to (x,\alpha)$,
    the local closedness of $\dom \Upsilon$ yields $(x,\alpha) \in \dom \Upsilon$.  
\end{proof}

\begin{lemma}\label{lem:Psi_vs_M}
	Let $\varphi\colon\R^n\times\R^m\to\overline\R$ be a proper, lower semicontinuous function,
	consider the associated mappings defined in \eqref{eq:marginal_function} and
	\eqref{eq:M_and_S}, and fix $\bar x\in\dom \Psi$ as well as $u \in \R^n$.
	If $\Psi$ is inner semicompact at $\bar x$ (inner calm* at $\bar x$ in direction $u$),
	then $\Upsilon$ is inner semicompact at $(\bar x,\alpha)$ w.r.t.\ $\dom \Upsilon$ for all $\alpha \in \R$ 
	(inner calm* at $(\bar x,\alpha)$ in direction $(u,\mu)$ w.r.t.\ $\dom \Upsilon$ for all
	$\alpha, \mu \in \R$).
\end{lemma}
\begin{proof}
 	Assume that $\Psi$ is inner semicompact at $\bar x$ and consider
 	$\alpha \in \R$
 	as well as a sequence $(x_k,\alpha_k)\to(\bar x,\alpha)$ 
 	such that $(x_k,\alpha_k)\in\dom \Upsilon$ for each $k\in\N$.
	By passing to a subsequence (without relabeling),
	inner semicompactness of $\Psi$ yields the existence of $\bar y\in\R^m$ and, for each $k\in\N$, $y_k \in \Psi(x_k) = \Upsilon(x_k,h(x_k))$ 
	such that $y_k\to\bar y$.
	Since $(x_k,\alpha_k) \in \dom \Upsilon$, we have $\alpha_k \geq h(x_k)$ and, thus,  $y_k \in \Upsilon(x_k,\alpha_k)$ for each $k\in\N$,
	i.e., inner semicompactness of $\Upsilon$ follows.

	Assume now that $\Psi$ is inner calm* at $\bar x$ in direction $u$ and consider
 	$\alpha \in \R$
	as well as sequences $t_k \downarrow 0$ and $(u_k,\mu_k) \to (u,\mu)$
	with $(\bar x ,\alpha) + t_k (u_k,\mu_k) \in \dom \Upsilon$ for each $k\in\N$.
	By passing to a subsequence (without relabeling), the directional inner calmness* of $\Psi$ 
	yields the existence of $\bar y\in\R^m$, $\kappa>0$, and, for each $k\in\N$, $y_k \in \Psi(x_k) = \Upsilon(x_k,h(x_k))$
	such that $y_k\to\bar y$ and
	\[
		\norm{y_k - \bar y} \leq t_k \kappa \norm{u_k} \leq  t_k \kappa \norm{(u_k,\mu_k)}\qquad\forall k\in\N.
	\]
	For each $k\in\N$, $\alpha + t_k \mu_k \geq h(\bar x + t_k u_k)$ follows as before, and we obtain $y_k \in \Upsilon(\bar x + t_k u_k,\alpha + t_k \mu_k)$,
	showing the claimed inner calmness* of $\Upsilon$ at $(\bar x,h(\bar x))$
	in direction $(u,\mu)$.
\end{proof}

Due to the above lemmas, inner semicompactness of $\Upsilon$ (or $\Psi$) ensures lower semicontinuity of $h$
and so it
will be used as a standing assumption in our analysis.

It is worth noting that the assumptions of \cref{lem:Psi_vs_M}
imposed on $\Psi$ are more coarse, since they apply
whenever $x_k \to \bar x$ (from direction $u$)
regardless of what happens with $h(x_k)$.
They imply the corresponding assumptions on $\Upsilon$
for any $\alpha \in \R$.
Note also that the statements are trivially satisfied for $\alpha < h(\bar x)$ since there are no sequences 
$(x_k,\alpha_k) \in \dom \Upsilon \subset \epi h$ approaching such points $(\bar x,\alpha)$.
Similar arguments apply to the issue of directional convergence
in the image space.
Since
\[
	T_{\dom \Upsilon}(\bar x,h(\bar x)) 
	\subset 
	T_{\epi h}(\bar x,h(\bar x))
	= 
	\epi {\dr}h(\bar x),
\]
only for $\mu \geq {\dr}h(\bar x)(u)$, the inner calmness* of $\Upsilon$
is not trivially satisfied.
In \cref{sec:calculus}, we propose two analogous results,
the coarser one based on $\Psi$ and the finer one based on $\Upsilon$.

Observe that \cref{lem:Psi_vs_M} assumes inner semicompactness/inner calmness* of $\Psi$
which might be stronger than the restrictions of these properties to
$\dom\Psi$. Exemplary, the following example illustrates that inner semicompactness of 
$\Psi$ w.r.t.\ $\dom\Psi$ is not enough to guarantee lower semicontinuity of the 
marginal function $h$ at the point of interest, i.e., 
this assumption is too weak in order to imply the requirements of \cref{lem:isc_of_M}.
\begin{example}\label{ex:non_lower_semicontinuous_marginal_function}
	Consider $\varphi\colon\R\times\R\to\R$ given by $\varphi(x,y):=e^{xy}$ for all $x,y\in\R$.
	Then we find
	\[
		h(x)
		=
		\begin{cases}
			1	& x=0,\\
			0	& x\neq 0,
		\end{cases}
		\qquad
		\Psi(x)
		=
		\begin{cases}
			\R	&	x=0,\\
			\varnothing	&	x\neq 0
		\end{cases}
		\qquad
		\forall x\in\R.
	\]
	Particularly, $\Psi$ is inner semicompact at $\bar x:=0$ w.r.t.\ $\dom \Psi=\{\bar x\}$, but $h$ fails
	to be lower semicontinuous at $\bar x$.\\
	Note that the associated mapping $\Upsilon$ indeed fails to be inner semicompact at $(\bar x,1)$
	w.r.t.\ $\dom \Upsilon$.
	In order to see this, choose $x_k:=-1/k^2$ and $\alpha_k:=1-1/k$ for each $k\in\N$.
	Then we have $(x_k,\alpha_k)\to(\bar x,1)$ and $\Upsilon(x_k,\alpha_k)=[x_k^{-1}\ln\alpha_k,\infty)$
	for each $k\in\N$. Due to $x_k^{-1}\ln\alpha_k\to\infty$ as $k\to\infty$, $\Upsilon$ is not inner
	semicompact at $(\bar x,1)$.
\end{example}

Finally, we note that lower semicontinuity of $h$ is also implied
by a related (restricted) inf-compactness assumption, see e.g.\
\cite[Hypothesis~6.5.1]{Clarke1983} and \cite[Definition~3.8]{GuoLinYeZhang2014}.
Without stating precise definitions of these properties, let us mention that
inf-compactness is clearly implied by inner semicompactness of $\Psi$,
see also \cite[Section~4.1]{BaiYe2021}.
On the other hand, our refined semicompactness assumption imposed on
$\Upsilon$ turns out to be milder than (perhaps equivalent to) restricted inf-compactness.
Indeed, for $\bar x$ with $|h(\bar x)|<\infty$,
let a sequence $(x_k,\alpha_k) \to (\xb,h(\xb))$ satisfy $(x_k,\alpha_k)\in\dom\Upsilon$ for each $k\in\N$. 
Then, for each $k\in\N$ and by definition of $\Upsilon$, we get the existence of $y_k$ with
$h(x_k) \leq \varphi(x_k,y_k) \leq \alpha_k$. 
Combining this with $\alpha_k\to h(\bar x)$ and taking the restricted inf-compactness into consideration
yields the existence of a nonempty, compact set $A$ and, for large enough $k\in\N$, some
$\hat y_k \in \Psi(x_k) \cap A \subset \Upsilon(x_k,\alpha_k) \cap A$.
Thus, $\Upsilon$ is inner semicompact at $(\xb,h(\xb))$.

\subsection{Second subderivative}\label{sec:second_subderivative}

The following definition of the second subderivative of a function is taken from \cite[Definition~13.3]{RoWe98}.
\begin{definition}\label{def:second_subderivative}
	Let $h\colon \mathbb{R}^n \rightarrow \overline \R$ be a lower semicontinuous function
	and fix $\zb\in\R^n$ with $|h(\bar z)|<\infty$ as well as $z^*\in\R^n$.
	The second subderivative of $h$ at $\zb$ for $z^*$ is
	the function ${\dr}^2 h(\bar z;z^*)\colon\R^n\to \overline \R$ defined by
	\[
 		{\dr}^2h(\bar z;z^*)(w):= \liminf\limits_{t\downarrow 0,\,w'\to w} \frac{h(\zb+tw')-h(\zb)-t\langle z^*, w'\rangle}{\frac{1}{2}t^2}
 		\qquad 
 		\forall w \in  \mathbb{R}^n.
	\]
\end{definition}

Right from \cref{def:second_subderivative}, we readily obtain the homogeneity properties
\begin{equation}\label{eq:homogeneity_property}
	\dr ^2(\alpha h)(\bar z;z^*)(w)
	=
	\alpha \dr ^2 h(\bar z;z^*/\alpha)(w)
	\qquad
	\forall \alpha\in(0,\infty)
\end{equation}
as well as the relations
\begin{equation}\label{eq:trivial_infinity_of_second_subderivative}
 	\begin{aligned}
        {\dr}h(\bar z)(w) > \skalp{z^*,w}
        & \quad\Longrightarrow\quad
        {\dr}^2h(\bar z;z^*)(w) = \infty,\\
        {\dr}h(\bar z)(w) < \skalp{z^*,w}
        & \quad\Longrightarrow\quad
        {\dr}^2h(\bar z;z^*)(w) = -\infty.
     \end{aligned}
\end{equation}
In \cref{sec:directional_proximal_subdifferential}, we provide some more details regarding
when ${\dr}^2h(\bar z;z^*)(w)$ is finite.
Observe that whenever $h$ is twice differentiable and $z^* = \nabla h(\zb)$, 
then the limit inferior can be replaced by a limit 
in \cref{def:second_subderivative} since
\[
	h(\zb+tw')-h(\zb) = t \nabla h(\zb) w' + \tfrac 12 t^2\nabla^2 h(\zb)(w',w') + \oo(t^2).
\]

The following lemma will be important for some proofs.
\begin{lemma}\label{lem:special_sequences}
	Let $h\colon \mathbb{R}^n \rightarrow \overline \R$ be a lower semicontinuous function
	and fix $\zb\in\R^n$ with $|h(\bar z)|<\infty$, $z^*\in\R^n$, and $w\in\R^n$.
	Then there are sequences $t_k\downarrow 0$ and $w_k\to w$ such that
	\begin{subequations}\label{eq:subderivative_via_convergence}
		\begin{align}
			\label{eq:first_subderivative_via_convergence}
			\dr h(\bar z)(w)
			&=
			\lim\limits_{k\to\infty}\frac{h(\bar z+t_kw_k)-h(\bar z)}{t_k},\\
			\label{eq:second_subderivative_via_convernece}
			\dr^2h(\bar z;z^*)(w)
			&=
			\lim\limits_{k\to\infty}\frac{h(\bar z+t_kw_k)-h(\bar z)-t_k\skalp{z^*,w_k}}{\tfrac12t_k^2}.
		\end{align}
	\end{subequations}
\end{lemma}
\begin{proof}
    Fix sequences $t_k\downarrow 0$ and $w_k\to w$ which satisfy
	\eqref{eq:second_subderivative_via_convernece}.
    By passing to a subsequence (without relabeling) we may assume
    that
    $\mu_k :=(h(\bar z+t_kw_k)-h(\bar z))/t_k \to \mu \in \overline \R$
    with $\mu \geq \dr h(\bar z)(w)$.
    If $\mu = \dr h(\bar z)(w)$, the sequences $t_k$ and $w_k$ have the derised properties.
    If $\mu > \dr h(\bar z)(w)$, consider sequences $\tilde t_k\downarrow 0$ and $\tilde w_k\to w$ 
    such that $\tilde\mu_k:=(h(\bar z+\tilde t_k\tilde w_k)-h(\bar z))/\tilde t_k\to\dr h(\bar z)(w)$.
	Note that, for some $\varepsilon>0$ and all sufficiently large $k,\ell\in\N$, $\tilde\mu_k+\varepsilon< \mu_\ell$ holds.  
	We define $\ell(1):=\min\{\ell\in\N\,|\,t_\ell\leq \tilde t_1\}$
	and $\ell(k+1):=\min\{\ell\in\N\,|\,t_\ell\leq \tilde t_{k+1},\,\ell>\ell(k)\}$ for each $k\in\N$.
	Thus, we have $t_{\ell(k)}\leq \tilde t_k$ for all $k\in\N$.
	Furthermore, due to $\skalp{z^*,w_{\ell(k)}-\tilde w_k}\to 0$, 
	$\tilde \mu_k-\skalp{z^*,\tilde w_k}+\varepsilon/2\leq \mu_{\ell(k)}-\skalp{z^*,w_{\ell(k)}}$
	is valid for large enough $k\in\N$.
	Together, this gives
    \[
		\dr^2 h(\bar z;z^*)(w)
		=
		\lim\limits_{k\to\infty}\frac{\mu_{\ell(k)} -\skalp{z^*,w_{\ell(k)}}}{\tfrac12t_{\ell(k)}}
		\geq
		\lim\limits_{k\to\infty}\frac{\tilde\mu_k -\skalp{z^*,\tilde w_k}}{\tfrac12 \tilde t_k}
		\geq
		\dr^2 h(\bar z;z^*)(w),
	\]
	showing $(h(\bar z+\tilde t_k\tilde w_k)-h(\bar z)-\tilde t_k\skalp{z^*,\tilde w_k})/(\frac12\tilde t_k^2)\to\dr^2 h(\bar z;z^*)(w)$
	and completing the proof.
\end{proof}

Sequences $t_k\downarrow 0$ and $w_k\to w$ in the sense of \cref{lem:special_sequences} will be said to
recover $\dr h(\bar z)(w)$ and $\dr^2h(\bar z;z^*)(w)$ simultaneously.

In the proposition below, we summarize some elementary sum rules which address the second subderivative.
\begin{proposition}\label{prop:sum_rules}
\leavevmode
	\begin{enumerate}
		\item\label{item:sum_with_smooth_function} For a twice differentiable function $f_0\colon\R^n\to\R$, a lower semicontinuous function
			$f\colon\R^n\to\overline\R$, some point $\zb\in\R^n$ with $|f(\bar z)|<\infty$, 
			and $z^*\in\R^n$, we have
			\[
				\dr^2(f_0+f)(\bar z;z^*)(w)
				=
				\nabla^2 f_0(\zb)(w,w)+\dr^2 f(\bar z;z^*-\nabla f_0(\bar z))(w)\qquad\forall w\in\R^n.
			\]
		\item\label{item:separable_sum} For two lower semicontinuous functions $h_1\colon \R^{n_1} \to \overline\R$ and $h_2\colon \R^{n_2} \to \overline\R$,
			$h\colon \R^n \to \overline\R$ for $n:=n_1 + n_2$ defined by $h(z) := h_1(z_1) + h_2(z_2)$ for all $z:=(z_1,z_2)\in\R^n$,
			some point $\zb:=(\bar z_1,\bar z_2)\in\R^n$ with $|h(\bar z)|<\infty$, 
			and $z^*:=(z_1^*,z_2^*)\in\R^n$, we have
			\[
				\dr^2 h(\bar z;z^*)(w)
				\geq
				\dr^2 h_1(\bar z_1;z^*_1)(w_1)+\dr^2 h_2(\bar z_2;z^*_2)(w_2)\qquad\forall w:=(w_1,w_2)\in\R^n
			\]
			provided the right-hand side is not $\infty - \infty$.
	\end{enumerate}
\end{proposition}
\begin{proof}
	\begin{enumerate}
		\item 
			Exploiting a second-order Taylor expansion of $f_0$ at $\bar z$, we find
			\begin{align*}
				\dr^2(f_0&+f)(\zb;z^*)(w)
				\\
				&=
				\liminf\limits_{t\downarrow 0,\,w'\to w}
				\left(
				\frac{\nabla f_0(\zb)w'+\tfrac 12t\nabla^2f_0(\zb)(w',w')+\oo(t)}{\tfrac12t}
				+
				\frac{f(\zb+tw')-f(\zb)-t\langle z^*,w'\rangle}{\tfrac12t^2}
				\right)
				\\
				&=
				\nabla^2f_0(\zb)(w,w)
				+
				\liminf\limits_{t\downarrow 0,\,w'\to w}
				\frac{f(\zb+tw')-f(\zb)-t\langle z^*-\nabla f_0(\zb),w'\rangle}{\tfrac12t^2}
				\\
				&=
				\nabla^2f_0(\zb)(w,w)+\dr^2f(\bar z;z^*-\nabla f_0(\bar z))(w)
			\end{align*}
			which yields the claim.
		\item This follows immediately from the definition of the second subderivative,
			taking into account that $\liminf_{k\to\infty} (a_k + b_k) \geq \liminf_{k\to\infty} a_k + \liminf_{k\to\infty} b_k$ holds
			except for the indeterminate case $\infty - \infty$.
	\end{enumerate}
\end{proof}

Applying \cref{prop:sum_rules}\,\ref{item:sum_with_smooth_function} with $f$ being constantly zero, we find
\begin{equation}\label{eq:second_subderivative_smooth_map}
	\dr^2 f_0(\zb;z^*)(w)
	=
	\begin{cases}
		\nabla^2f_0(\bar z)(w,w)	&	z^*=\nabla f_0(\bar z),\\
		\infty						&   \langle z^*-\nabla f_0(\bar z),w\rangle<0,\\
		-\infty						&	\langle z^*-\nabla f_0(\bar z),w\rangle\geq 0,\,z^*\neq\nabla f_0(\bar z)
	\end{cases}
	\quad
	\forall w\in\R^n
\end{equation}
for each twice differentiable function $f_0\colon\R^n\to\R$, $\zb\in\R^n$, and $z^*\in\R^n$,
and this recovers \cite[Exercise~13.8]{RoWe98}.

We end this section by providing formulas for the second subderivative of two convex functions which we will need later,
namely the Euclidean norm and the maximum function.
In both cases, we only consider $z^*$ and $w$
such that $\skalp{z^*,w}$ equals the subderivative and
$z^*$ belongs to the subdifferential at the reference point, respectively, for otherwise
the second subderivative attains only the values $\pm \infty$,
see \eqref{eq:trivial_infinity_of_second_subderivative},
\cref{prop:basic_properties_second_subderivative,lem:geometric_vs_standard_proximal_subdifferential}, or \cite[Proposition 13.5]{RoWe98}.
Moreover, noting that the norm is twice differentiable at all non-zero points and keeping \eqref{eq:second_subderivative_smooth_map}
in mind, we restrict ourselves to the origin in the next lemma.
\begin{lemma}\label{lem:second_subderivative_of_euclidean_norm}
	For each $w\in\R^n$, we have $\dr\norm{\cdot}(0)(w) = \norm{w}$,
	and for each $z^*\in \partial \norm{\cdot}(0) = \mathbb B_1(0)$ satisfying $\skalp{z^*,w}=\norm{w}$, we find
	$\dr^2\norm{\cdot}(0;z^*)(w) = 0$.
\end{lemma}
\begin{proof}
	The formula for the subderivative is straightforward.
	Let us note that
	\[
		\dr^2\norm{\cdot}(0;z^*)(w)
		=
		\liminf\limits_{t\downarrow 0,\,w'\to w}
		\frac{2(\norm{w'}-\skalp{z^*,w'})}{t}.
	\]
	Choosing $w' := w$, we recover $0$, showing that $\dr^2\norm{\cdot}(0;z^*)(w) \leq 0$.
	Since $\norm{z^*}\leq 1$, the Cauchy--Schwarz inequality gives the estimate
	$\norm{w'}-\skalp{z^*,w'}\geq\norm{w'}(1-\norm{z^*})\geq 0$, i.e., $\dr^2\norm{\cdot}(0;z^*)(w) = 0$.
\end{proof}

\begin{lemma}\label{prop:vecmax}
	Let $\vecmax\colon\R^n\to\R$ be the function given by
	\[
		\vecmax(z):=\max\{z_1,\ldots,z_n\}\qquad\forall z\in\R^n.
	\]
	For fixed $\bar z\in\R^n$ and $w\in\R^n$, we define
	\[
		I(\bar z)
		:=
		\{i\in \{1,\ldots,n\}\,|\,\bar z_i=\vecmax(\bar z)\}.
	\]
	Then we have
	$\dr\vecmax(\bar z)(w)= \max\{w_i\,|\,i\in I(\bar z)\}$,
	and for each $z^*\in \partial \vecmax(\bar z)$ satisfying $\langle z^*,w\rangle=\max\{w_i\,|\,i\in I(\bar z)\}$,
	we find $\dr^2\vecmax(\bar z;z^*)(w) = 0$.
\end{lemma}
\begin{proof}
	The formula for the first subderivative is well known and can be distilled 
	from \cite[Exercise~8.31]{RoWe98}.
	Thus, let us pick $\bar z,w\in\R^n$ and
	\[
		z^*
		\in
		\partial \vecmax(\bar z)
		=
		\left\{z^* \in \R_{+}^n \,\middle|\, \mathsmaller\sum\nolimits_{i = 1}^n z_i^*=1, \, z^*_{i}=0 \, \forall i \notin I(\bar z) \right\}
	\]
	satisfying $\langle z^*,w\rangle=\max\{w_i\,|\,i\in I(\bar z)\}$.
	For an arbitrary index $i_0\in I(\bar z)$, we find
	\begin{align*}
		\dr^2\vecmax(\bar z;z^*)(w)
		&=
		\liminf\limits_{w'\to w,\,t\downarrow 0}
			\frac{\vecmax(\bar z+tw')-\bar z_{i_0}-t\langle z^*,w'\rangle}{\tfrac12t^2}
		\\
		&=
		\liminf\limits_{w'\to w,\,t\downarrow 0}
			\frac{\max\{tw'_i\,|\,i\in I(\bar z)\}-t\langle z^*,w'\rangle}{\tfrac12t^2}
		%\\
		%&=
		=\liminf\limits_{w'\to w,\,t\downarrow 0}
			\frac{\max\{w'_i\,|\,i\in I(\bar z)\}-\langle z^*,w'\rangle}{\tfrac12t}.
	\end{align*}
	Choosing $w' := w$, we recover $0$, showing that $\dr^2\vecmax(\bar z;z^*)(w) \leq 0$.
	On the other hand, we have
	\begin{align*}
		\max\{w'_i\,|\,i\in I(\bar z)\} - \langle z^*,w'\rangle
		&=
		\max\{w'_i\,|\,i\in I(\bar z)\} - \mathsmaller\sum\nolimits_{i\in I(\bar z)}z^*_i w'_i
		\\
		&\geq
		\max\{w'_i\,|\,i\in I(\bar z)\} \left(1 - \mathsmaller\sum\nolimits_{i\in I(\bar z)}z^*_i\right)
		=
		0,
	\end{align*}
	and $\dr^2\vecmax(\bar z;z^*)(w) = 0$ follows.
\end{proof}

\subsection{Second-order optimality conditions}\label{sec:second_order_conditions}

The following result, which can be found in \cite[Proposition~3.100]{BonSh00}, provides second-order necessary and sufficient
optimality conditions for the unconstrained minimization of proper functions.
\begin{proposition}\label{prop:SOC_unconstrained}
	Given a proper, lower semicontinuous function $h\colon\R^n \to \overline \R$ 
	and $\zb\in\R^n$ with $|h(\bar z)|<\infty$, the following statements hold.
	\begin{enumerate}
 		\item If $\zb$ is a local minimizer of $h$, then $0 \in \widehat\partial h(\zb)$ and 
 			$\dr^2h(\zb;0)(w) \geq 0$ for all $w\in\R^n$.
 		\item\label{item:SOC_unconstrained_sufficient}
 			Having $\dr^2h(\zb;0)(w) > 0$ for all $w\in\R^n\setminus\{0\}$ is equivalent to
			having the existence of $\varepsilon > 0$ and $\delta > 0$ such that
			\begin{equation}\label{eq:quadratic_growth_unconstrained}
				h(z) \geq h(\zb) + \varepsilon \norm{z - \zb}^2 \qquad\forall z\in\mathbb B_\delta(\zb).
			\end{equation}
			Particularly, $\bar z$ is a strict local minimizer of $h$.
	\end{enumerate}
\end{proposition}

Note that the assumptions of \cref{prop:SOC_unconstrained}\,\ref{item:SOC_unconstrained_sufficient}
also imply that $\dr h(\zb)(w) \geq \langle 0,w \rangle = 0$ holds for each $w\in\R^n\setminus\{0\}$,
see \cite[Proposition~13.5]{RoWe98}.
From \eqref{eq:trivial_infinity_of_second_subderivative}, $\dr h(\zb)(w) > 0$ for some $w\in\R^n$ yields $\dr^2h(\zb;0)(w) = \infty$.
In \cref{sec:directional_proximal_subdifferential}, we obtain a stronger statement, namely
that $\dr^2h(\zb;0)(w) > 0$ implies $0 \in \widetilde \partial^p h(\bar z;w)$ (which in turn yields $\dr h(\zb)(w) \geq 0$),
where $\widetilde \partial^p h(\bar z;w)$ denotes the directional proximal pre-subdifferential of $h$ at $\bar z$ in direction $w$,
see \cref{def:directional_proximal_subdifferentials} below.

Let us mention that whenever \eqref{eq:quadratic_growth_unconstrained} holds for some $\varepsilon>0$ and $\delta>0$,
then we say that $h$ satisfies a second-order growth condition at $\bar z$.

Consider now the problem of minimizing a twice differentiable function $f_0\colon\R^n \to \R$ 
over a closed set $S \subset \R^n$ and set $h := f_0 + \delta_S$ for the indicator function $\delta_S$ of $S$.
By closedness of $S$, $\delta_S$ is lower semicontinuous, and, obviously, $\delta_S$ is proper.
Taking into account \cref{prop:sum_rules}, we find
\[	
	\dr^2 h(\bar z;z^*)(w) = \nabla^2 f_0(\bar z)(w,w) + \dr^2 \delta_S(\bar z;z^* - \nabla f_0(\bar z))(w)
\]
for all $\bar z\in S$, $z^*\in\R^n$, and $w\in\R^n$.
\cref{prop:SOC_unconstrained} thus yields the following result.

\begin{proposition}\label{prop:SOC_constrained}
	Given a twice differentiable function $f_0\colon\R^n \to\R$ 
	and a closed set $S\subset\R^n$,
	the following statements hold.
		\begin{enumerate}
 			\item If $\zb \in \R^n$ is a local minimizer of $f_0$ over $S$, then $0 \in \nabla f_0(\zb) + \widehat N_S(\zb)$ and
				\[
					\nabla^2 f_0(\zb)(w,w) + \dr^2 \delta_S(\zb;-\nabla f_0(\zb))(w) \geq 0 \quad \forall w\in\R^n.
				\]
 			\item\label{item:SOC_constrained_sufficient} 
 				Having 
				\[ 
					\nabla^2 f_0(\zb)(w,w) + \dr^2 \delta_S(\zb;-\nabla f_0(\zb))(w) > 0 \quad\forall w\in\R^n\setminus\{0\}
				\]
				is equivalent to having the existence of $\varepsilon > 0$ and $\delta > 0$ such that
				\begin{equation}\label{eq:quadratic_growth}
					f_0(z) \geq f_0(\zb) + \varepsilon \norm{z - \zb}^2 \quad \forall z\in S\cap\mathbb B_{\delta}(\zb).
				\end{equation}
				Particularly, $\bar z$ is a strict local minimizer of $f_0$ over $S$.
		\end{enumerate}
\end{proposition}

Note again that $0 \in \widehat\partial(f_0+\delta_S)(\bar z)=\nabla f_0(\zb) + \widehat N_S(\zb)$ 
is implied by the requirements of 
\cref{prop:SOC_constrained}\,\ref{item:SOC_constrained_sufficient} 
and \cite[Exercise~8.8]{RoWe98}.
Moreover, \cref{prop:basic_properties_second_subderivative_indicator} from below, equalling \cite[Proposition~2.18]{BeGfrYeZhouZhang1}, in this case
yields that, for each non-zero $w\in\R^n$, either $w \notin T_S(\zb)$ or $0 \in \nabla f_0(\zb) + \Np_S(\bar z;w)$, which also gives $\skalp{\nabla f_0(\zb), w} \geq 0$.
Here, $\Np_S(\bar z;w)$ stands for the proximal pre-normal cone to $S$ at $\bar z$ in direction $w$, 
see \cref{def:directional_proximal_normal_cone} below.
Following the proof of \cref{prop:SOC_unconstrained}, we even find $\dr^2\delta_S(\bar z;-\nabla f_0(\bar z))(0)=0$ which
gives $0\in\nabla f_0(\bar z)+\widehat N^p_S(\bar z)$, see \cref{prop:basic_properties_second_subderivative_indicator} again.

Again, whenever there are $\varepsilon>0$ and $\delta>0$ such that \eqref{eq:quadratic_growth} holds, we say that the optimization problem
$\min\{f_0(z)\,|\,z\in S\}$ satisfies the second-order growth condition at $\bar z$.

\section{Directional proximal normal cones and subdifferentials}\label{sec:proximal_stuff}

\subsection{Second subderivative of the indicator function and the directional proximal normal cone}\label{sec:indicator_function}

Motivated by the second-order optimality conditions for constrained optimization problems from \cref{prop:SOC_constrained}, 
we look deeper into the second subderivative of the indicator function 
$\delta_S\colon\R^n\to\overline\R$ for a closed set $S \subset \R^n$.
As noted in \cite{BeGfrYeZhouZhang1}, given $\zb\in S$ and $z^* \in \R^n$, we get 
\begin{equation}\label{eq:indicatorf}
	{\dr}^2\delta_{S}(\zb;z^*)(w)
 	=
 	\liminf\limits_{t\downarrow 0,\,w'\to w}\frac{\delta_{S}(\bar z+tw')-\delta_{S}(\bar z)-t\langle z^*, w' \rangle }{\frac{1}{2}t^2}
	=
	\liminf\limits_{\substack{t\downarrow 0,\,w'\to w,\\ \zb+tw'\in S}}\frac{-2\langle z^*, w' \rangle }{t}
	\quad
	\forall w\in T_S(\bar z)
\end{equation}
directly from \cref{def:second_subderivative}.
It is immediate that $w\notin T_S(\zb)$ or $\langle z^*,w\rangle<0$ implies $\dr^2\delta_S(\zb;z^*)(w)=\infty$, 
see \cref{prop:basic_properties_second_subderivative_indicator} as well.
Additionally, for $\bar z\in\inn S$, we obviously have
\[
	\dr^2\delta_S(\bar z;z^*)(w)
	=
	\begin{cases}
		\infty		& \skalp{z^*,w}<0,\\
		0			& z^*=0,\\
		-\infty		& \text{otherwise}
	\end{cases}
	\qquad
	\forall z^*\in\R^n,\,\forall w\in\R^n.
\]
Moreover, since $\delta_S = \alpha \delta_S$ for $\alpha > 0$,
the homogeneity property \eqref{eq:homogeneity_property}
can be restated as
\begin{equation}\label{eq:homogeneity_property_2}
	\dr^2\delta_S(\bar z;\alpha z^*)(w)
	=
	\alpha \dr^2\delta_S(\bar z;z^*)(w)
	\qquad
	\forall \alpha\in(0,\infty).
\end{equation}

Let us investigate a simple example where the second subderivative of an indicator function at some
boundary point can be calculated easily.
\begin{example}\label{ex:second_subderivative_of_indicator_of_negative_real_line}
	A direct calculation gives
	\[
		\dr^2\delta_{\R_-}(0;z^*)(w)
		=
		\begin{cases}
			\infty	&	w>0\,\text{ or }\,z^*w<0,\\
			0		&	z^*\geq 0,\,w\leq 0,\,z^*w=0,\\
			-\infty	&	\text{otherwise}
		\end{cases}
		\qquad\forall z^*\in\R,\,\forall w\in\R.
	\]
\end{example}

The following simple calculus rules are consequences of \eqref{eq:indicatorf}.
\begin{lemma}\label{prop:second_subderivative_of_indicator_of_unions}
	Let $S_1,\ldots,S_\ell\subset\R^n$ be closed sets, set $S:=\bigcup_{i=1}^\ell S_i$, 
	and fix $\bar z\in S$ as well as $z^*\in\R^n$.
	Then we have
	\[
		\dr^2\delta_S(\bar z;z^*)(w)
		=
		\inf\limits_{i\in J(\bar z;w)}\dr^2\delta_{S_i}(\bar z;z^*)(w)\qquad\forall w\in\R^n
	\]
	where $J(\bar z;w):=\{i\in\{1,\ldots,\ell\}\,|\,\bar z\in S_i,
	w\in\ T_{S_i}(\bar z)\}$.
	If $S_1,\ldots,S_\ell$ are convex polyhedral sets, we have
	\[
		\dr^2\delta_S(\bar z;z^*)(w)
		=
		\inf\limits_{i\in J(\bar z;w)}\delta_{K_{S_i}(\bar z;z^*)}(w)\qquad\forall w\in\R^n
	\]
	for all $z^*\in\bigcap_{i\in J(\bar z;w)}N_{S_i}(\bar z;w)$,
	where $K_{S_i}(\bar z;z^*):=T_{S_i}(\bar z)\cap\{z^*\}^\perp$ is used for each $i\in J(\bar z;w)$.
\end{lemma}
\begin{proof}
	Let us define $J(\bar z):=\{i\in\{1,\ldots,\ell\}\,|\,\bar z\in S_i\}$.
	Then, from \eqref{eq:indicatorf}, we find
	\[
		\dr^2\delta_S(\bar z;z^*)(w)
		=
		\inf\limits_{i\in J(\bar z)}\dr^2\delta_{S_i}(\bar z;z^*)(w)\quad\forall w\in\R^n.
	\]
	Note that for each $i\in J(\bar z)\setminus J(\bar z;w)$, we
	have $w\notin T_{S_i}(\bar z)$ and, thus, $\dr^2\delta_{S_i}(\bar z;z^*)(w)=\infty$,
	i.e., the indices from $J(\bar z)\setminus J(\bar z;w)$ can be discarded
	from the infimum.

	In order to prove the formula for the union of convex polyhedral sets,
	we fix $w\in\R^n$ arbitrarily.
	Applying \cite[Exercise~13.17]{RoWe98} gives
	\[
		\dr^2\delta_{S_i}(\bar z;z^*)(w)
		=
		\delta_{K_{S_i}(\bar z;z^*)}(w)
	\]
	for all $z^*\in N_{S_i}(\bar z)$ and $i\in J(\bar z;w)$.
	Due to $N_{S_i}(\bar z;w)\subset N_{S_i}(\bar z)$
	for all $i\in J(\bar z;w)$, the stated formula is valid for all
	$z^*\in\bigcap_{i\in J(\bar z;w)}N_{S_i}(\bar z;w)$ which gives the claim.
\end{proof}

A related result can be found in \cite[Proposition~3.2]{ThinhChuongAnh2021}.
Let us point out that \cref{prop:second_subderivative_of_indicator_of_unions} 
together with \cref{prop:basic_properties_second_subderivative_indicator} from below shows 
that the second subderivative of an indicator associated with a union 
of finitely many convex polyhedral sets takes only value $0$ as soon as it is finite.

\begin{lemma}\label{lem:second_subderivative_of_indicator_of_product}
	For $\ell\in\N$ and $n_1,\ldots,n_\ell\in\N$, we fix closed sets $S_i\subset\R^{n_i}$ and
	points $\bar z_i\in S_i$, $i=1,\ldots,\ell$.
	Then, for arbitrary $z_i^*,w_i\in\R^{n_i}$, $i=1,\ldots,\ell$, we have
	\[
		\dr^2\delta_{S_1\times\ldots\times S_\ell}
		((\bar z_1,\ldots,\bar z_\ell);(z^*_1,\ldots,z^*_\ell))
		(w_1,\ldots,w_\ell)
		\geq 
		\sum_{i=1}^\ell\dr^2\delta_{S_i}(\bar z_i;z^*_i)(w_i)
	\]
	provided the right-hand side does not contain summands 
	of type $-\infty$ and $\infty$ simultaneously.
\end{lemma}

It turns out that there is a connection between the second subderivative of a set indicator
and a so-called directional proximal normal cone to the same set, which has been defined
in \cite[Definition~2.8]{BeGfrYeZhouZhang1}.

\begin{definition}\label{def:directional_proximal_normal_cone}
	Given a closed set $S\subset\mathbb R^n$, a point $\bar z \in S$, and a direction $w\in T_S(\bar z)$, 
	we define the proximal pre-normal cone to $S$ in direction $w$ at $\bar z$ as
	\begin{align*}
		\Np_S(\bar z;w)
		&:=
		\left\{
			z^*\in\R^n\,\middle|\,
				\begin{aligned}
					&\exists\gamma>0,\,\forall t_k\downarrow 0,\,\forall w_k\to w\text{ such that }\bar z+t_kw_k\in S\,\forall k\in\N\colon\\
					&\qquad\skalp{z^*,w_k}\leq\gamma t_k\norm{w_k}^2\text{ for all sufficiently large }k\in\N
				\end{aligned}
		\right\},
	\end{align*}
	and the proximal normal cone to $S$ at $\bar z$ in direction $w$ as
	\[
		\widehat N^p_S(\bar z;w)
		:= 
		\Np_S(\bar z;w)\cap \{w\}^\perp.
	\]
	In case where $w\not\in T_S(\bar z)$, we set $\Np_S(\bar z;w):= \widehat N^p_S(\bar z;w):=\varnothing$.
\end{definition}

In \cite{BeGfrYeZhouZhang1}, the proximal pre-normal cone has been defined in a different way without
sequences in order to visualize its close relationship to the standard proximal normal cone.
In the lemma below, we show that both definitions are equivalent.
\begin{lemma}\label{lem:proximal_pre-normal_cone}
	Given a closed set $S\subset\R^n$, a point $\bar z\in S$, and a direction $w\in T_S(\bar z)$, we have
	\[
		\Np_S(\bar z;w)
		=
		\{
			z^* \in \R^n \,|\, \exists\gamma,\delta,\rho>0\colon\, \skalp{z^*,z-\bar z}\leq \gamma \|z-\bar z\|^2 \ 
			\ \forall z\in S\cap (\bar z+\mathbb V_{\delta,\rho}(w))
		\}
	\]
	where, for each $\delta,\rho>0$, $\mathbb V_{\delta,\rho}(w)\subset\R^n$ is a so-called directional
	neighborhood of $w$ and given by means of
	\[
		\mathbb V_{\delta, \rho}(w) 
		:= 
		\left \{ 
			w'\in \mathbb{B}_\delta(0)\,\middle | \,
 			\big\| \| w \|w' -\|w'\| w \big\|\leq \rho \|w'\|\| w \|  
 		\right\}.
	\]
\end{lemma}
\begin{proof}
	For $w:=0$, the equivalence follows easily since the directional neighborhood coincides with a conventional
	neighborhood. In this case, both formulas recover the classical proximal normal cone to $S$ at $\bar z$.
		Indeed, the inclusion $\supset$ can be easily shown by a direct calculation.
		On the other hand, suppose that, for some $z^*\in\Np_S(\bar z;0)$, there is a sequence $z_k\to\bar z$
		such that $z_k\in S$ and $\skalp{z^*,z_k-\bar z}>k\norm{z_k-\bar z}^2$ for all $k\in\N$.
		Then $z_k\neq\bar z$ for all $k\in\N$, and we can set $t_k:=\norm{z_k-\bar z}^{1/2}$ and $w_k:=(z_k-\bar z)/t_k$
		for each $k\in\N$ in order to get $t_k\downarrow 0$, $w_k\to 0$, as well as $\bar z+t_kw_k\in S$ and
		$\skalp{z^*,w_k}>kt_k\norm{w_k}^2$ for each $k\in\N$ contradicting $z^*\in\Np_S(\bar z;0)$.
	Thus, let us assume that $w\neq 0$.
	 
	Fix $z^*\in\Np_S(\bar z;w)$. 
	Suppose that for each $k\in\N$, there is $z_k\in S\cap(\bar z+\mathbb V_{1/k,1/k}(w))$ such that
	$z_k\to\bar z$ and $\skalp{z^*,z_k-\bar z}>k\norm{z_k-\bar z}^2$.
	This is only possible if $z_k\neq \bar z$ for each $k\in\N$, so we can set $t_k:=\norm{z_k-\bar z}/\norm{w}$ and
	$w_k:=\norm{w}(z_k-\bar z)/\norm{z_k-\bar z}$ which gives $z_k=\bar z+t_kw_k$ for each $k\in\N$.
	Furthermore, we find $t_k\downarrow 0$ and $w_k\to w$ by definition of the directional neighborhood. 
	By construction, we have $t_k\skalp{z^*,w_k}>k\,t_k^2\norm{w_k}^2$ for each $k\in\N$ contradicting the
	definition of the proximal pre-normal cone.
	Thus, the inclusion $\subset$ has been shown. The proof of the converse inclusion is analogous.
\end{proof}

It is clear from \cref{lem:proximal_pre-normal_cone} that the proximal pre-normal cone is, in general,
larger than the classical proximal normal cone, i.e., $\widehat{N}_S^p(\bar z) \subset \Np_S(\bar z;w)$ 
holds for each $\bar z \in S$ and $w\in T_S(\bar z)$.
Moreover, for any such $\bar z$ and $w$ as well as each vector $z^*\in\R^n$ satisfying $\langle z^*, w\rangle <0$,
\cref{def:directional_proximal_normal_cone} yields $z^*\in\Np_S(\bar z;w)$.
In fact, we even have the estimates
\begin{equation}\label{eq:bounds_proximal_prenormal_cone}
	\{z^*\in\R^n\,|\, \langle z^*, w\rangle <0\} 
	\subset 
	\Np_S(\bar z;w) 
	\subset 
	\{z^*\in\R^n\,|\,\langle z^*, w\rangle \leq 0\}.
\end{equation}
Consequently, if a vector $z^*\in\Np_S(\bar z;w)$ satisfies $\langle z^*, w\rangle <0$, it does not provide much useful information.
Hence, it is natural to intersect the directional proximal pre-normal cone with the annihilator of $w$ in order to 
define the directional proximal normal cone.

As shown in \cite[Proposition~2.9]{BeGfrYeZhouZhang1},
both $\Np_S(\bar z;w)$ and $\widehat N^p_S(\bar z;w)$ are convex cones for $\bar z\in S$ and $w\in T_S(\bar z)$.
Furthermore, one has the estimates
\begin{equation}\label{EqPropDirProxNormalCone}
	\widehat N^p_S(\bar z)\cap\{w\}^\perp
	\subset 
	\widehat N^p_S(\bar z;w)
	\subset 
	\widehat N_{T_S(\bar z)}(w)
	\subset 
	N_{T_S(\bar z)}(w) 
	\subset 
	N_S(\bar z;w).
\end{equation}
In particular, when $S$ is a closed convex set, we get
\begin{equation}\label{EqDirProxNormalConeConvex}
	\widehat N^p_S(\bar z;w)
	=
	N_S(\bar z;w)
	= 
	N_S(\bar z)\cap \{w\}^\perp
	=
	N_{T_S(\bar z)}(w)
\end{equation}
for any such $\bar z$ and $w$.

In \cite[Proposition~2.18]{BeGfrYeZhouZhang1}, the following result was proven. 

\begin{proposition}\label{prop:basic_properties_second_subderivative_indicator}
 	Consider a closed set $S\subset\R^n$, $\bar z \in S$, $z^*\in\R^n$, and $w\in\R^n$.
 	Then the following statements hold.
	\begin{enumerate}
  		\item\label{item:equal_pos_infty} 
  			If $w\not\in T_S(\bar z)$ or $\langle z^*,w\rangle<0$, then $\dr ^2\delta_S(\bar z;z^*)(w)=\infty$.
  		\item\label{item:greater_neg_infty} 
  			For $w\in T_S(\bar z)$, we have $\dr^2\delta_S(\bar z;z^*)(w)>-\infty$ if and only if $z^*\in \Np_S(\bar z;w)$.
  		\item If $\dr^2\delta_S(\bar z;z^*)(w)$ is finite, then $z^*\in \widehat N_S^p(\bar z;w)$.
	\end{enumerate}
\end{proposition}

\subsection{The directional proximal subdifferential}\label{sec:directional_proximal_subdifferential}

In this subsection, we interrelate the second subderivative of a lower semicontinuous
function with a new directional proximal subdifferential which is introduced below.

\begin{definition}\label{def:directional_proximal_subdifferentials}
	Given a lower semicontinuous function $h\colon\R^n \to \overline \R$, a point $\bar z \in \R^n$ such that $|h(\bar z)|<\infty$,
	and a direction $w\in \R^n$ such that $|\dr h(\bar z)(w)| < \infty$,
	we define the proximal pre-subdifferential of $h$ at $\bar z$ in direction $w$ as
	\begin{align*}
		\widetilde\partial^p h(\bar z;w)
		:=
		\{z^*\in\R^n\,|\,\dr^2h(\bar z;z^*)(w)>-\infty\},
	\end{align*}
	and the proximal subdifferential of $h$ at $\bar z$ in direction $w$ as
	\[
		\widehat \partial^p h(\bar z;w)
		:= 
		\widetilde\partial^p h(\bar z ;w)\cap \{z^*\in\R^n\,|\,\dr h(\bar z)(w) = \langle z^*,w\rangle\}.
	\]
	If $|\dr h(\bar z)(w)| = \infty$, we set $\widetilde\partial^p h(\bar z ;w) := \widehat \partial^p h(\bar z;w) := \varnothing$.
	Finally, for some $\omega\in\R$, the sets
    \begin{align*}
        \widetilde \partial^p_g h(\bar z;(w,\omega))
        & :=
        \{z^*\in\R^n\,|\,(z^*,-1)\in \Np_{\epi h}((\bar z,h(\bar z));(w,\omega))\},\\
        \widehat \partial^p_g h(\bar z;(w,\omega))
        & :=
        \{z^*\in\R^n\,|\,(z^*,-1)\in\widehat N^p_{\epi h}((\bar z,h(\bar z));(w,\omega))\}
    \end{align*}
	are referred to as the geometric proximal pre-subdifferential and subdifferential of $h$ at $\bar z$ in direction $(w,\omega)$, respectively.
\end{definition}

Note that in case $\dr h(\bar z)(w) = -\infty$, we have $\skalp{z^*,w}>\dr h(\bar z)(w)$, so that \eqref{eq:trivial_infinity_of_second_subderivative}
gives $\dr^2h(\bar z;z^*)(w)=-\infty$, i.e., $\widetilde\partial^ph(\bar z;w)=\varnothing$ would also follow from the
defining relation of the proximal pre-subdifferential. 
On the other hand, in case $\dr h(\bar z)(w) = \infty$, we find $\dr^2h(\bar z;z^*)(w)=\infty$ for each $z^*\in\R^n$ 
from \eqref{eq:trivial_infinity_of_second_subderivative},
and the defining relation of the proximal pre-subdifferential would give us the whole space,
but we still define the set $\widetilde\partial^p h(\bar z;w)$ to be empty.
Note that this parallels the definition of the proximal pre-normal cone in a direction which is not tangent, see \cref{def:directional_proximal_normal_cone}.
Particularly, if $\dr h(\bar z)(w) = \infty$, we get $\widetilde \partial^p_g h(\bar z;(w,\omega)) = \widehat \partial^p_g h(\bar z;(w,\omega)) = \varnothing$
for all $\omega \in \R$ since $T_{\epi h}(\bar z,h(\bar z)) = \epi \dr h(\bar z)$, i.e., there is no $\omega \in \R$ 
with $(w,\omega) \in T_{\epi h}(\bar z,h(\bar z))$. 

Note that if $|\dr h(\bar z)(w)| < \infty$, then \eqref{eq:trivial_infinity_of_second_subderivative} yields
\begin{equation}\label{eq:trivial_estimates_proximal_pre-subdifferential}
	\{z^*\in\R^n\,|\, {\dr}h(\bar z)(w) > \skalp{z^*,w}\}
	\subset
	\widetilde\partial^p h(\bar z;w)
	\subset
	 \{z^*\in\R^n\,|\, {\dr}h(\bar z)(w) \geq \skalp{z^*,w}\},
\end{equation}
which motivates the definition of the proximal subdifferential of $h$ at $\bar z$ in direction $w$. 

Similar as in \cref{lem:proximal_pre-normal_cone}, we can show that the proximal pre-subdifferential
admits an alternative representation via directional neighborhoods which is free of sequences.
\begin{lemma}\label{lem:proximal_pre_subdifferential}
	Given a lower semicontinuous function $h\colon\R^n\to\overline\R$, a point $\bar z\in\R^n$ such that
	$|h(\bar z)|<\infty$, and a direction $w\in\R^n$ such that $|\dr h(\bar z)(w)|<\infty$, we have
	\[	
		\widetilde\partial ^ph(\bar z;w)
		=\{z^* \in \R^n \,|\, \exists \gamma,\delta,\rho>0\colon\,
		h(z) \geq h(\bar z) + \langle z^*,z-\bar z\rangle - \gamma \|z-\bar z\|^2\ \
		\forall z\in \bar z+\mathbb V_{\delta,\rho}(w)\}.
	\]
\end{lemma}
\begin{proof}
		For $w:=0$, both definitions 
		are equivalent by definition of the second subderivative
		since the directional neighborhood reduces to a classical neighborhood.
		In this case, both definitions recover the classical proximal subdifferential 
		of $h$ at $\bar z$.
		Indeed, the inclusion $\supset$ can be easily shown by a direct calculation.
		On the other hand, suppose that, for some $z^*\in\widetilde\partial ^ph(\bar z;0)$, 
		there is a sequence $z_k\to\bar z$
		such that $h(z_k) - h(\bar z) - \langle z^*,z_k-\bar z\rangle < - k \|z_k-\bar z\|^2$ 
		for all $k\in\N$.
		Then $z_k\neq\bar z$ for all $k\in\N$, and we can set either
		$t_k:=\norm{z_k-\bar z}^{1/2}$ and $w_k:=(z_k-\bar z)/t_k$ provided 
		$k \|z_k-\bar z\| \to \infty$ holds at least along a subsequence
		or
		$t_k:=k^{1/3}\norm{z_k-\bar z}$ and $w_k:=(z_k-\bar z)/t_k$ provided 
		the expression $k \|z_k-\bar z\|$ remains bounded.		
		In both cases, $t_k\downarrow 0$, $w_k\to 0$, as well as
		$(h(\bar z + t_k w_k) - h(\bar z) - t_k \langle z^*,w_k\rangle)/t_k^2 < - k \|z_k-\bar z\|^2/t_k^2$ for each $k\in\N$.
		Since $k \|z_k-\bar z\|^2/t_k^2 \to \infty$ in both cases (at least along a subsequence), 
		this contradicts $z^*\in\widetilde\partial ^ph(\bar z;0)$.
	Thus, let us assume that $w\neq 0$.
	
	Fix $z^*\in\widetilde\partial^ph(\bar z;w)$. Suppose that for each $k\in\N$, there is 
	$z_k\in \bar z+\mathbb V_{1/k,1/k}(w)$
	such that $z_k\to\bar z$ and $h(z_k)<h(\bar z)+\skalp{z^*,z_k-\bar z}-k\norm{z_k-\bar z}^2$.
	This gives $z_k\neq\bar z$ for each $k\in\N$, so we can set $t_k:=\norm{z_k-\bar z}/\norm{w}$ and
	$w_k:=\norm{w}(z_k-\bar z)/\norm{z_k-\bar z}$ which gives $z_k=\bar z+t_kw_k$ for each $k\in\N$.
	By construction, we find
	\[
		\frac{h(\bar z+t_kw_k)-h(\bar z)-t_k\skalp{z^*,w_k}}{t_k^2\norm{w_k}^2}<-k
	\]
	for each $k\in\N$. Multiplying this with $2\norm{w_k}^2$ while noting that $\norm{w_k}\to\norm{w}\neq 0$,
	taking the limit $k\to\infty$ gives $\dr^2h(\bar z;z^*)(w)=-\infty$ which contradicts the definition of
	the proximal pre-subdifferential. 
	Thus, the inclusion $\subset$ has been shown. The converse relation follows in analogous fashion.
\end{proof}

By construction, the directional proximal pre-normal and normal cone from \cref{def:directional_proximal_normal_cone}
to some closed set correspond to the directional proximal pre-subdifferential and subdifferential from 
\cref{def:directional_proximal_subdifferentials} of the associated indicator function, respectively.
\begin{example}\label{ex:proximal_subdifferential_of_indicator}
	Let $S\subset\R^n$ be a closed set and fix $\bar z\in S$.
	One can easily check that for arbitrary $w\in\R^n$, we have
	\[
		|\dr\delta_S(\bar z)(w)|<\infty
		\quad\Longleftrightarrow\quad
		\dr\delta_S(\bar z)(w)=0
		\quad\Longleftrightarrow\quad
		w\in T_S(\bar z),
	\]
	and $\{z^*\in\R^n\,|\,\dr\delta_S(\bar z)(w)=\skalp{z^*,w}\}=\{w\}^\perp$ follows for each $w\in T_S(\bar z)$.
	Putting this together with \cref{lem:proximal_pre-normal_cone,lem:proximal_pre_subdifferential}, we immediately see
	\[
		\widetilde\partial^p\delta_S(\bar z;w)=\Np_S(\bar z;w),\qquad
		\widehat\partial^p\delta_S(\bar z;w)=\widehat N^p_S(\bar z;w)
		\qquad\forall w\in T_S(\bar z),
	\]
	and $\widetilde\partial^p\delta_S(\bar z;w)=\widehat\partial^p\delta_S(\bar z;w)=\varnothing$ if $w\notin T_S(\bar z)$.
\end{example}

Subsequently, we state an analogy of \cref{prop:basic_properties_second_subderivative_indicator} which is an
immediate consequence of \eqref{eq:trivial_infinity_of_second_subderivative} and \cref{def:directional_proximal_subdifferentials}.

\begin{proposition}\label{prop:basic_properties_second_subderivative}
	Given a lower semicontinuous function $h\colon\R^n \to \overline \R$, a point $\bar z \in \R^n$ such that $|h(\bar z)|<\infty$,
	$z^*\in\R^n$, and a direction $w\in \R^n$, the following statements hold.
	\begin{enumerate}
  		\item\label{item:+infty} 
  			If ${\dr}h(\bar z)(w) > \skalp{z^*,w}$, then ${\dr}^2h(\bar z;z^*)(w) = \infty$.
  		\item\label{item:-infty}
			If ${\dr}h(\bar z)(w) < \infty$, we have ${\dr}^2h(\bar z;z^*)(w)>-\infty$ if and only if $z^*\in \widetilde\partial^p h(\bar z;w)$.
  		\item If ${\dr}^2h(\bar z;z^*)(w)$ is finite, then $z^*\in \widehat \partial^p h(\bar z;w)$.
	\end{enumerate}
\end{proposition}

Observe that \cref{ex:proximal_subdifferential_of_indicator} and \cref{prop:basic_properties_second_subderivative}
precisely recover \cref{prop:basic_properties_second_subderivative_indicator}.

The upcoming results discuss the calculus of second subderivatives of indicator
functions associated with graphs and epigraphs of single-valued mappings.
Furthermore, we carve out the role of the directional proximal subdifferential
in this context.

\begin{proposition}\label{prop:second_derivative_mapping}
	For a continuous function $F\colon\R^n \to \R^m$, $(\bar x,\bar y) \in \gph F$, 
	and pairs $(x^*,y^*), (u,v) \in \R^n \times \R^m$, we have
	\[
		\dr^2\delta_{\gph F}((\bar x,\bar y);(x^*,y^*))(u,v)
		=
		\dr^2\delta_{\gph F^{-1}}((\bar y,\bar x);(y^*,x^*))(v,u)
		\geq
		\dr^2  \langle -y^*,F\rangle (\bar x;x^*)(u),
	\]
	and the last estimate holds as equality	for some $v \in DF(\bar x)(u)$
	with $\dr \skalp{ -y^*,F } (\bar x)(u) = \skalp{-y^*,v}$
	whenever $F$ is calm at $\bar x$ in direction $u$.
\end{proposition}
\begin{proof}
	The first equality is trivial.
	Observing that $\bar y=F(\bar x)$, consider sequences $t_k \downarrow 0$ and $(u_k,v_k) \to (u,v)$ with $F(\bar x) + t_k v_k = F(\bar x + t_k u_k)$ 
	such that
	\begin{align*}
		\dr^2\delta_{\gph F}((\bar x,\bar y);(x^*,y^*))(u,v)
		=
		\lim_{k \to \infty} \frac{-2\langle (x^*,y^*), (u_k, v_k) \rangle }{t_k}.
	\end{align*}
	Then we have
	\begin{align*}
		\dr^2\delta_{\gph F}((\bar x,\bar y);(x^*,y^*))(u,v)
		& = 
		\lim_{k \to \infty} \frac{\skalp{ -y^*,F }(\bar x + t_k u_k) - \skalp{-y^*,F }(\bar x)  - t_k\skalp{ x^*,u_k }}{\frac{1}{2} t_k^2}
		\\
		&
		\geq
		\dr^2  \skalp{ -y^*,F } (\bar x;x^*)(u).
	\end{align*}
	If such sequences $t_k\downarrow 0$ and $(u_k,v_k)\to(u,v)$ do not exist, 
	$(u,v)\notin T_{\gph F}(\bar x,\bar y)$ is valid and the inequality holds trivially.
	
	To prove the converse inequality, consider $t_k \downarrow 0$ and $u_k \to u$ 
	which recover $\dr\skalp{-y^*,F}(\bar x)(u)$ and $\dr^2\skalp{-y^*,F}(\bar x;x^*)(u)$ simultaneously,
	see \cref{lem:special_sequences},
	and set $v_k := (F(\bar x + t_k u_k) - F(\bar x))/t_k$ for each $k\in\N$.
	By the assumed calmness of $F$ at $\bar x$ in direction $u$,
	we know that there is $v \in DF(\bar x)(u)$ satisfying $\dr \skalp{ -y^*,F } (\bar x)(u) = \skalp{-y^*,v}$
	such that by passing to a subsequence (without relabeling) we may assume $v_k \to v$.
	Thus,
	\begin{align*}
		\dr^2  \skalp{ -y^*,F } (\bar x;x^*)(u)
		=
		\lim_{k \to \infty} \frac{-2\langle (x^*,y^*), (u_k, v_k) \rangle }{t_k}
		\geq
		\dr^2\delta_{\gph F}((\bar x,\bar y);(x^*,y^*))(u,v)
	\end{align*}
	and the proof is completed.
\end{proof}

Given the connections between second subderivatives and directional proximal normal cones as well as subdifferentials,
see \cref{prop:basic_properties_second_subderivative_indicator} as well as
\cref{def:directional_proximal_subdifferentials} and \cref{prop:basic_properties_second_subderivative},
estimates for second subderivatives automatically contain certain estimates for the 
directional proximal normals and subdifferentials.
We demonstrate this in the following corollary.

\begin{corollary}
    In the setting of \cref{prop:second_derivative_mapping},
    let $v \in DF(\bar x)(u)$.
    We have
	\begin{align*}
        x^* \in \widetilde\partial^p \skalp{ -y^*,F } (\bar x;u)
        & \quad\Longrightarrow\quad
        (x^*,y^*) \in \Np_{\gph F}((\bar x,\bar y);(u,v))\\
        %\quad \forall\, v \in DF(\bar x)(u)\\
        x^* \in \widehat \partial^p \skalp{ -y^*,F } (\bar x;u)
        & \quad\Longrightarrow\quad
        (x^*,y^*) \in \widehat N^p_{\gph F}((\bar x,\bar y);(u,v))
        %\quad \forall\, v \in DF(\bar x)(u)
	 \textrm{ provided } \dr \skalp{ -y^*,F } (\bar x)(u) = \skalp{-y^*,v}
    \end{align*}
	and the reverse implications also hold for some $v \in DF(\bar x)(u)$ 
	with $\dr \skalp{ -y^*,F } (\bar x)(u) = \skalp{-y^*,v}$
	whenever $F$ is calm at $\bar x$ in direction $u$.
\end{corollary}
\begin{proof}
	For $x^* \in \widetilde\partial^p \skalp{ -y^*,F } (\bar x;u)$,
    \cref{def:directional_proximal_subdifferentials} and \cref{prop:second_derivative_mapping}
	imply
	\[
        	- \infty
        	<
        	\dr^2  \skalp{ -y^*,F } (\bar x;x^*)(u)
        	\leq
        	\dr^2\delta_{\gph F}((\bar x,\bar y);(x^*,y^*))(u,v)
	\]
	for all $v \in \R^m$.
	Consequently, \cref{prop:basic_properties_second_subderivative_indicator}
    	yields $(x^*,y^*) \in \Np_{\gph F}((\bar x,\bar y);(u,v))$ for all $v \in DF(\bar x)(u)$.
	If, additionally, $\skalp{x^*,u} = \dr \skalp{ -y^*,F } (\bar x)(u)$, i.e., $x^* \in \widehat\partial^p \skalp{ -y^*,F } (\bar x;u)$, then for $v$
	satisfying the relation $\dr \skalp{ -y^*,F } (\bar x)(u) = \skalp{-y^*,v}$, we get
	$\skalp{(x^*,y^*),(u,v)} = \skalp{x^*,u} - \skalp{-y^*,v} = 0$
	and $(x^*,y^*) \in \widehat N^p_{\gph F}((\bar x,\bar y);(u,v))$ follows.

	Next, let $F$ be calm at $\bar x$ in direction $u$.
	Consider $(x^*,y^*) \in \Np_{\gph F}((\bar x,\bar y);(u,v))$ with $v \in DF(\bar x)(u)$ and $\dr \skalp{ -y^*,F } (\bar x)(u) = \skalp{-y^*,v}$
	satisfying
	\[
		\dr^2\delta_{\gph F}((\bar x,\bar y);(x^*,y^*))(u,v)
		=
		\dr^2 \skalp{ -y^*,F } (\bar x;x^*)(u).
	\]
	Such $v$ exists by \cref{prop:second_derivative_mapping}.
	Both of these second subderivatives are greater than $- \infty$ by \cref{prop:basic_properties_second_subderivative_indicator},
	and \cref{prop:basic_properties_second_subderivative} implies
	$x^* \in \widetilde\partial^p \skalp{ -y^*,F } (\bar x;u)$.
	Again, if $\skalp{(x^*,y^*),(u,v)} = 0$, i.e., $(x^*,y^*)\in\widehat N_{\gph F}^p((\bar x,\bar y);(u,v))$, we get
	$\skalp{x^*,u} = \skalp{-y^*,v} = \dr \skalp{ -y^*,F } (\bar x)(u)$
	and $x^* \in \widehat \partial^p \skalp{ -y^*,F } (\bar x;u)$ follows.
\end{proof}

The above corollary can be viewed as a form of the scalarization formula for coderivatives of single-valued mappings, see \cite[Proposition~9.24]{RoWe98}.

\begin{proposition}\label{pro:second_derivative_via_indicator_to_epigraph}
    Let $h\colon \mathbb{R}^n \rightarrow \overline \R$ be a lower semicontinuous function
	and fix $\zb\in\R^n$ with $|h(\bar z)|<\infty$ as well as $z^*\in\R^n$ and $(w,\omega) \in \R^{n}\times\R$.
	We have
	\[
        {\dr}^2 \delta_{\epi h}((\bar z,h(\bar z));(z^*,-1))(w,\omega)
        \geq
 		{\dr}^2h(\bar z;z^*)(w),
	\]
	and the estimate holds as equality for each $\omega$ satisfying
	\[
		\omega \in 
		\begin{cases}
			\R & \textrm{if } {\dr}h(\bar z)(w) = \infty,\\
			(-\infty, \langle z^*, w \rangle) & \textrm{if } {\dr}h(\bar z)(w) = -\infty,\\
			\{{\dr}h(\bar z)(w)\} & \textrm{if } {\dr}h(\bar z)(w) \in \R.
		\end{cases}
	\]
\end{proposition}
\begin{proof}
	Consider sequences $t_k \downarrow 0$ and $(w_k,\omega_k) \to (w,\omega)$ with $(\bar z + t_k w_k,h(\bar z) + t_k \omega_k) \in \epi h$ 
	for each $k\in\N$ such that
	\[
		{\dr}^2 \delta_{\epi h}((\bar z,h(\bar z));(z^*,-1))(w,\omega)
		=
		\lim_{k \to \infty} \frac{-2\langle (z^*,-1), (w_k, \omega_k) \rangle }{t_k}
		=
		\lim_{k \to \infty} \frac{\omega_k -\langle z^*, w_k \rangle }{\frac{1}{2}t_k}.
	\]
	Due to $\omega_k\geq(h(\bar z+t_kw_k)-h(\bar z))/t_k$ for each $k\in\N$, we have
	\begin{align*}
		{\dr}^2 \delta_{\epi h}((\bar z,h(\bar z));(z^*,-1))(w,\omega)
		\geq 
		\lim_{k \to \infty} \frac{h(\bar z + t_k w_k) - h(\bar z) - t_k\langle z^*, w_k \rangle }{\frac{1}{2}t_k^2}
		\geq
		\dr^2 h (\bar z;z^*)(w).
	\end{align*}
	If such sequences $t_k\downarrow 0$ and $(w_k,\omega_k)\to(w,\omega)$ do not exist, 
	$(w,\omega)\notin T_{\epi h}(\bar z,h(\bar z))$ is valid and the inequality holds trivially.
	
	The converse relation will be shown by a distinction of cases.
	If ${\dr}h(\bar z)(w) = \infty$, we get $\dr h(\bar z)(w)>\skalp{z^*,w}$ and
	\[
		{\dr}^2 \delta_{\epi h}((\bar z,h(\bar z));(z^*,-1))(w,\omega)
		\geq
		\dr^2 h (\bar z;z^*)(w)
		=
		\infty
	\]
	for each $\omega \in \R$
	from the first part of the proof and \eqref{eq:trivial_infinity_of_second_subderivative}.
	If ${\dr}h(\bar z)(w) = -\infty$, we find sequences $t_k'\downarrow 0$ and $w_k'\to w$ such that
	$\omega_k':=(h(\bar z+t_k'w_k')-h(\bar z))/t_k'\to-\infty$.
	Thus, for each $\omega<\skalp{z^*,w}$, we have $\omega_k'<\omega$ and, thus, $(\bar z+t_k'w_k',h(\bar z)+t_k'\omega)\in\epi h$ for
	large enough $k\in\N$.
	Together with \eqref{eq:trivial_infinity_of_second_subderivative}, this gives
	\begin{align*}
		{\dr}^2 \delta_{\epi h}((\bar z,h(\bar z));(z^*,-1))(w,\omega)
		&\leq 
		\liminf_{k \to \infty} \frac{\omega -\langle z^*, w_k' \rangle }{\frac{1}{2}t_k'}
		=
		-\infty
		=
		\dr^2h(\bar z;z^*)(w).
	\end{align*}
	Finally, if ${\dr}h(\bar z)(w) \in \R$
	let us pick sequences $t_k\downarrow 0$ and $w_k\to w$ 
	recovering $\dr h(\bar z)(w)$ and $\dr^2 h(\bar z;z^*)(w)$ simultaneously, see \cref{lem:special_sequences}.
	Setting $\omega_k:=(h(\bar z+t_kw_k)-h(\bar z))/t_k$ for each $k\in\N$, we obtain
	\[
		\dr^2 h (\bar z;z^*)(w)
		=
		\lim_{k \to \infty} \frac{-2\langle (z^*,-1), (w_k, \omega_k) \rangle }{t_k}
		\geq
		{\dr}^2 \delta_{\epi h}((\bar z,h(\bar z));(z^*,-1))(w,{\dr}h(\bar z)(w)),
	\]
	and this completes the proof.
\end{proof}

Taking into account \cref{prop:basic_properties_second_subderivative_indicator},
the above result actually clarifies the close relationship between the directional proximal subdifferential
of a function and its geometric counterpart.
\begin{corollary}\label{lem:geometric_vs_standard_proximal_subdifferential}
	Given a lower semicontinuous function $h\colon\R^n\to\overline \R$, a point $\bar z\in\R^n$ such that $|h(\bar z)|<\infty$,
	and a direction $w\in\R^n$ such that $|\dr h(\bar z)(w)| < \infty$, we have
	\[
		\widetilde \partial^p h(\bar z;w)
		=
		\widetilde \partial^p_g h(\bar z;(w,{\dr}h(\bar z)(w))),
		\quad
		\widehat \partial^p h(\bar z;w)
		=
		\widehat \partial^p_g h(\bar z;(w,{\dr}h(\bar z)(w))).
	\]
	If $h$ is convex, we have
	\[
		\widehat \partial^p h(\bar z;w)
		=
		\partial h(\bar z) \cap \{z^* \in \R^n \mv \skalp{z^*,w} = {\dr}h(\bar z)(w)\}.
	\]
\end{corollary}
\begin{proof}
	Since
	\[
		\widehat \partial^p_g h(\bar z;(w,{\dr}h(\bar z)(w)))
		=
		\widetilde \partial^p_g h(\bar z;(w,{\dr}h(\bar z)(w)))
		\cap
		\{z^*\in\R^n\,|\,\dr h(\bar z)(w) = \langle z^*,w\rangle\},
	\]
	the second statement follows immediately from the first one.

	Let $z^* \in \widetilde \partial^p h(\bar z;w)$.
	Since $\dr h(\bar z)(w) \in \R$,
	\cref{def:directional_proximal_subdifferentials} and \cref{pro:second_derivative_via_indicator_to_epigraph} yield the relations
	\[
		-\infty 
		< 
		{\dr}^2h(\bar z;z^*)(w) 
		= 
		{\dr}^2 \delta_{\epi h}((\bar z,h(\bar z));(z^*,-1))(w,\dr h(\bar z)(w)),
	\]
	and since $(w,\dr h(\bar z)(w)) \in \epi \dr h(\bar z) = T_{\epi h}(\bar z,h(\bar z))$,  \cref{prop:basic_properties_second_subderivative_indicator}
	implies $z^* \in \widetilde \partial^p_g h(\bar z;(w,{\dr}h(\bar z)(w)))$.

	The opposite inclusion follows from the same arguments.
		The convex case follows from the definition
		of $\widehat \partial^p_g h(\bar z;(w,{\dr}h(\bar z)(w)))$ and \eqref{EqDirProxNormalConeConvex}.
\end{proof}

\section{Calculus for second subderivatives}\label{sec:calculus}

This section is devoted to the calculus of second subderivatives.
First, we propose two very general calculus rules for second subderivatives, namely a chain rule,
i.e., the rule for compositions, and a rule for marginal functions.
Afterwards, we apply these rules to derive some other calculus principles for the second subderivative
of set indicators where the set is given as the smooth image or pre-image of a closed set.

\subsection{Compositions and marginal functions}

Let us start with the consideration of a very general chain rule.

\begin{theorem}\label{The:Composition_nonsmooth}
Consider a lower semicontinuous function $g\colon\R^m \to \overline \R$
and a continuous mapping $F\colon \R^n \to \R^m$,
and, for $h\colon\R^n \to \overline \R$ given by $h := g \circ F$,
let $\bar x \in \R^n$ be chosen such that $|h(\bar x)|<\infty$.
If $F$ is calm at $\bar x$ in direction $u\in\R^n$,
then there exist $v \in DF(\bar x)(u)$ such that
for each $x^*\in\R^n$, one has
\begin{align*}
	\dr h(\bar x)(u)
    &  \geq
    \dr \skalp{y^*,F}(\bar x)(u)
    +
    \dr g(F(\bar x))(v)
    -
    \skalp{y^*,v},
    \\
    {\dr}^2h(\bar x;x^*)(u)
    &\geq
    {\dr}^2 \skalp{y^*,F}(\bar x;x^*)(u)
    +
    {\dr}^2 g(F(\bar x);y^*)(v)
\end{align*}
for all $y^* \in \R^m$ such that the right-hand side does not contain the summands $-\infty$ and $\infty$ simultaneously.
\end{theorem}
\begin{proof}
 By \cref{lem:special_sequences}, there are sequences $t_k \downarrow 0$ and $u_k \to u$ which recover
 $\dr h(\bar x)(u)$ and $\dr^2h(\bar x;x^*)(u)$ simultaneously. For each $k\in\N$, let us set
 $v_k := (F(\bar x+t_ku_k) - F(\bar x))/t_k$. Due to the postulated calmness assumption, we may assume
 $v_k\to v\in DF(\bar x)(u)$. Thus, we find
 \begin{align*}
 	\dr h(\bar x)(u)
 	&=
 	\lim\limits_{k\to\infty}\frac{g(F(\bar x+t_ku_k))-g(F(\bar x))}{t_k}
 	\\
 	&=
 	\lim\limits_{k\to\infty}
 		\left(
 			\frac{g(F(\bar x)+t_kv_k)-g(F(\bar x))}{t_k}
 			+
 			\frac{\skalp{y^*,F}(\bar x+t_ku_k)-\skalp{y^*,F}(\bar x)}{t_k}
 			-
 			\skalp{y^*,v_k}
 		\right)
 	\\
 	&\geq
 	\dr g(F(\bar x))(v)+\dr\skalp{y^*,F}(\bar x)(u)-\skalp{y^*,v}
 \end{align*}
 and
 \begin{align*}
        &{\dr}^2h(\bar x;x^*)(u) 
         = 
        \lim_{k \to \infty} \frac{g(F(\bar x+t_ku_k))-g(F(\bar x))-t_k \langle x^*, u_k\rangle}{\frac{1}{2}t_k^2}\\
        &\qquad
         =
        \lim_{k \to \infty}
        \left(
        \frac{g(F(\bar x)+t_k v_k)-g(F(\bar x))-t_k \langle y^*, v_k\rangle}{\frac{1}{2}t_k^2}
        +
        \frac{\skalp{y^*,F}(\bar x+t_ku_k)-\skalp{y^*,F}(\bar x)-t_k \langle x^*, u_k\rangle}{\frac{1}{2}t_k^2}
        \right)
        \\
        &\qquad\geq
        \dr^2g(F(\bar x);y^*)(v)
        +
        \dr^2\skalp{y^*,F}(\bar x;x^*)(u)
 \end{align*}
    where the last inequality follows since $\liminf_{k\to\infty} (a_k + b_k) \geq \liminf_{k\to\infty} a_k + \liminf_{k\to\infty} b_k$ holds
    except for the indeterminate case $\infty - \infty$, respectively.
\end{proof}

\begin{corollary}\label{Cor:Chian_rule_directional_proximal_subdifferentials}
	In the setting of \cref{The:Composition_nonsmooth},
	assume that $\dr h(\bar x)(u)<\infty$.
	Then there exists $v \in DF(\bar x)(u)$ such that for each $x^*\in\R^n$ and $y^*\in \R^m$, we have
	the following implications:
	\begin{align*}
		y^* \in \widetilde\partial^p g(F(\bar x);v),\
        x^* \in \widetilde\partial^p \skalp{y^*,F } (\bar x;u)
        & \quad\Longrightarrow\quad
        x^* \in \widetilde\partial^p h (\bar x;u),
        \\
        \dr h(\bar x)(u)=\dr\skalp{y^*,F}(\bar x)(u),\
        y^* \in \widetilde\partial^p g(F(\bar x);v),\
        x^* \in \widehat\partial^p \skalp{y^*,F } (\bar x;u)
        & \quad\Longrightarrow\quad
        x^* \in \widehat\partial^p h (\bar x;u).
    	\end{align*}
\end{corollary}
\begin{proof}
	For $y^* \in \widetilde\partial^p g(F(\bar x);v)$ and 
	$x^* \in \widetilde\partial^p \skalp{y^*,F } (\bar x;u)$,
	we get the estimates ${\dr}^2 \skalp{y^*,F}(\bar x;x^*)(u) > -\infty$ and ${\dr}^2 g(F(\bar x);y^*)(v) > - \infty$
	as well as $|\dr\skalp{y^*,F}(\bar x)(u)|<\infty$ and $|\dr g(F(\bar x))(v)|<\infty$
	by \cref{def:directional_proximal_subdifferentials}.
	Consequently, the estimates from \cref{The:Composition_nonsmooth} apply and yield
	that $\dr h(\bar x)(u)>-\infty$ and ${\dr}^2h(\bar x;x^*)(u) > -\infty$.
	This in turn gives
	$x^* \in \widetilde\partial^p h (\bar x;u)$ since we assumed ${\dr}h(\bar x)(u) < \infty$.
	
	In case where $\dr h(\bar x)(u)=\dr\skalp{y^*,F}(\bar x)(u)$, $y^* \in \widetilde\partial^p g(F(\bar x);v)$, 
	and $x^* \in \widehat\partial^p \skalp{y^*,F } (\bar x;u)$, we can deduce $x^*\in\widetilde\partial^ph(\bar x;u)$ as
	above. Furthermore, due to $\dr\skalp{y^*,F}(\bar x)(u)=\skalp{x^*,u}$, we have $\dr h(\bar x)(u)=\skalp{x^*,u}$,
	and $x^*\in\widehat\partial^p h(\bar x;u)$ follows.
\end{proof}

Under the assumptions of \cref{Cor:Chian_rule_directional_proximal_subdifferentials},
the second implication also gives
\[
        \left.
        \begin{aligned}
        &\dr h(\bar x)(u)=\dr\skalp{y^*,F}(\bar x)(u)+\dr g(F(\bar x))(v)-\skalp{y^*,v},\\
        &y^* \in \widehat\partial^p g(F(\bar x);v),\
        x^* \in \widehat\partial^p \skalp{y^*,F } (\bar x;u)
        \end{aligned}
        \right\}
         \quad\Longrightarrow\quad
        x^* \in \widehat\partial^p h (\bar x;u),
\]
which demonstrates the close relationship to the subderivative chain rule from 
\cref{The:Composition_nonsmooth}.

If $F$ is twice continuously differentiable, we get the following corollary from \cref{The:Composition_nonsmooth} by taking into account
\[
	DF(\bar x)(u)=\{\nabla F(\bar x)u\},
	\quad
	\nabla \skalp{y^*,F}(\bar x)
	=
	\nabla F(\bar x)^\top y^*,
	\quad
	\nabla^2\skalp{y^*,F}(\bar x)(u,u)
	=
	\skalp{y^*,\nabla ^2F(\bar x)(u,u)},
\]
local Lipschitzness of $F$ around $\bar x$, and \eqref{eq:second_subderivative_smooth_map}.

\begin{corollary}\label{The:Composition}
Consider a lower semicontinuous function $g\colon\R^m \to \overline \R$
and a twice continuously differentiable mapping $F\colon \R^n \to \R^m$,
and, for $h\colon\R^n \to \overline \R$ given by $h := g \circ F$,
let $\bar x \in \R^n$ be chosen such that $|h(\bar x)|<\infty$.
Then for each $x^*\in\R^n$ and $u\in\R^n$, one has
\begin{align*}
	\dr h(\bar x)(u)
	&\geq
	\dr g(F(\bar x))(\nabla F(\bar x)u),
	\\
    {\dr}^2h(\bar x;x^*)(u)
    &\geq
    \sup_{\nabla F(\bar x)^{\top} y^* = x^*}
    \left( \skalp{y^*,\nabla^2 F(\bar x)(u,u)} 
    		+ 
    		{\dr}^2 g(F(\bar x);y^*)(\nabla F(\bar x)u) 
    \right).
\end{align*}
If $\nabla F(\bar x)$ possesses full row rank, for each $y^*\in\R^m$ and $u\in\R^n$, one has
\begin{align*}
	\dr h(\bar x)(u)
	&=
	\dr g(F(\bar x))(\nabla F(\bar x)u),
	\\
	{\dr}^2h(\bar x;\nabla F(\bar x)^\top y^*)(u)
    &=
    \skalp{y^*,\nabla^2 F(\bar x)(u,u)} 
    + 
    {\dr}^2 g(F(\bar x);y^*)(\nabla F(\bar x)u).
\end{align*}
\end{corollary}
\begin{proof}
	The lower estimates are a direct consequence of \cref{The:Composition_nonsmooth}.
	In order to show the second statement, choose sequences
	$t_k\downarrow 0$ and $v_k\to\nabla F(\bar x)u$ which recover $\dr g(F(\bar x))(\nabla F(\bar x)u)$
	and $\dr ^2 g(F(\bar x);y^*)(\nabla F(\bar x)u)$ simultaneously, see \cref{lem:special_sequences}.
	Furthermore, set $x^*:=\nabla F(\bar x)^\top y^*$.
	Since $\nabla F(\bar x)$ possesses full row rank,
	\cite[Exercise~9.44]{RoWe98} implies the existence of a constant $\kappa>0$
	and neighborhoods $U\subset\R^n$ of $\bar x$ and $V\subset\R^m$ of $F(\bar x)$ such that
	\[
		\dist(x,F^{-1}(y))\leq\kappa\dist(y,F(x))\qquad\forall x\in U,\,\forall y\in V,
	\]
	i.e., $F$ is so-called metrically regular at $(\bar x,F(\bar x))$. 
	For sufficiently large $k\in\N$, we may apply this estimate with
	$x:=\bar x+t_ku$ and $y:=F(\bar x)+t_kv_k$ in order to find $x_k\in\R^n$
	such that $F(x_k)=F(\bar x)+t_kv_k$ and 
	$\norm{x_k-\bar x-t_ku}\leq\kappa\norm{F(\bar x+t_ku)-F(\bar x)-t_kv_k}$.
	Let us set $u_k:=(x_k-\bar x)/t_k$ for each $k\in\N$ sufficiently large.
	Then we find
	\begin{align*}
		\norm{u_k-u}
		\leq
		\kappa\left\Vert\frac{F(\bar x+t_ku)-F(\bar x)}{t_k}-v_k\right\Vert
		\to
		0
	\end{align*}
	from $v_k\to\nabla F(\bar x)u$, i.e., $u_k\to u$.
	Thus, we can exploit
	$x_k=\bar x+t_ku_k$ for all sufficiently large $k\in\N$ in order to find
	\begin{align*}
		\dr g(F(\bar x))(\nabla F(\bar x)u)
		=
		\lim\limits_{k\to\infty}\frac{g(F(\bar x+t_ku_k))-g(F(\bar x))}{t_k}
		\geq
		\dr h(\bar x)(u)
	\end{align*}
	and
	\begin{align*}
		\dr^2g(F(\bar x);y^*)(\nabla F(\bar x)u)
		&=
		\lim\limits_{k\to\infty}
		\frac{g(F(\bar x+t_ku_k))-g(F(\bar x))-t_k\skalp{y^*,v_k}}{\frac12t_k^2}
		\\
		&=
		\lim\limits_{k\to\infty}
		\left(
			\frac{g(F(\bar x+t_ku_k))-g(F(\bar x))-t_k\skalp{x^*,u_k}}{\frac{1}{2}t_k^2}
			+
			\frac{\skalp{x^*,t_ku_k}-\skalp{y^*,t_kv_k}}{\frac{1}{2}t_k^2}
		\right)
		\\
		&\geq
		\dr^2h(\bar x;x^*)(u)
		+
		\lim\limits_{k\to\infty}
		\frac{\skalp{y^*,\nabla F(\bar x)(x_k-\bar x)+F(\bar x)-F(x_k)}}{\frac12t_k^2}
		\\
		&=
		\dr^2h(\bar x;x^*)(u)
		+
		\lim\limits_{k\to\infty}
		\skalp{y^*,-\nabla^2F(\bar x)(u_k,u_k)}
		\\
		&=
		\dr^2h(\bar x;x^*)(u)-\skalp{y^*,\nabla^2F(\bar x)(u,u)},
	\end{align*}
	where we used a second-order Taylor expansion of $F$ at $\bar x$ in the last but one equality.
	Noting that the converse relations hold due to the general estimates, the proof
	is complete.
\end{proof}

Let us note that the result in \cref{The:Composition} is essentially different from
the chain rule which can be found in \cite[Theorem~13.14]{RoWe98}. Therein,
the authors exploit a less restrictive qualification condition than the full row rank
of $\nabla F(\bar x)$ in order to derive a general lower estimate of the second
subderivative, and in order to
get equality, they additionally assume that $g$ is a so-called fully amendable function locally
around $F(\bar x)$, i.e., the composition of a twice continuously differentiable inner
and a piecewise linear-quadratic outer function. In contrast, \cref{The:Composition} 
yields a general lower estimate even in the absence of a qualification condition.
Since an upper estimate is not so important for our purposes,
we show equality only under the comparatively strong full rank condition without focusing on minimal assumptions.

\begin{corollary}\label{Cor:Sum_rule}
	For two lower semicontinuous functions $h_1,h_2\colon \R^{n} \to \overline\R$,
	$h\colon \R^n \to \overline\R$ defined by $h(x) := h_1(x) + h_2(x)$ for all $x\in\R^n$,
	some point $\bar x \in\R^n$ with $|h(\bar x)|<\infty$,
	and $x^*, u \in\R^n$, we have
	\[
		\dr^2 h(\bar x;x^*)(u)
		\geq
		\dr^2 h_1(\bar x;x^*_1)(u)+\dr^2 h_2(\bar x;x^*_2)(u)
	\]
	for all $x^*_1,x^*_2 \in \R^n$ such that
	$x^*_1 + x^*_2 = x^*$
	provided the right-hand side is not $\infty - \infty$.
\end{corollary}
\begin{proof}
	For the proof, we first apply the chain rule from \cref{The:Composition} 
	to $F(x) := (x,x)$, $x\in\R^n$, and $g(x_1,x_2) := h_1(x_1) + h_2(x_2)$, $x_1,x_2\in\R^n$,
	and then \cref{prop:sum_rules}~\ref{item:separable_sum}.
\end{proof}

Next, we will study a marginal function rule for the second subderivative.
Recall from \cref{sec:svm} that we choose a proper, lower semicontinuous function 
$\varphi\colon\R^n\times\R^m\to\overline\R$
and consider the associated marginal function $h\colon\R^n\to\overline\R$ 
given as in \eqref{eq:marginal_function}
together with the set-valued mappings $\Upsilon\colon\R^n\times\R\tto\R^m$
and $\Psi\colon\R^n\tto\R^m$ from \eqref{eq:M_and_S}.
Recall that $\Psi(x):=\varnothing$ is used for each $x\in\R^n$ such that
$\varphi(x,y)=\infty$ holds for all $y\in\R^m$.
Furthermore, we would like to mention again that inner semicompactness
of $\Upsilon$ or $\Psi$ yields lower semicontinuity of $h$ so that the
consideration of subderivatives is reasonable, see \cref{lem:isc_of_M,lem:Psi_vs_M}.

In the next theorem, we address the second subderivative of marginal functions with the aid of the solution mapping $\Psi$.

\begin{theorem}\label{The:Marginal}
    Consider a proper, lower semicontinuous function $\varphi\colon\R^n \times \R^m \to \overline \R$
    and fix $\bar x\in\dom \Psi$ for the mapping $\Psi\colon\R^n\tto\R^m$ 
    given in \eqref{eq:M_and_S}, and let $\Psi$ be inner semicompact at $\bar x$.
    Then for each $x^*\in\R^n$ and $u\in\R^n$, one has
    \begin{subequations}\label{eq:forward_rule}
	\begin{align}
			\label{eq:forward_rule_subderivative}
    		{\dr} h(\bar x)(u)
        	& \leq
        	\inf_{y \in \Psi(\bar x), v \in \R^m}
        	{\dr} \varphi(\bar x,y)(u,v),\\
        	\label{eq:forward_rule_second_subderivative}
		{\dr}^2 h(\bar x;x^*)(u)
        	& \leq
        	\inf_{y \in \Psi(\bar x), v \in \R^m}
        	{\dr}^2 \varphi((\bar x,y);(x^*,0))(u,v)
    	\end{align}
    \end{subequations}
    where $h\colon\R^n \to \overline \R$ is the marginal function
    defined in \eqref{eq:marginal_function}.
    
    On the other hand, suppose that $\Psi$ is inner calm* at $\bar x$ in direction $u$.
    Then the estimates \eqref{eq:forward_rule} hold as equalities,
    and whenever $\dr h(\bar x)(u)$ is finite,
    both infima therein are attained at some pair $(y,v) \in \R^m \times \R^m$ 
    with $y \in \Psi(\bar x)$ and $v \in D\Psi(\bar x,y)(u)$.
\end{theorem}
\begin{proof}
	For arbitrary $y\in \Psi(\bar x)$ and $v \in \R^m$,
	by \cref{lem:special_sequences},
    	there are sequences $t_k \downarrow 0$ and $(u_k,v_k) \to (u,v)$ satisfying
    	\begin{align*}
	{\dr} \varphi(\bar x,y)(u,v)
        & = 
        \lim_{k \to \infty} \frac{\varphi((\bar x,y)+t_k(u_k,v_k))-\varphi(\bar x,y)}{t_k}\\
        & \geq 
        \lim_{k \to \infty} \frac{h(\bar x + t_k u_k)-h(\bar x)}{t_k}
        \geq
        {\dr}h(\bar x)(u),\\
        {\dr}^2 \varphi((\bar x,y);(x^*,0))(u,v)
        & = 
        \lim_{k \to \infty} \frac{\varphi((\bar x,y)+t_k(u_k,v_k))-\varphi(\bar x,y)-t_k \langle (x^*,0), (u_k,v_k)\rangle}{\frac{1}{2}t_k^2}\\
        & \geq 
        \lim_{k \to \infty} \frac{h(\bar x + t_k u_k)-h(\bar x)-t_k \langle x^*, u_k\rangle}{\frac{1}{2}t_k^2}
        \geq
        {\dr}^2h(\bar x;x^*)(u)
    	\end{align*}
    	since $h(\bar x) = \varphi(\bar x,y)$. This shows the first claim.

    To prove the second claim, 
	consider sequences
    $t_k \downarrow 0$ and $u_k \to u$ which recover ${\dr}h(\bar x)(u)$ and ${\dr}^2h(\bar x;x^*)(u)$ simultaneously, see \cref{lem:special_sequences}.
	By passing to a subsequence (without relabeling), 
	the assumed inner calmness* of $\Psi$ in direction $u$ yields the existence of $\kappa > 0$, a point $y\in\R^m$,
    and a sequence $y_k\to y$ such that
    \begin{subequations}\label{eq:consequence_of_ic*}
    	\begin{align}
    		\label{eq:Sol_Map_Seq}
					&y_k \in\Psi(\bar x+t_ku_k),\\
			\label{eq:estimate_ic*}
					&\norm{y_k - y} \leq t_k \kappa \norm{u_k}
    	\end{align}
    \end{subequations}
    hold for each $k\in\N$.
    Due to \eqref{eq:estimate_ic*}, $v_k:=(y_k - y)/t_k$ remains bounded as $k\to\infty$, 
    and we may assume $v_k\to v$ for some $v \in \R^m$ 
    with $\norm{v} \leq \kappa \norm{u}$.
    
	Again, by passing to a subsequence (without relabeling), we may assume that $h(\bar x+t_ku_k) \to \alpha \geq h(\bar x)$,
	taking into account lower semicontinuity of $h$, see
	\cref{lem:isc_of_M,lem:Psi_vs_M}.
	Moreover, if $\alpha > h(\bar x)$, we get $\dr h(\bar x)(u) = \dr^2h(\bar x;x^*)(u) = \infty$ and the converse inequalities in \eqref{eq:forward_rule} are trivial.
	Thus, we assume that $\varphi(\bar x+t_ku_k,y_k) = h(\bar x+t_ku_k) \to h(\bar x)$, and lower semicontinuity of $\varphi$ gives
    $\varphi(\bar x,y) \leq h(\bar x)$, i.e., $y \in \Psi(\bar x)$.
	Moreover, \eqref{eq:Sol_Map_Seq} yields $v \in D\Psi(\bar x,y)(u)$
	since $y_k = y + t_k v_k$ holds for all $k\in\N$.
    Finally, we obtain
    \begin{align*}
	{\dr} h(\bar x)(u)
        & = 
        \lim_{k \to \infty} \frac{h(\bar x+t_ku_k)-h(\bar x)}{t_k}
        = 
        \lim_{k \to \infty} \frac{\varphi(\bar x+t_ku_k,y + t_k v_k)-\varphi(\bar x,y)}{t_k}
        \geq 
        {\dr} \varphi(\bar x,y)(u,v),\\
        {\dr}^2h(\bar x;x^*)(u)
        & = 
        \lim_{k \to \infty} \frac{h(\bar x+t_ku_k)-h(\bar x)-t_k \langle x^*, u_k\rangle}{\frac{1}{2}t_k^2} \\
        & = 
        \lim_{k \to \infty} \frac{\varphi(\bar x+t_ku_k,y + t_k v_k)-\varphi(\bar x,y)-t_k \langle (x^*,0), (u_k,v_k)\rangle}{\frac{1}{2}t_k^2}
        \geq 
        {\dr}^2 \varphi((\bar x,y);(x^*,0))(u,v),
    \end{align*}
    which completes the proof.
\end{proof}

\begin{corollary}\label{cor:hidden_calculus_marginal_functions}
    In the setting of \cref{The:Marginal}, let $|\dr h (\bar x)(u)| < \infty$.
    For all $y \in \Psi(\bar x)$ and $v \in \R^m$, one has
    \begin{align*}
        x^* \in \widetilde\partial^p h (\bar x;u)
        & \quad\Longrightarrow\quad
        (x^*,0) \in \widetilde\partial^p \varphi ((\bar x,y);(u,v))
        \quad \textrm{provided } \dr \varphi (\bar x,y)(u,v) < \infty\\
        %\quad \forall\, v \in DF(\bar x)(u)\\
        x^* \in \widehat \partial^p h (\bar x;u)
        & \quad\Longrightarrow\quad
        (x^*,0) \in \widehat\partial^p \varphi ((\bar x,y);(u,v))
        \quad \textrm{provided } \dr \varphi (\bar x,y)(u,v) = \dr h (\bar x)(u).
    \end{align*}
    Moreover, the reverse implications hold for some
    $y \in \Psi(\bar x)$ and $v \in D\Psi(\bar x,y)(u)$ with $\dr \varphi (\bar x,y)(u,v) = \dr h (\bar x)(u)$
    provided $\Psi$ is inner calm* at $\bar x$ in direction $u$.
\end{corollary}
\begin{proof}
	For $x^* \in \widetilde\partial^p h (\bar x;u)$,
    \cref{def:directional_proximal_subdifferentials,The:Marginal}
	imply
	\[
        	- \infty
        	<
        	{\dr}^2h(\bar x;x^*)(u)
        	\leq
        	{\dr}^2 \varphi((\bar x,y);(x^*,0))(u,v)
	\]
	for all $y \in \Psi(\bar x)$ and $v \in D\Psi(\bar x,y)(u)$,
	and this also gives $\dr \varphi(\bar x,y)(u,v)>-\infty$
	due to \eqref{eq:trivial_infinity_of_second_subderivative}.
	Consequently, $(x^*,0) \in \widetilde\partial^p \varphi ((\bar x,y);(u,v))$ if $\dr \varphi (\bar x,y)(u,v) < \infty$.
	If, additionally, $\skalp{x^*,u} = \dr h (\bar x)(u)$, i.e.,
	$x^*\in\widehat\partial^p h(\bar x;u)$, then for $y\in\Psi(\bar x)$ and
	$v\in\R^m$  satisfying $\dr \varphi (\bar x,y)(u,v) = \dr h (\bar x)(u)$, we get
	\begin{equation}\label{eq:4_equalities}
        \skalp{(x^*,0),(u,v)}
        =
        \skalp{x^*,u}
        =
        \dr h (\bar x)(u)
        =
        \dr \varphi (\bar x,y)(u,v)        
	\end{equation}
	and $(x^*,0) \in \widehat\partial^p \varphi ((\bar x,y);(u,v))$ follows.

	Let us now argue for the reverse implications.
	Consider now $(x^*,0) \in \widetilde\partial^p \varphi ((\bar x,y);(u,v))$
	with $y \in \Psi(\bar x)$ and $v \in D\Psi(\bar x,y)(u)$ satisfying
	\[
		\dr \varphi (\bar x,y)(u,v)
		=
		{\dr}h(\bar x)(u),
		\qquad
		{\dr}^2 \varphi((\bar x,y);(x^*,0))(u,v)
		=
		{\dr}^2h(\bar x;x^*)(u).
	\]
	Such $y$, $v$ exist by \cref{The:Marginal}.
	Both second subderivatives are thus greater than $- \infty$ by \cref{def:directional_proximal_subdifferentials},
	and so it also gives
	$x^* \in \widetilde\partial^p h (\bar x;u)$.
	If $\dr \varphi (\bar x,y)(u,v) = \skalp{(x^*,0),(u,v)}$, we again get the equalities
	\eqref{eq:4_equalities} and $x^* \in \widehat \partial^p h (\bar x;u)$ follows.
\end{proof}

To better demonstrate \cref{The:Marginal}, we use it to compute the second subderivative of the distance function to a closed set $\Omega \subset \R^n$,
given via the representation from \eqref{eq:distance_function_as_marginal_function}.
\begin{corollary}\label{Cor:Distance_function}
    For a closed set $\Omega \subset \R^n$, $\xb \in \Omega$, 
    $x^* \in \R^n$ with $\norm{x^*} \leq 1$, and $u\in T_\Omega(\xb)$, we have
    \[
        {\dr}^2 \dist(\cdot,\Omega)(\xb;x^*)(u) \geq {\dr}^2 \delta_\Omega (\xb;x^*)(u).
    \]
\end{corollary}
\begin{proof}
    First, we claim that
    \[
        {\dr}^2 \dist(\cdot,\Omega)(\xb;x^*)(u) = {\dr}^2 (\varphi_1 + \varphi_2)((\xb,\xb);(x^*,0))(u,u)
    \]
    for $\varphi_1(x,y) := \norm{y - x}$ and $\varphi_1(x,y) := \delta_\Omega(y)$,
    $(x,y)\in\R^n\times\R^n$.
    By \cref{Prop:Ic*Proj}, the mapping $\Psi := P_\Omega$ is inner calm* at $\xb$, and so we can apply \cref{The:Marginal}.
    Since, $P_\Omega(\xb) = \xb$, there is no need for the infimum over $\yb \in P_\Omega(\xb)$.
    Let us now argue why we can omit also the infimum over directions $v \in \R^n$ and take just $v=u$.
    Note that the proof of \cref{The:Marginal} yields that we can take $v$ with $\norm{v} \leq \kappa \norm{u}$
    for some $\kappa \geq 0$, which covers the case $u=0$.
    If $u\neq 0$, the proof of \cref{The:Marginal} uses \cref{lem:special_sequences} in order to only consider
    the sequences $t_k$ and $u_k$ which recover also ${\dr} \dist(\cdot,\Omega)(\xb)(u)$.
    This means that
    \[
    \dist(\xb + t_k u_k,\Omega)/(t_k \norm{u_k}) \to {\dr} \dist(\cdot,\Omega)(\xb)(u) / \norm{u} = 0
    \]
    due to $u\in T_\Omega(\xb)$ and \cite[Example 8.53]{RoWe98}.
    Thus, the second statement of \cref{Prop:Ic*Proj} yields that the sequence $v_k$ from the proof of \cref{The:Marginal}
    indeed converges to $u$ along a subsequence.
    
    Next, in order to estimate ${\dr}^2 (\varphi_1 + \varphi_2)((\xb,\xb);(x^*,0))(u,u)$, 
    we apply the sum rule from \cref{Cor:Sum_rule} first and then the chain rule
    from \cref{The:Composition} twice, taking into account $\varphi_1 = \norm{\cdot} \circ h_1$ and $\varphi_2 = \delta_\Omega \circ h_2$
    for the simple smooth mappings $h_1(x,y) := x - y$ and $h_2(x,y) := y$,
    $(x,y)\in\R^n\times\R^n$, with vanishing second derivatives.
    In the sum rule, we choose to split $(x^*,0)$ into $(x^*,-x^*) + (0,x^*)$ 
    due to the structure of $h_1$ and $h_2$.
    Consequently, taking into account \cref{lem:second_subderivative_of_euclidean_norm}, 
    we conclude
    \begin{align*}
        {\dr}^2 (\varphi_1 + \varphi_2)((\xb,\xb);(x^*,0))(u,u)
        & \geq 
        {\dr}^2 \varphi_1((\xb,\xb);(x^*,-x^*))(u,u) + {\dr}^2 \varphi_2((\xb,\xb);(0,x^*))(u,u)\\
        & = 
        {\dr}^2 \norm{\cdot}(0;x^*)(0) + {\dr}^2 \delta_\Omega(\xb;x^*)(u)
        =
        {\dr}^2 \delta_\Omega(\xb;x^*)(u)
    \end{align*}
	to complete the proof.
\end{proof}

	Let us now turn our focus on the general marginal function rule again.
With the additional assumption ${\dr} h(\bar x)(u) \geq \skalp{x^*,u}$,
we can derive estimates, similar to those ones in \cref{The:Marginal}, 
under the milder inner semicompactness/inner calmness* requirements imposed on $\Upsilon$, 
see \cref{lem:Psi_vs_M}.

\begin{theorem}\label{The:Marginal_2}
    Consider a proper, lower semicontinuous function $\varphi\colon\R^n \times \R^m \to \overline \R$
    and fix $\bar x \in \dom \Psi$ for the mapping $\Psi\colon\R^n\tto\R^m$ 
    given in \eqref{eq:M_and_S}.
    Suppose that $\Upsilon\colon\R^n\times\R\tto\R^m$ given in \eqref{eq:M_and_S}
	is inner semicompact at $(\bar x,h(\bar x))$ w.r.t.\ $\dom\Upsilon$.
    Then for each $x^*\in\R^n$ and $u\in\R^n$, the estimates \eqref{eq:forward_rule}
    hold.
    
    On the other hand, suppose that ${\dr} h(\bar x)(u) \geq \skalp{x^*,u}$.
    Then the following statements hold.
    \begin{enumerate}
    	\item If $\dr h(\bar x)(u)>\skalp{x^*,u}$, then equality holds in
    		\eqref{eq:forward_rule_second_subderivative}.
    	\item If ${\dr} h(\bar x)(u) = \skalp{x^*,u}$ and if
    		$\Upsilon$ is inner calm* at $(\bar x,h(\bar x))$ in direction 
    		$(u,\skalp{x^*,u})$ w.r.t.\ $\dom \Upsilon$,
    		then the estimates \eqref{eq:forward_rule} hold as equalities,
    		and whenever $\dr^2h(\bar x;x^*)(u)$ is finite,
    		both infima therein are attained at some pair $(y,v) \in \R^m \times \R^m$ 
    		with $y \in \Psi(\bar x)$ 
    		and $v \in D\Upsilon((\bar x,h(\bar x)),y)(u,\skalp{x^*,u})$.
    \end{enumerate}
\end{theorem}
\begin{proof}
	The first claim follows by the same arguments as in the previous theorem.

	To prove the second claim, consider sequences
    $t_k \downarrow 0$ and $u_k \to u$ 
    recovering the subderivatives $\dr h(\bar x)(u)$ and $\dr^2 h(\bar x;x^*)(u)$ simultaneously,
    see \cref{lem:special_sequences}.
	By passing to a subsequence (without relabeling) and taking into account 
	lower semicontinuity of $h$ again, see
	\cref{lem:isc_of_M}, we may assume that $h(\bar x+t_ku_k) \to \alpha \geq h(\bar x)$,
	and that
	\[
		\mu_k:=(h(\bar x+t_ku_k)-h(\bar x))/t_k \to {\dr} h(\bar x)(u) \geq \skalp{x^*,u}
	\]
	due to the postulated assumptions.
	If $\alpha > h(\bar x)$, we get $\dr h(\bar x)(u)=\dr^2h(\bar x;x^*)(u) = \infty$
	and the converse inequalities in \eqref{eq:forward_rule} follow trivially.
	If $\dr h(\bar x)(u)>\skalp{x^*,u}$, we find
	$\dr^2h(\bar x;x^*)(u)=\infty$, and the converse estimate in 
	\eqref{eq:forward_rule_second_subderivative} is trivial.
	Thus, we assume that $\alpha = h(\bar x)$ and ${\dr} h(\bar x)(u) = \skalp{x^*,u}$.
	
	The postulated inner semicompactness of $\Upsilon$ at $(\bar x,h(\bar x))$ yields that, locally around that point, $\epi h = \dom \Upsilon$
	and so $(\bar x+t_ku_k,h(\bar x+t_ku_k)) \in \dom \Upsilon$ for sufficiently large $k\in\N$,
	see \cref{lem:isc_of_M}.
	The assumed inner calmness* of $\Upsilon$ at $(\bar x,h(\bar x))$ in direction $(u,\skalp{x^*,u})$ then yields the existence of $\kappa > 0$,
	a point $y\in\R^m$, and a sequence $y_k\to y$ such that
    \begin{align*}
    	&y_k \in \Upsilon(\bar x+t_ku_k,h(\bar x)+t_k \mu_k))=\Psi(\bar x+t_ku_k),\\
		&\norm{y_k - y} \leq t_k \kappa \norm{(u_k,\mu_k)}
    \end{align*}
   	hold for each $k\in\N$. 
   	The remainder of the proof follows the same arguments as used to show
   	\cref{The:Marginal}.
\end{proof}

Similar as in \cref{cor:hidden_calculus_marginal_functions}, one can also derive hidden calculus results
on the directional proximal (pre-) subdifferential of $h$ in the setting of \cref{The:Marginal_2}.
For brevity, however, we leave this straightforward task to the interested reader.

Let us recall that the inner calmness* of $\Psi$ at $\bar x$ (or $\Upsilon$ at $(\bar x,h(\bar x))$) in \cref{The:Marginal} (or \cref{The:Marginal_2})
is inherent whenever $\Psi$ (or $\Upsilon$) is isolatedly calm at all points $(\bar x,y)\in\gph\Psi$ (or $((\bar x,h(\bar x)),y)\in\gph\Upsilon)$
since we already claimed that inner semicompactness holds at $\bar x$ (or $(\bar x,h(\bar x))$),
and isolated calmness can be checked in terms of the associated graphical derivative, i.e., in terms of problem data,
see e.g. \cite{GfrOut14a}.

\subsection{Other calculus rules}\label{sec:image_pre-image_rule}

In order to apply our results to interesting optimization problems, let us derive some additional calculus rules for second subderivatives.
More precisely, we look into the term ${\dr}^2\delta_{S}(x;x^*)(u)$ assuming that $S$ has a pre-image or image structure.

Applying \cref{The:Composition} with $h:=\delta_{C}$ for a closed set $C \subset \R^m$ immediately yields the following pre-image rule.

\begin{proposition}\label{prop:second_subderivative_pre-image}
	Consider a closed set $C \subset \R^m$ and a twice continuously differentiable mapping $F\colon \R^n \to \R^m$,
	and let $\bar x \in S:=F^{-1}(C)$.
	Then for each $x^*\in\R^n$ and $u\in\R^n$, one has
	\[
		{\dr}^2\delta_{S}(\bar x;x^*)(u)
		\geq
		\sup_{\nabla F(\bar x)^{\top} y^* = x^*} 
			\left( \skalp{y^*,\nabla^2 F(\bar x)(u,u)} + {\dr}^2\delta_{C}(F(\bar x);y^*)(\nabla F(\bar x)u) \right).
	\]
	If $\nabla F(\bar x)$ possesses full row rank,
	for each $y^*\in\R^m$ and $u\in\R^n$, one has
	\[
		{\dr}^2\delta_{S}(\bar x;\nabla F(\bar x)^\top y^*)(u)
		=
		\skalp{y^*,\nabla^2 F(\bar x)(u,u)} + {\dr}^2\delta_{C}(F(\bar x);y^*)(\nabla F(\bar x)u).
	\]
\end{proposition}

On the basis of \cref{The:Marginal_2}, we now aim to derive an image rule.
In order to do that, we assume that $S:=G(Q):=\{x\in\R^n\,|\,\exists y\in Q\colon\,G(y)=x\}$
holds for a twice continuously differentiable mapping $G\colon\R^m\to\R^n$ and a closed set
$Q\subset\R^m$. Associated with the data is the set-valued mapping $\Phi\colon\R^n\tto\R^m$ given by
\begin{equation}\label{eq:der_Phi}
	\Phi(x):= Q\cap G^{-1}(x)\qquad\forall x\in\R^n.
\end{equation}
Note that $\dom\Phi=S$.

\begin{proposition}\label{prop:second_subderivative_image}
	Consider a closed set $Q \subset \R^m$ and a twice continuously differentiable mapping 
	$G\colon \R^m \to \R^n$, and let $\bar x \in S:=G(Q)$ be chosen
	such that $\Phi$ is inner semicompact at $\bar x$ w.r.t.\ $S$.
	Then for each $x^*\in\R^n$ and $u\in\R^n$ as well as each pair $(y,v) \in \R^m \times \R^m$ with 
	$G(y) = \bar x$ and $\nabla G(y) v = u$, one has
	\begin{equation}\label{eq:forward_rule_image}
		{\dr}^2\delta_{S}(\bar x;x^*)(u)
		\leq
		- \skalp{x^*,\nabla^2 G(y)(v,v)} + {\dr}^2\delta_{Q}(y;\nabla G(y)^{\top} x^*)(v),
	\end{equation}
	and the converse holds always true if, additionally, $u\notin T_S(\bar x)$ or $\skalp{x^*,u}<0$
	in which cases both sides in \eqref{eq:forward_rule_image} equal $\infty$.
	
	Given $u \in T_S(\bar x)$ such that $\skalp{x^*,u}= 0$,
	assume that the mapping $\Phi$ is
	inner calm* at $\bar x$ in direction $u$ w.r.t.\ $S$.
	Then we have
	\[
		{\dr}^2\delta_{S}(\bar x;x^*)(u)
		=
		\inf\limits_{G(y)=\bar x,\,\nabla G(y) v = u}
		\bigl(-\skalp{x^*,\nabla^2 G(y)(v,v)} + {\dr}^2\delta_{Q}(y;\nabla G(y)^{\top} x^*)(v)\bigr).
	\]
	Furthermore, if ${\dr}^2\delta_{S}(\bar x;x^*)(u)$ is finite,
	then there exists a pair $(y,v) \in \R^m \times \R^m$ with $G(y) = \bar x$ and $\nabla G(y) v = u$
	such that the estimate \eqref{eq:forward_rule_image} holds as equality, i.e., the 
	infimum is attained.
\end{proposition}
\begin{proof}
    Setting $\varphi(x,y) := \delta_{(\R^n \times Q) \cap \gph G^{-1}}(x,y)$ for all $x\in\R^n$ and $y\in\R^m$ yields
    $\delta_{S}(x) = \inf_{y \in \R^m} \varphi(x,y)$, i.e., we may apply \cref{The:Marginal_2} in order to verify the claim
    (note that $\varphi$ is obviously proper, lower semicontinuous, and lower bounded).
    First, let us show that
    \begin{equation}\label{eq:second_subderivative_of_special_indicator}
		{\dr}^2\varphi((\bar x,y);(x^*,0))(u,v)
		=
		- \skalp{x^*,\nabla^2 G(y)(v,v)} + {\dr}^2\delta_{Q}(y;\nabla G(y)^{\top} x^*)(v)
	\end{equation}
    holds for $(y,v)\in\R^m\times\R^m$ satisfying $G(y)=\bar x$ and $\nabla G(y) v = u$.
    
    Taking into account \eqref{eq:indicatorf}, consider sequences $t_k \downarrow 0$ and $(u_k,v_k) \to (u,v)$ such that $y + t_k v_k \in Q$,
	\[
		\bar x + t_k u_k = G(y + t_k v_k) = G(y) + t_k\nabla G(y)v_k + t_k^2/2 \nabla^2 G(y)(v,v) + \oo(t_k^2),
	\]
	i.e., $u_k = \nabla G(y)v_k + t_k/2 \nabla^2 G(y)(v,v) + \oo(t_k) \to \nabla G(y)v$,
	and
    \begin{align*}
		{\dr}^2\varphi((\bar x,y);(x^*,0))(u,v)
		& = 
		\lim_{k \to \infty} -2 \skalp{(x^*,0),(u_k,v_k)}/t_k.
	\end{align*}
	Then we immediately find
	\begin{align*}
		\dr ^2\varphi((\bar x,y);(x^*,0))(u,v)
		& = 
		-\skalp{x^*,\nabla^2 G(y)(v,v)} + \lim_{k \to \infty} -2 \skalp{x^*,\nabla G(y)v_k}/t_k 
		\\
		& \geq 
		- \skalp{x^*,\nabla^2 G(y)(v,v)} + {\dr}^2\delta_{Q}(y;\nabla G(y)^{\top} x^*)(v).
	\end{align*}
		We note that this estimate holds trivially whenever such sequences $t_k\downarrow 0$ and
		$(u_k,v_k)\to(u,v)$ do not exist since this gives
		$(u,v)\notin T_{(\R^n\times Q)\cap\gph G^{-1}}(\bar x,y)$ and, thus,
		$\dr^2\varphi((\bar x,y);(x^*,0))(u,v)=\infty$.
	
	To show the opposite inequality in \eqref{eq:second_subderivative_of_special_indicator}, 
	consider now sequences $t_k \downarrow 0$ and $v_k \to v$ such that $y + t_k v_k \in Q$
	as well as
  	\[
		{\dr}^2\delta_{Q}(y;\nabla G(y)^{\top} x^*)(v) = \lim_{k \to \infty} -2 \skalp{\nabla G(y)^{\top} x^*, v_k}/t_k.
	\]
	Particularly, $G(y)=\bar x \in S$ and $G(y + t_k v_k) =  \bar x + t_k u_k \in S$
	for $u_k\in\R^n$ given by
	\[
		u_k:= (G(y + t_k v_k) - G(y))/t_k = \nabla G(y)v_k + t_k^2/2 \nabla^2 G(y)(v,v) + \oo(t_k^2).
	\]
	Similar computations as above thus yield the opposite inequality, taking into account
	$\skalp{\nabla G(y)^{\top} x^*, v_k} = \skalp{x^*, \nabla G(y) v_k}$ and $u_k \to \nabla G(y) v = u$.
		Clearly, the converse inequality in \eqref{eq:second_subderivative_of_special_indicator} holds
		trivially if there are no such sequences $t_k\downarrow 0$ and $v_k\to v$ since this gives
		$v\notin T_Q(y)$ and, thus, $\dr^2\delta_Q(y;\nabla G(y)^\top x^*)(v)=\infty$.
	    
	In order to apply \cref{The:Marginal_2}, we interrelate the inner
	semicompactness and inner calmness* of $\Upsilon\colon\R^n\times\R\tto\R^m$ given 
	as in \eqref{eq:M_and_S} with the respective properties of $\Phi$ from
	\eqref{eq:der_Phi}.
	By definition of $\varphi$, we find
	\[
		\Upsilon(x,\alpha)
		=
		\begin{cases}
			\{y\in Q\,|\,G(y)=x\}	&	\alpha\geq 0,\\
			\varnothing				&	\alpha<0
		\end{cases}
		\qquad
		\forall x\in\R^n,\,\forall \alpha\in\R.
	\]
	Note that we have $\dom \Upsilon=S\times\R_+$. 
	A simple calculation reveals that inner semicompactness of $\Phi$ at $\bar x$ w.r.t.\ $S$
	gives inner semicompactness of $\Upsilon$ at $(\bar x,0)$ w.r.t.\ $\dom \Upsilon$ (and vice versa).
	Furthermore, inner calmness* of $\Phi$ at $\bar x$ w.r.t.\ $\dom \Phi$ in direction $u$ gives 
	inner calmness* of $\Upsilon$ at $(\bar x,0)$ w.r.t.\ $\dom \Upsilon$ in direction $(u,\skalp{x^*,u})$.
	
	Now, by means of \cref{The:Marginal_2}, the general upper estimate in \eqref{eq:forward_rule_image}
	is valid for each pair $(y,v)\in\R^m\times\R^m$ with $G(y)=\bar x$ and $\nabla G(y)v=u$
	since $\Upsilon$ from above is assumed to be
	inner semicompact at $(\bar x,0)$ w.r.t.\ $\dom\Upsilon$.
	If $u\notin T_S(\bar x)$, we find $v\notin T_Q(y)$ from \cite[Theorem~6.43]{RoWe98}
	and \cref{prop:basic_properties_second_subderivative_indicator}\,\ref{item:equal_pos_infty}
	gives that both sides of the estimate equal $\infty$.
	In case $\skalp{x^*,u}<0$, we also have $\skalp{\nabla G(y)^\top x^*,v}=\skalp{x^*,\nabla G(y)v}=\skalp{x^*,u}<0$
	and, again, \cref{prop:basic_properties_second_subderivative_indicator}\,\ref{item:equal_pos_infty}
	yields that both sides of the estimate equal $\infty$.
	
	Now, we address the converse relation.
	Observing that the marginal function associated with $\varphi$ satisfies $h:=\delta_S$ in our present setting,
	we find $\dr h(\bar x)(u)=0$ for each $u\in T_S(\bar x)$,
	see \cref{ex:proximal_subdifferential_of_indicator}.
	Hence, the final statements follow from \cref{The:Marginal_2} as well.
\end{proof}

\begin{remark}\label{rem:isc_and_ic*_of_Phi}
	We exploit the notation from \cref{prop:second_subderivative_image}.
	In \cite[Section~5.1.3]{BenkoMehlitz2022a}, it has been pointed out that $\Phi\colon\R^n\tto\R^m$
	from \eqref{eq:der_Phi} is
	inner semicompact at some point $\bar x\in S$ w.r.t.\ its domain $S$ 
	whenever there is a neighborhood $U\subset\R^n$ of $\bar x$
	such that $\Phi(U)$ is bounded, and this is trivially satisfied whenever $Q$ is compact.
	
	Additionally, $\Phi$ is isolatedly calm at $(\bar x,y)$ for some $y\in \Phi(\bar x)$
	whenever we have
	\[
		\nabla G(y)v=0,\quad v\in T_Q(y)\quad\Longrightarrow\quad v=0,
	\]
	see \cite[Section~5.1.3]{BenkoMehlitz2022a} again.
	Whenever this is satisfied for all $y\in \Phi(\bar x)$, and if $\Phi$ is inner semicompact
	at $\bar x$ w.r.t.\ $S$, then $\Phi$ is inner calm* at $\bar x$ w.r.t.\ $S$, see
	\cite[Lemma~4.3(ii)]{BenkoMehlitz2022a}.
\end{remark}

\section{Second-order sufficient conditions in constrained optimization}\label{sec:SOSC_geometric}

In this section, we aim to derive second-order sufficient optimality conditions
for the constrained optimization problem
\begin{equation}\label{eq:constraint_problem}\tag{P}
	\min\{f_0(x)\,|\,F(x)\in C\},
\end{equation}
for twice continuously differentiable mappings $f_0\colon\R^n\to\R$ and $F\colon\R^n\to\R^m$
as well as a closed set $C\subset\R^m$.
The feasible set of \eqref{eq:constraint_problem} will be denoted by $S:=F^{-1}(C)$
throughout the section.
Clearly, this could be done on the basis of \cref{prop:SOC_constrained} and the pre-image rule
stated in \cref{prop:second_subderivative_pre-image}.
This approach improves \cite[Theorem~7.1(ii) and Proposition~7.3]{MohammadiMordukhovichSarabi2021}
as it drops a constraint qualification and/or structural assumptions on $C$ while yielding the quadratic growth condition \eqref{eq:quadratic_growth}
associated with \eqref{eq:constraint_problem}
under the same second-order condition (corresponding to \eqref{eq:second_order_condition_for_geometric_constraints} below with $\alpha = 1$).
In \cite[Theorem~3.3]{BeGfrYeZhouZhang1}, however, a result was shown, which uses the milder
second-order condition (see, again, \eqref{eq:second_order_condition_for_geometric_constraints} below), yet gives a stronger statement in terms
of the following notion introduced by Penot in \cite{Penot1998}.
 \begin{definition}\label{def:essential_local_minimizer}
 A point $\bar{x}\in\R^n$ is said to be an essential local minimizer of second order for problem \eqref{eq:constraint_problem} if $\bar x$ is feasible and there exist $\varepsilon>0$ and $\delta>0$ such that
 \begin{equation}\label{f_def_e-l-m-s-o}
    f(x):=\max\{f_0(x)-f_0(\bar{x}), \dist(F(x),C) \}
 \geq
 \varepsilon \norm{x-\bar{x}}^2 \quad \forall x\in \mathbb{B}_\delta(\bar{x}).
 \end{equation}
 \end{definition}
 In the presence of a mild constraint qualification, \eqref{f_def_e-l-m-s-o} is equivalent
 to the quadratic growth condition at $\bar x$, see \cite[Lemma~3.5]{BeGfrYeZhouZhang1}.
 In general, \eqref{f_def_e-l-m-s-o} is stronger;
 in fact, it is equivalent to the quadratic growth of the function $f$ (which vanishes at $\xb$)
 and thus fully characterized by its second subderivative by
\cref{prop:SOC_unconstrained}.
Hence, we estimate the second subderivative of $f$,
showing versatility of the calculus rules in the process,
and fully recover \cite[Theorem~3.3]{BeGfrYeZhouZhang1}.
To this end, for each constant $\alpha\geq 0$, let us introduce the associated Lagrangian function
$L^\alpha\colon\R^n\times\R^m\to\R$ given by
\begin{equation}\label{eq:Lagrangian}
	L^\alpha(x,\lambda):=\alpha f_0(x)+\skalp{\lambda,F(x)}
	\qquad
	\forall x\in\R^n,\,\forall\lambda\in\R^m.
\end{equation}
Furthermore, for each $\bar x\in S$, we introduce the critical
cone to $S$ at $\bar x$ by means of
\begin{equation}\label{eq:critical_cone}
	\mathcal C(\bar x)
	:=
	\{u\in\R^n\,|\,\nabla F(\bar x)u\in T_C(F(\bar x)),\skalp{\nabla f_0(\bar x),u}\leq 0\}.
\end{equation}
We use the multiplier set $\Lambda^\alpha(\bar x)$ defined by
\begin{equation}\label{eq:multiplier_set}
	\Lambda^\alpha(\bar x)
	:=
	\{\lambda\in\R^m\,|\,\nabla_xL^\alpha(\bar x,\lambda)=0\}.
\end{equation}
\begin{theorem}\label{thm:SOSC_pre-image}
	Let $\bar x\in\R^n$ be a feasible point of \eqref{eq:constraint_problem}.
	Assume that for each $u\in\mathcal C(\bar x)\setminus\{0\}$, there are
	$\alpha\geq 0$ and 
	a multiplier $\lambda\in\Lambda^\alpha(\bar x)$ such that
	\begin{equation}\label{eq:second_order_condition_for_geometric_constraints}
		\nabla^2_{xx}L^\alpha(\bar x,\lambda)(u,u)
		+
		\dr^2\delta_C(F(\bar x);\lambda)(\nabla F(\bar x)u)
		>
		0.
	\end{equation}
	Then $\bar{x}$ is an essential local minimizer of second order for \eqref{eq:constraint_problem}.
\end{theorem}
\begin{proof}
    The statement follows from \cref{prop:SOC_unconstrained}
    once we show $\dr^2 f(\xb;0)(u) > 0$ for all $u\in\R^n\setminus\{0\}$,
    where $f$ is defined in \eqref{f_def_e-l-m-s-o}.
	Note that $f=\vecmax \circ G$ for $G\colon \R^n \to \R^2$ given by 
	$G(x):= (G_1(x),G_2(x)) := (f_0(x)-f_0(\bar{x}), \dist(F(x),C))$
	for all $x\in\R^n$.

    If $\dr f(\xb)(u) > \skalp{0,u} = 0$, we get $\dr^2 f(\xb;0)(u) = \infty$ by
    \eqref{eq:trivial_infinity_of_second_subderivative}.
    Thus, let us assume that $\dr f(\xb)(u)\leq 0$ holds for $u\neq 0$.
    The chain rule \cref{The:Composition_nonsmooth}
    yield the existence of $v \in DG(\xb)(u)$
    such that
    \begin{subequations}\label{eq:chain_rule_technical}
		\begin{align}
			\label{eq:chain_rule_technical_first_order}
				\dr f(\bar x)(u)&\geq\dr \skalp{y^*,G}(\bar x)(u)
					+\dr\vecmax((0,0))(v)-\skalp{y^*,v},\\
			\label{eq:chain_rule_technical_second_order}
				\dr^2 f(\xb;0)(u)
				& \geq 
				\dr^2 \skalp{y^*,G}(\xb;0)(u)
					+
				\dr^2 \vecmax ((0,0);y^*)(v)
		\end{align}
	\end{subequations}
	holds \emph{for each} $y^* \in \R^2$
	(which gives well-defined expressions).
	Now, $v=(v_1,v_2) \in DG(\xb)(u)$ means $v_1 = \nabla f_0(\xb) u$ and
	\begin{equation}\label{v_2_estimate}
		v_2 
		\geq 
		\dr \dist(F(\cdot),C)(\xb)(u) 
		\geq 
		\dr \dist(\cdot,C)(F(\xb))(\nabla F(\xb)u) 
		= 
		\dist\bigl(\nabla F(\xb)u,T_C(F(\xb))\bigr) 
		\geq 
		0
	\end{equation}
	by the chain rule from \cref{The:Composition_nonsmooth} (with $y^*:=0$),
	\cite[Example~8.53]{RoWe98}, and the following argument:
	From $v_2\in DG_2(\bar x)(u)$, we find sequences $u_k\to u$, $v_{2k}\to v_2$, and
	$t_k\downarrow 0$ such that $v_{2k}=\dist(F(\bar x+t_ku_k),C)/t_k$ for all
	$k\in\N$, so taking the limit as $k\to\infty$ while respecting the definition of the
	subderivative gives the first estimate.
    Using \eqref{eq:chain_rule_technical_first_order}
    with $y^*=0$, \cref{prop:vecmax}, 
    and Lipschitzianity of the distance function,
    we obtain
	\[
		0 \geq \dr f(\xb)(u) \geq \dr \vecmax ((0,0))(v) = \max \{\nabla f_0(\xb) u,v_2\},
    \]
    which together with \eqref{v_2_estimate} gives $v_2=0$ and $u\in\mathcal C(\bar x)$.
    Now, our assumptions guarantee the existence of $\alpha\geq 0$ 
	and $\lambda\in\Lambda^\alpha(\bar x)$ satisfying
	\eqref{eq:second_order_condition_for_geometric_constraints},
	which also implies that $\alpha$ and $\lambda$ cannot be zero simultaneously.
	Moreover, if $\nabla f_0(\xb) u < 0$,
	\eqref{eq:second_order_condition_for_geometric_constraints}
	together with $\lambda\in\Lambda^\alpha(\bar x)$ also imply $\alpha = 0$, for otherwise
	$\skalp{\lambda, \nabla F(\xb)u} > 0$ 
	and $\dr^2\delta_C(F(\bar x);\lambda)(\nabla F(\bar x)u) = -\infty$ by
	\cref{prop:basic_properties_second_subderivative_indicator}~\ref{item:greater_neg_infty}. 
	Thus, we have $\alpha\,\nabla f_0(\bar x)u=0$.
	We will show that one has
	\begin{equation}\label{eq:main_estimate}
		\dr^2 f(\xb;0)(u)
		\geq
		\frac{1}{\alpha + \norm{\lambda}} 
		\big(
			\nabla^2_{xx}L^\alpha(\bar x,\lambda)(u,u)
			+
			\dr^2\delta_C(F(\bar x);\lambda)(\nabla F(\bar x)u) 
		\big).
	\end{equation}
	Set $\hat\alpha := \alpha/(\alpha + \norm{\lambda})$, 
	note that $1 - \hat\alpha = \norm{\lambda}/(\alpha + \norm{\lambda})$,
	and apply \eqref{eq:chain_rule_technical_second_order} with $y^*:= (\hat\alpha,1 - \hat\alpha)$ 
	as well as \cref{prop:vecmax} to obtain
	\begin{align*}
		\dr^2 f(\xb;0)(u)
		& \geq 
		\dr^2 (\hat\alpha G_1 + (1 - \hat\alpha) G_2)(\xb;0)(u)
		+
		\dr^2 \vecmax ((0,0);(\hat\alpha,1 - \hat\alpha))(\nabla f_0(\xb) u,0)\\
		& = 
		\dr^2 (\hat\alpha G_1 + (1 - \hat\alpha) G_2)(\xb;0)(u).
	\end{align*}
	If $\lambda = 0$, we get $\alpha>0$, $\hat\alpha = 1$, and $\nabla f_0(\xb)u = 0$
	from above, 
	and \eqref{eq:main_estimate} follows from
	\eqref{eq:second_subderivative_smooth_map}.
	Otherwise, if $\lambda\neq 0$, we set $\hat\lambda := \lambda/\norm{\lambda}$ 
	and continue with the sum rule from 
	\cref{prop:sum_rules}~\ref{item:sum_with_smooth_function}
	as well as the homogeneity property stated in \eqref{eq:homogeneity_property}.
	Additionally, applying the chain rule from \cref{The:Composition_nonsmooth} 
	with $y^*:=\hat\lambda$ as well as \cref{Cor:Distance_function} then yields
	\begin{align*}
		\dr^2 (\hat\alpha G_1 + (1 - \hat\alpha) G_2)(\xb;0)(u)
		& \geq 
		\hat\alpha \nabla^2 f_0(\xb)(u,u) 
		+ 
		(1-\hat\alpha) 
			\dr^2 G_2\left(\xb;-\frac{\hat\alpha}{1-\hat\alpha}\nabla f_0(\xb)\right)(u)\\
		& \geq 
		\hat\alpha \nabla^2 f_0(\xb)(u,u) 
		+  
		(1-\hat\alpha) \nabla^2 \skalp{\hat\lambda,F}(\xb)(u,u)\\
		&\quad
		+  
		 (1-\hat\alpha) \dr^2 \dist(\cdot,C)(F(\xb);\hat\lambda)(\nabla F(\bar x)u)\\
		& \geq 
		\frac{1}{\alpha + \norm{\lambda}} \bigl(\nabla^2_{xx}L^\alpha(\bar x,\lambda)(u,u)
		+
		\dr^2\delta_C(F(\bar x);\lambda)(\nabla F(\bar x)u) \bigr),
	\end{align*}
	noting that $\nabla F(\xb)^{\top} \hat\lambda = - \alpha/\norm{\lambda} \nabla f_0(\xb) 
	= - \hat\alpha/(1-\hat\alpha) \nabla f_0(\xb)$
	and $(1-\hat\alpha) \hat\lambda = \lambda/ (\alpha + \norm{\lambda})$
	and taking into account \eqref{eq:homogeneity_property_2}.
	This completes the proof.
\end{proof}

In \cite[Theorem~3.3]{BeGfrYeZhouZhang1},
this result has been proven via a classical contradiction argument, i.e., via the standard approach
to obtain second-order sufficient optimality conditions,
while our proof is direct, relying on the calculus rules.
In general, the calculus-based approach is certainly very convenient
as one just applies the formulas,
but it may not always be the best since it uses artificial steps, which can be accompanied with artificial assumptions.
In case of second subderivatives, however, it seems that the calculus
works very well  since it can handle even the complicated structure of the function $f$ from \eqref{f_def_e-l-m-s-o} 
without adding superfluous requirements or loosing valuable information.

As pointed out in \cite[Proposition~3.4]{BeGfrYeZhouZhang1}, there is some additional information
about the multipliers available in \cref{thm:SOSC_pre-image} which can be distilled from
the estimate \eqref{eq:second_order_condition_for_geometric_constraints}.
\begin{remark}\label{rem:additional_information_multiplier}
	Let the assumptions of \cref{thm:SOSC_pre-image} be valid.
	Then, for each $u\in\mathcal C(\bar x)\setminus\{0\}$, 
	\eqref{eq:second_order_condition_for_geometric_constraints} yields
	$\dr^2\delta_C(F(\bar x);\lambda)(\nabla F(\bar x)u)>-\infty$
	for some $\lambda\in\Lambda^\alpha(\bar x)$ where $\alpha\geq 0$.
	By means of \cref{prop:basic_properties_second_subderivative_indicator}, 
	this immediately yields
	$\lambda\in\Np_C(F(\bar x);\nabla F(\bar x)u)$.
	Due to \eqref{eq:bounds_proximal_prenormal_cone}, we have
	$\skalp{\lambda,\nabla F(\bar x)u}\leq 0$.
	On the other hand, the definitions of $\mathcal C(\bar x)$ and $\Lambda^\alpha(\bar x)$ give 
	\begin{align*}
		0
		\leq
		\skalp{-\alpha\nabla f_0(\bar x),u}
		=
		\skalp{\nabla F(\bar x)^\top\lambda,u}
		=
		\skalp{\lambda,\nabla F(\bar x)u},
	\end{align*}
	and $\skalp{\lambda,\nabla F(\bar x)u}=0$ follows.
	This gives the additional information 
	$\lambda\in\widehat N_C^p(F(\bar x);\nabla F(\bar x)u)$.
\end{remark}

In the following subsections, we discuss various constrained optimization problems
of the form \eqref{eq:constraint_problem}.
More precisely, we investigate three settings differing from each other by the
particular structure of the set $C$:
\begin{itemize}[leftmargin=6em]
 	\item[\cref{sec:disjunctive_programs}:] 
 		$C$ is polyhedral, i.e., it has no curvature 
 		(standard nonlinear or disjunctive programs);
 	\item[\cref{sec:SOSC_second_order_cone}:] 
 		$C$ is curved, but simple (nonlinear second-order cone programs);
 	\item[\cref{sec:structured_geometric_constraints}:] 
 		$C$ is an image of a pre-image of a simple set (bilevel programs, programs with (quasi-) variational inequality constraints).
\end{itemize}

\subsection{Disjunctive programs}\label{sec:disjunctive_programs}

Here, we apply \cref{thm:SOSC_pre-image} to so-called disjunctive programs
where $C:=\bigcup_{i=1}^\ell P_i$ holds for convex polyhedral sets 
$P_1,\ldots,P_\ell\subset\R^m$. Then $S=F^{-1}(C)$ can be used to
represent feasible sets modeled via
complementarity-, cardinality-, switching-, or vanishing-type constraints,
exemplary, but also standard nonlinear optimization problems,
see e.g.\ \cite{BenkoCervinkaHoheisel2019,FlegelKanzowOutrata2007,Mehlitz2019b}
for an introduction to disjunctive programming and suitable references for
more information on the aforementioned subclasses.
For $\bar y\in C$, we make use of the index set
$J(\bar y):=\{i\in\{1,\ldots,\ell\}\,|\,\bar y\in P_i\}$.
Then, for some $\bar x\in S$, we find
\[
	\mathcal C(\bar x)
	=
	\left\{
		w\in\R^n\,\middle|\,
		\nabla F(\bar x)u\in\bigcup\nolimits_{i\in J(F(\bar x))}T_{P_i}(F(\bar x)),
		\skalp{\nabla f_0(\bar x),u}\leq 0
	\right\}
\]
e.g.\ from \cite[Table~4.1]{AubinFrankowska2009}.
Furthermore, for each $u\in\mathcal C(\bar x)$, we make use of the index set
\[
	J(\bar x;u):=\{i\in J(F(\bar x))\,|\,\nabla F(\bar x)u\in T_{P_i}(F(\bar x))\}.
\]
Based on \cref{thm:SOSC_pre-image}, we find the following 
second-order sufficient optimality condition for the associated
optimization problem \eqref{eq:constraint_problem}.

\begin{theorem}\label{cor:SOSC_disjunctive_programs}
	Let $\bar x\in\R^n$ be a feasible point of \eqref{eq:constraint_problem} 
	where $C:=\bigcup_{i=1}^\ell P_i$ holds for convex polyhedral
	sets $P_1,\ldots,P_\ell\subset\R^m$.
	Assume that for each $u\in\mathcal C(\bar x)\setminus\{0\}$, 
	there are $\alpha\geq 0$
	and a multiplier $\lambda\in\Lambda^\alpha(\bar x)\cap 
	\bigcap_{i\in J(\bar x;u)}N_{T_{P_i}(F(\bar x))}(\nabla F(\bar x)u)$ 
	such that
	\begin{equation}\label{eq:second_order_condition_for_disjunctive_constraints}
		\nabla^2_{xx}L^\alpha(\bar x,\lambda)(u,u)
		>
		0.
	\end{equation}
	Then $\bar x$ is an essential local minimizer of second order 
	for \eqref{eq:constraint_problem}.
\end{theorem}
\begin{proof}
	\Cref{thm:SOSC_pre-image} shows that the assertion of the
	corollary is true whenever for each $u\in\mathcal C(\bar x)\setminus\{0\}$,
	we find $\alpha\geq 0$ and some $\lambda\in\Lambda^\alpha(\bar x)$
	satisfying \eqref{eq:second_order_condition_for_geometric_constraints}.
	The assumptions guarantee that, for each $u\in\mathcal C(\bar x)\setminus\{0\}$,
	we find $\alpha\geq 0$ and $\lambda\in	\Lambda^\alpha(\bar x)\cap
	\bigcap_{i\in J(\bar x;u)}N_{T_{P_i}(F(\bar x))}(\nabla F(\bar x)u)$
	with \eqref{eq:second_order_condition_for_disjunctive_constraints}.
	\Cref{prop:second_subderivative_of_indicator_of_unions} shows
	$\dr^2\delta_C(F(\bar x);\lambda)(\nabla F(\bar x)u)\in\{0,\infty\}$
	in this case.
	Thus, \eqref{eq:second_order_condition_for_disjunctive_constraints} implies
	\eqref{eq:second_order_condition_for_geometric_constraints}.
\end{proof}

Let us note that under some additional assumptions and in case where 
$\ell:=1$ and $\alpha:=1$, a second-order
sufficient condition similar to the one from \cref{cor:SOSC_disjunctive_programs} has been
obtained in \cite[Example~13.25]{RoWe98}.
The sufficient conditions from \cite[Theorem~4.1]{ThinhChuongAnh2021} reduce to ours
when applied to the present situation. However, the authors present them in the presence
of an additional qualification condition.
Furthermore, \cref{cor:SOSC_disjunctive_programs} recovers the sufficient conditions
from \cite[Theorem~3.17]{Gfrerer2014}
and \cite[Theorem~6.1]{BeGfrYeZhouZhang1},
taking into account that we have
\begin{align*}
		\widehat N^p_C(F(\bar x);\nabla F(\bar x)u)
		&=
		\widehat N_{T_C(F(\bar x))}(\nabla F(\bar x)u)
		\\
		&
		=
		\widehat N_{\bigcup_{i\in J(F(\bar x))}T_{P_i}(F(\bar x))}(\nabla F(\bar x)u)
		=
		\bigcap\limits_{i\in J(\bar x;u)}N_{T_{P_i}(F(\bar x))}(\nabla F(\bar x)u)
\end{align*}
from \cite[Remark 5.2]{BeGfrYeZhouZhang1} and \cite[formula~(22)]{BenkoCervinkaHoheisel2019}.

Simplicity of the second-order sufficient conditions from \cref{cor:SOSC_disjunctive_programs}
is caused by the fact that, despite being variationally difficult and highly non-convex,
unions of finitely many convex polyhedral sets are not \emph{curved} causing the
second subderivative of the associated indicator function to be zero if finite, 
see \cref{prop:second_subderivative_of_indicator_of_unions}.
In the next section, we consider the situation where $C$ is an instance of
the well-known second-order cone which possesses curvature.

\subsection{Nonlinear second-order cone programming}\label{sec:SOSC_second_order_cone}

Let us take a closer look at a popular situation where the abstract set
$C$ in the setting of \eqref{eq:constraint_problem} is curved. 
Therefore, recall that for a given integer $s\geq 3$, the set 
\begin{equation*}
	\mathcal Q_s:=\{y\in\R^s\,|\,(y_2^2+\ldots+y_s^2)^{1/2}\leq y_1\}
\end{equation*}
is referred to as second-order or ice-cream cone in $\R^s$.
First, let us compute the second subderivative of the indicator function 
associated with a second-order cone.
Therefore, we coin some additional notation as follows:
\[
	\nnorm{v}:=(v_2^2+\ldots+v_s^2)^{1/2},
	\qquad
	\nskalp{y,v}:=\sum_{i=2}^s y_iv_i
	\qquad
	\forall v,y\in\R^s.
\]

\begin{lemma}\label{lem:second_subderivative_indicator_second_order_cone}
	Fix an integer $s\geq 3$. 
	For each $\bar y\in\mathcal Q_s$, $y^*\in\R^s$, and $v\in\R^s$,
	the following assertions hold.
	\begin{enumerate}
		\item For $\bar y\in\inn\mathcal Q_s$, we have
			\[
				\dr^2\delta_{\mathcal Q_s}(\bar y;y^*)(v)
				=
				\begin{cases}
					\infty		& \skalp{y^*,v}<0,\\
					0			& y^*=0,\\
					-\infty		& \text{otherwise.}
				\end{cases}
			\]
		\item For $\bar y:=0$, we have
			\[
				\dr^2\delta_{\mathcal Q_s}(\bar y;y^*)(v)
				=
				\begin{cases}
					\infty	&\skalp{y^*,v}<0 \text{ or } v\notin \mathcal Q_s,\\
					0		&y^*\in-\mathcal Q_s\cap\{v\}^\perp,\,v\in\mathcal Q_s,\\
					-\infty &\text{otherwise.}
				 \end{cases}
			\]
		\item For $\bar y\in\bdry\mathcal Q_s\setminus\{0\}$, we have
		\[
			\dr^2\delta_{\mathcal Q_s}(\bar y;y^*)(v)
			\geq
			\begin{cases}
				\frac{\norm{y^*}}{\norm{\bar y}}(\nnorm{v}^2-v_1^2)	&	y^*=\beta q(\bar y), \beta\geq 0,\,\nskalp{\bar y,v}\leq \bar y_1v_1,\,\skalp{y^*,v}=0,\\
				\infty					&	y^*=\beta q(\bar y) \text{ and }(\nskalp{\bar y,v}>\bar y_1v_1\text{ or }\skalp{y^*,v}<0),\\
				-\infty					&	\text{otherwise,}
			\end{cases}
		\]
		where $q(\bar y):=\nnorm{\bar y}^{-1}\sum_{i=2}^s\bar y_i\mathtt e_i-\mathtt e_1$.
		Furthermore, equality holds in the first two cases.
		\item Whenever $y^*\in N_{\mathcal Q_s}(\bar y)\cap\{v\}^\perp$ holds, we find
		\[
			\dr^2\delta_{\mathcal Q_s}(\bar y;y^*)(v)
			=
			\begin{cases}
				0				&	\bar y\in \inn\mathcal Q_s,\\
				0				&	\bar y=0,\,v\in\mathcal Q_s,\\
				\frac{\norm{y^*}}{\norm{\bar y}}(\nnorm{v}^2-v_1^2)	&	\bar y\in\bdry\mathcal Q_s\setminus\{0\},\,\nskalp{\bar y,v}\leq \bar y_1v_1,\\
				\infty			&	v\notin T_{\mathcal Q_s}(\bar y).
			\end{cases}
		\]
	\end{enumerate}
\end{lemma}
\begin{proof}
	The assertion of the first statement is trivial.
	
	For the proof of the second statement, we observe that, since $\mathcal Q_s$ is a cone,
	we find
	\[
		\dr^2\delta_{\mathcal Q_s}(\bar y;y^*)(v)
		=
		\liminf_{t\downarrow 0,\,v'\to v,\,v'\in\mathcal Q_s}
		\frac{-2\skalp{y^*,v'}}{t}.
	\]
	The cases where $\skalp{y^*,v}< 0$ or $v\notin \mathcal Q_s=T_{\mathcal Q_s}(\bar y)$ and $\skalp{y^*,v}>0$ are clear, so let us assume
	$\skalp{y^*,v}=0$ and $v\in\mathcal Q_s$. In case where $y^*\in-\mathcal Q_s$, we find
	$\skalp{y^*,v'}\leq 0$ for all $v'\in\mathcal Q_s$, and the second
	subderivative obviously vanishes (here, we used that the polar cone 
	of $\mathcal Q_s$ is $-\mathcal Q_s$). Otherwise, there is $\tilde v\in\mathcal Q_s$
	with $\skalp{y^*,\tilde v}>0$. By convexity of $\mathcal Q_s$, we have
	$v+\tilde v/k\in\mathcal Q_s$ for each $k\in\N$, and
	\[
		-2\skalp{y^*,v+\tilde v/k}/(1/k^2)=-2k\skalp{y^*,\tilde v}\to-\infty
	\]
	shows that the second subderivative is $-\infty$.
	
	For the proof of the third statement, we introduce a function $\phi\colon\R^s\to\R$ by
	means of $\phi(y):=\nnorm{y}-y_1$ for all $y\in\R^s$ and observe that, due to
	$\bar y\in\bdry\mathcal Q_s\setminus\{0\}$, we have $\phi(\bar y)=0$ and $\phi$ is twice
	continuously differentiable at $\bar y$ with non-vanishing gradient $q(\bar y)$. 
	Furthermore, we have
	$\mathcal Q_s=\phi^{-1}(\R_-)$, i.e., \cref{prop:second_subderivative_pre-image}
	can be applied to get a formula for the second subderivative.
	More precisely, using
		\[
			Q(\bar y)
			:=
			\frac{1}{\nnorm{\bar y}}
			\diag(\mathtt e_2+\ldots+\mathtt e_s)
			-
			\frac{1}{\nnorm{\bar y}^3}
			\begin{pmatrix}
				0 	&	0	&	\ldots	& 0
				\\
				0	&	\bar y_2^2	&	\ldots	& \bar y_2\bar y_s
				\\
				\vdots	& \vdots	&\ddots	&\vdots
				\\
				0	&	\bar y_s\bar y_2	&	\ldots	&	\bar y_s^2
			\end{pmatrix},
		\]
	we find
	\[
		\dr^2\delta_{\mathcal Q_s}(\bar y;y^*)(v)
		\geq
		\begin{cases}
			\beta Q(\bar y)(v,v)+\dr^2\delta_{\R_-}(0;\beta)(\skalp{q(\bar y),v})	&	y^*=\beta q(\bar y),\\
			-\infty	& \text{otherwise,}
		\end{cases}
	\]
	and equality holds in the first case.
	Taking \cref{ex:second_subderivative_of_indicator_of_negative_real_line} into account, this gives
	\[
		\dr^2\delta_{\mathcal Q_s}(\bar y;y^*)(v)
		\geq
		\begin{cases}
			\beta Q(\bar y)(v,v)&	y^*=\beta q(\bar y), \beta\geq 0,\,\skalp{q(\bar y),v}\leq 0,\,\skalp{y^*,v}=0,\\
			\infty				 &  y^*=\beta q(\bar y) \text{ and }(\skalp{q(\bar y),v}>0\text{ or }\skalp{y^*,v}<0),\\
			-\infty				 & \text{otherwise,}
		\end{cases}
	\]
	and equality holds in the first two cases.
	Finally, let us simplify the expression $\beta Q(\bar y)(v,v)$ in the first of the appearing cases.
	Thus, let us assume that $y^*=\beta q(\bar y)$ for some $\beta\geq 0$ such that $\skalp{q(\bar y),v}\leq 0$
	and $\skalp{y^*,v}=0$. The case where $\beta=0$ is trivial.
	Therefore, let us assume $\beta>0$.
	Then $\skalp{y^*,v}=0$ gives $\nnorm{\bar y}v_1=\sum_{i=2}^s\bar y_iv_i$.
	From $\bar y\in\bdry\mathcal Q_s\setminus\{0\}$, we have $\norm{\bar y}^2=\bar y_1^2+\nnorm{\bar y}^2=2\nnorm{\bar y}^2$.
	Additionally, $y^*=\beta q(\bar y)$ gives $\norm{y^*}^2=\beta^2\norm{q(\bar y)}^2=\beta^2(1+\nnorm{y}^{-2}\sum_{i=2}^s\bar y_i^2)=2\beta^2$.
	Thus, we end up with
	\begin{align*}
		\beta Q(\bar y)(v,v)
		&=
		\beta\left(\frac{\nnorm{v}^2}{\nnorm{\bar y}}-\frac{1}{\nnorm{\bar y}^3}\left(\sum_{i=2}^s\bar y_iv_i\right)^2\right)
		\\
		&=
		\beta\left(\frac{\nnorm{v}^2}{\nnorm{\bar y}}-\frac{v_1^2}{\nnorm{\bar y}}\right)
		=
		\frac{\sqrt 2\beta}{\sqrt 2\nnorm{\bar y}}(\nnorm{v}^2-v_1^2)
		=
		\frac{\norm{y^*}}{\norm{\bar y}}(\nnorm{v}^2-v_1^2),
	\end{align*}
	and the assertion follows from
	\begin{align*}
		\skalp{q(\bar y),v}\leq 0
		\quad\Longleftrightarrow\quad
		-\nnorm{\bar y}v_1+\sum_{i=2}^s\bar y_iv_i\leq 0
		\quad\Longleftrightarrow\quad
		\nskalp{\bar y,v}\leq \bar y_1v_1
	\end{align*}
	since $\bar y_1=\nnorm{\bar y}$. Note that the above equivalent expressions provide  representations of the condition
	$v\in T_{\mathcal Q_s}(\bar y)$ in the given situation.
	
	The final statement follows from the first three.
\end{proof}

In \cite[Theorem~3.1]{HaMoSa17}, it has been shown that $\delta_{\mathcal Q_s}$ is
already so-called twice epi-differentiable, and that its second epi-derivative can be
calculated by means of the formula stated in the final statement of the above lemma. 
However, the proof of \cite[Theorem~3.1]{HaMoSa17} is much
more technical than our proof of \cref{lem:second_subderivative_indicator_second_order_cone}
which, particularly for the difficult case $\bar y\in\bdry\mathcal Q_s\setminus\{0\}$, exploits
the pre-image rule from \cref{prop:second_subderivative_pre-image}.

Note that the simple pre-image structure $\mathcal Q_s=\phi^{-1}(\R_-)$ of the set $\mathcal Q_s$ for $\phi(y)=\nnorm{y}-y_1$ is
clearly valid for all points in $\mathcal Q_s$, but $\phi$ is not differentiable at 
the origin. 

For given integers $m_1,\ldots,m_\ell\geq 3$ such that $m:=m_1+\ldots+m_\ell$,
we consider 
\begin{equation}\label{eq:product_of_second_order_cones}
	C:=\prod_{i=1}^\ell\mathcal Q_{m_i}
\end{equation}
in \eqref{eq:constraint_problem}
with twice continuously differentiable data functions $f_0\colon\R^n\to\R$ and $F\colon\R^n\to\R^m$. 
Furthermore, let $F_i\colon\R^n\to\R^{m_i}$, $i=1,\ldots,\ell$, be the component mappings of $F$ such that
\[
	F(x)\in C
	\quad\Longleftrightarrow\quad
	F_i(x)\in\mathcal Q_{m_i}\quad\forall i\in\{1,\ldots,\ell\}.
\]
Based on \cref{thm:SOSC_pre-image}, we find the following second-order sufficient
optimality conditions for the associated problem \eqref{eq:constraint_problem}.
\begin{theorem}\label{thm:SOSC_second_order_cone_programming}
	Let $\bar x\in\R^n$ be a feasible point of \eqref{eq:constraint_problem} where $C$ is given
	as in \eqref{eq:product_of_second_order_cones}.
	Assume that for each $u\in \mathcal C(\bar x)\setminus\{0\}$, there 
	are $\alpha\geq 0$ and
	multipliers $\lambda_i\in N_{\mathcal Q_{m_i}}(F_i(\bar x))\cap\{\nabla F_i(\bar x)u\}^\perp$, $i=1,\ldots,\ell$, such that
	\begin{align*}
		0&=\alpha\nabla f_0(\bar x)+\sum_{i=1}^\ell\nabla F_i(\bar x)^\top\lambda_i,
		\\
		0&<\left(\alpha\nabla^2f_0(\bar x)+\sum_{i=1}^\ell\sum_{j=1}^{m_i}(\lambda_i)_j\nabla^2(F_i)_j(\bar x)\right)(u,u)
		+
		\sum\limits_{i\in\mathcal I(\bar x)}\frac{\norm{\lambda_i}}{\norm{F_i(\bar x)}}\bigl(\nnorm{\nabla F_i(\bar x)u}^2-(\nabla F_i(\bar x)u)_1^2\bigr),
	\end{align*}
	where $\mathcal I(\bar x):=\{i\in\{1,\ldots,\ell\}\,|\,F_i(\bar x)\in\bdry\mathcal Q_{m_i}\setminus\{0\}\}$.
	Then $\bar x$ is an essential local minimizer of second order 
	for \eqref{eq:constraint_problem}.
\end{theorem}
\begin{proof}
	We note that the convexity of the cones $\mathcal Q_{m_i}$, $i=1,\ldots,\ell$, yields
	\[
		T_C(F(\bar x))
		=
		\prod_{i=1}^\ell T_{\mathcal Q_{m_i}}(F_i(\bar x)).
	\]
	Thus, $u\in\mathcal C(\bar x)$ satisfies $\nabla F_i(\bar x)u\in T_{\mathcal Q_{m_i}}(F_i(\bar x))$, and we find
	\begin{align*}
		\widehat N^p_C(F(\bar x);\nabla F(\bar x)u)
		&=
		N_C(F(\bar x);\nabla F(\bar x)u)
		=
		N_C(F(\bar x))\cap\{\nabla F(\bar x)u\}^\perp
		\\
		&
		=
		\left(\prod_{i=1}^\ell N_{\mathcal Q_{m_i}}(F_i(\bar x))\right)\cap\{\nabla F(\bar x)u\}^\perp
		\supset
		\prod_{i=1}^\ell N_{\mathcal Q_{m_i}}(F_i(\bar x))\cap\{\nabla F_i(\bar x)u\}^\perp.
	\end{align*}
	Keeping \cref{lem:second_subderivative_of_indicator_of_product,lem:second_subderivative_indicator_second_order_cone} 
	as well as the discussion right after \cref{thm:SOSC_pre-image} in mind, 
	the assumptions of \cref{thm:SOSC_second_order_cone_programming} imply validity of the assumptions of \cref{thm:SOSC_pre-image}
	and the assertion follows.
\end{proof}

A related second-order condition, which comprises a non-vanishing term that incorporates the curvature of the second-order cone, 
has been obtained in \cite[Theorem~29]{BonnansRamirez2005}. The latter result, however, has been formulated in the presence 
of a constraint qualification and makes use of the same Lagrange multiplier for each non-vanishing critical direction. 
As shown above, this is not necessary when dealing with sufficient optimality conditions.
Similarly, our result enhances the sufficient optimality condition presented in
\cite[Proposition~2.1]{HangMordukhovichSarabi2022}.

\subsection{Structured geometric constraints}\label{sec:structured_geometric_constraints}

Observe that the set $C$ in \eqref{eq:constraint_problem} on its own could be the image or pre-image of another closed set.
More precisely, it has been mentioned in \cite[Section~1]{BeGfrYeZhouZhang1} that the setting
\begin{equation}\label{eq:image_of_pre-image}
	C:=H(Q),\qquad Q:=\{z\in\R^\ell\,|\,G(z)\in D\}
\end{equation}
for twice continuously differentiable mappings $G\colon\R^\ell\to\R^p$ and $H\colon\R^\ell\to\R^m$
as well as a simple, closed set $D\subset\R^p$ (e.g., a polyhedral set)
is of significant interest since it covers the special situations
where \eqref{eq:constraint_problem} is a bilevel optimization problem or an optimization problem
with (quasi-) variational inequality constraints, see e.g.\ \cite{Dempe2002,FacchneiPang2003,LuoPangRalph1996,OutrataKocvaraZowe1998}.
Note that essentially the same structure has also been recognized in \cite{BenkoMehlitz2021}, where the set $C$ corresponds to the graph of a set-valued
mapping which is a composition of two other mappings, and the intermediate variables were named \emph{implicit variables} therein.
It is well known that the graph of a composition possesses the structure \eqref{eq:image_of_pre-image}, 
see also \cref{ex:composition_of_set_valued_maps} below.

These considerations underline the need not only for the standard pre-image rule from \cref{prop:second_subderivative_pre-image},
but also for the image rule from \cref{prop:second_subderivative_image}, which is valid under the comparatively mild inner calmness* assumption.
The challenging setting from \eqref{eq:image_of_pre-image} will be investigated deeply in a larger context in a forthcoming paper 
by the authors of \cite{BeGfrYeZhouZhang1} so we do not provide many details here.
Instead, we just want to emphasize that the calculus for second subderivatives is so satisfying that it allows to handle
even such challenging structures with ease as long as second-order sufficient optimality conditions are under consideration.
 
Let us apply our results from \cref{sec:image_pre-image_rule} in order to find a lower estimate
for the second subderivative of the indicator function associated with $C$ from \eqref{eq:image_of_pre-image}.

\begin{lemma}\label{lem:second_subderivative_indicator_image_of_pre-image}
	Fix $\bar y\in C$ given in \eqref{eq:image_of_pre-image}, $y^*\in\R^m$, and $v\in\R^m$.
	Define $\Phi\colon\R^m\tto\R^\ell$ by means of 
	\begin{equation}\label{eq:definition_Phi_image_of_preimage}
		\Phi(y):=\{z\in\R^\ell\,|\,G(z)\in D,\,H(z)=y\}\qquad\forall y\in\R^m,
	\end{equation}
	and assume that $\Phi$ is inner semicompact at $\bar y$ w.r.t.\ $C=\dom\Phi$ and, if $\skalp{y^*,v}=0$, inner calm* at $\bar y$ 
	in direction $v$ w.r.t.\ $C$.
	Then we have
	\[
		\dr^2\delta_C(\bar y;y^*)(v)
		\geq
		\inf\limits_{z\in\Phi(\bar y),\nabla H(z)w=v}
		\sup\limits_{\nabla G(z)^\top\eta=\nabla H(z)^\top y^*}
		\vartheta(z,w,y^*,\eta)
	\]
	for the function $\vartheta\colon\R^\ell\times\R^\ell\times\R^m\times\R^p\to\R$ given by
	\begin{equation}\label{eq:definition_varphi_image_of_preimage}
		\vartheta(z,w,y^*,\eta)
		:=
		\skalp{\eta,\nabla^2G(z)(w,w)}-\skalp{y^*,\nabla^2H(z)(w,w)}+\dr^2\delta_D(G(z);\eta)(\nabla G(z)w)
	\end{equation}
	for arbitrary $z,w\in\R^\ell$, $y^*\in\R^m$, and $\eta\in\R^p$.
\end{lemma}
\begin{proof}
	Due to \cref{prop:second_subderivative_image}, the assumptions of the lemma guarantee
	\[
		\dr^2\delta_C(\bar y;y^*)(v)
		=
		\inf\limits_{z\in\Phi(\bar y),\,\nabla H(z)w=v}
		\left(
		-\skalp{y^*,\nabla^2H(z)(w,w)}
		+
		\dr^2\delta_Q(z;\nabla H(z)^\top y^*)(w)
		\right).
	\]
	Now, we can apply \cref{prop:second_subderivative_pre-image} in order to obtain the fully explicit lower estimate.
\end{proof}

Relying on \cref{rem:isc_and_ic*_of_Phi}, we obtain that in case $\skalp{y^*,v}=0$, the necessary
inner calmness* of the mapping $\Phi$ is inherent whenever the qualification condition
\begin{equation}\label{eq:CQ_image_of_pre-image}
	\nabla H(z)w=0,\quad \nabla G(z)w\in T_D(G(z))\quad\Longrightarrow\quad w=0\qquad\forall z\in\Phi(\bar y)
\end{equation}
is valid since we have $T_{Q}(z)=T_{G^{-1}(D)}(z)\subset\nabla G(z)^{-1}T_D(G(z))$ for each $z\in\Phi(\bar y)$
from \cite[Theorem~6.31]{RoWe98}.

Let us note that the assumptions of \cref{lem:second_subderivative_indicator_image_of_pre-image} hold trivially if
$H\colon\R^m\to\R^m$ is continuously invertible at $\bar z:=H^{-1}(\bar y)$ such that $\nabla H(\bar z)$ is regular.
In this case, the given lower estimate for the second subderivative simplifies to
\[
	\dr^2\delta_C(\bar y;y^*)(v)
	\geq
	\sup\limits_{\nabla G(\bar z)^\top\eta=\nabla H(\bar z)^\top y^*}
	\vartheta(\bar z,\bar w,y^*,\eta)
\]
where $\bar w:=\nabla H(\bar z)^{-1}v$. Particularly, whenever $H\colon\R^m\to\R^m$ is given by $H(z):=Az-b$ for all
$z\in\R^m$ where $A\in\R^{m\times m}$ is a regular matrix and $b\in\R^m$ is a vector, i.e., whenever $H$ models just
a change of coordinates, the above formula applies and the term involving the second derivative of $H$ vanishes. 
The situation where $H$ is the identity map corresponds to the setting where $C$ itself is just a pre-image.

\begin{remark}\label{rem:multiplier_mappings}
	Let us note that in the setting discussed in \cite[Section~1, formula~(3)]{BeGfrYeZhouZhang1},
	the mapping $\Phi\colon\R^m\tto\R^\ell$ from \eqref{eq:definition_Phi_image_of_preimage} is closely related to
	the Lagrange multiplier mapping of another given variational problem.
	Inner semicompactness and inner calmness* of such mappings have been shown to be valid under reasonably mild constraint
	qualifications,
	see \cite[Theorem~3.9]{Be19} for details.
\end{remark}

\begin{example}\label{ex:composition_of_set_valued_maps}
	For set-valued mappings $M_1\colon\R^{m_1}\tto\R^{m_2}$ and $M_2\colon\R^{m_2}\tto\R^{m_3}$ with a closed graph, 
	we consider the situation where $C:=\gph(M_2\circ M_1)$ holds.
	Here, $M_2\circ M_1\colon\R^{m_1}\tto\R^{m_3}$ is the composition of $M_1$ and $M_2$ given by
	\[
		(M_2\circ M_1)(y_1):=\bigcup\limits_{y_2\in M_1(y_1)}M_2(y_2)\qquad\forall y_1\in\R^{m_1}.
	\]
	Setting $\ell:=m_1+m_2+m_3$, $m:=m_1+m_3$, $p:=m_1+m_2+m_2+m_3$,
	\[
		H(z_1,z_2,z_3):=(z_1,z_3),\qquad G(z_1,z_2,z_3):=(z_1,z_2,z_2,z_3)\qquad\forall (z_1,z_2,z_3)\in\R^\ell,
	\]
	and $D:=\gph M_1\times\gph M_2$, we find $C=H(G^{-1}(D))$.
	In the present situation, we have $\Phi(y_1,y_3)=\{(y_1,y_2,y_3)\in\R^\ell\,|\,y_2\in M_1(y_1),\,y_3\in M_2(y_2)\}$, i.e.,
	$\Phi$ is closely related to the so-called intermediate mapping $\Theta\colon\R^{m}\tto\R^{m_2}$, given by
	$\Theta(y_1,y_3):=\{y_2\in M_1(y_1)\,|\,y_3\in M_2(y_2)\}$ for all $(y_1,y_3)\in\R^{m}$, 
	which is associated with the composition $M_2\circ M_1$,
	see \cite[Section~5.3]{BenkoMehlitz2022a}. Clearly, for given $(\bar y_1,\bar y_3)\in\dom\Phi$, $\Phi$ is
	inner semicompact at $(\bar y_1,\bar y_3)$ w.r.t.\ $\dom\Phi$ (inner calm* at $(\bar y_1,\bar y_3)$ in
	direction $(v_1,v_3)\in\R^m$ w.r.t.\ $\dom\Phi$) if and only if $\Theta$ is inner semicompact at $(\bar y_1,\bar y_3)$
	w.r.t.\ $\dom\Theta=\dom\Phi$ (inner calm* at $(\bar y_1,\bar y_3)$ in direction $(v_1,v_3)\in\R^m$ w.r.t.\ $\dom\Theta$),
	and the latter is trivially satisfied if $M_1$ is single-valued and continuous (single-valued and calm in direction $(v_1,v_3)$).
	Observe that the qualification condition \eqref{eq:CQ_image_of_pre-image} is implied by
	\[
		(0,w_2)\in T_{\gph M_1}(\bar y_1,y_2),\,(w_2,0)\in T_{\gph M_2}(y_2,\bar y_3)
		\quad\Longrightarrow\quad
		w_2=0
		\qquad\forall y_2\in\Theta(\bar y_1,\bar y_3),
	\]
	and this is obviously satisfied if, for each $y_2\in\Theta(\bar y_1,\bar y_3)$, $M_1$ is isolatedly calm at $(\bar y_1,y_2)$
	or $M_2^{-1}$ is isolatedly calm at $(\bar y_3,y_2)$. 
\end{example}

Finally, we apply \cref{lem:second_subderivative_indicator_image_of_pre-image} in order to find
second-order sufficient optimality conditions for \eqref{eq:constraint_problem} where $C$ is given as 
in \eqref{eq:image_of_pre-image}.

\begin{theorem}\label{thm:SOSC_image_of_preimage}
	Let $\bar x\in\R^n$ be a feasible point of \eqref{eq:constraint_problem} where $C$ is given as in
	\eqref{eq:image_of_pre-image}.
	Furthermore, let $\Phi\colon\R^m\tto\R^\ell$ defined in \eqref{eq:definition_Phi_image_of_preimage} be inner semicompact
	and inner calm* at $F(\bar x)$ w.r.t.\ $C$.
	Assume that for each $u\in\R^n\setminus\{0\}$ satisfying
	\[
		\nabla F(\bar x)u\in\bigcup\limits_{z\in\Phi(F(\bar x))}\{\nabla H(z)w\in\R^m\,|\,\nabla G(z)w\in T_D(G(z))\},
		\qquad
		\skalp{\nabla f_0(\bar x),u}\leq 0,
	\]
	there are $\alpha\geq 0$ and a multiplier $\lambda\in\Lambda^\alpha(\bar x)$ 
	such that,
	for each $z\in\Phi(F(\bar x))$ and $w\in\R^\ell$ satisfying
	$\nabla H(z)w=\nabla F(\bar x)u$, there is a multiplier
	$\eta\in\Sigma(z,w,\lambda):=\{\eta\in\widehat N^p_D(G(z);\nabla G(z)w)\,|\,\nabla G(z)^\top\eta=\nabla H(z)^\top\lambda\}$
	such that
	\[
		 \nabla^2_{xx}L^\alpha(\bar x,\lambda)(u,u)
		 +
		 \vartheta(z,w,\lambda,\eta)
		>
		0, 
	\]
	where $\vartheta$ has been defined in \eqref{eq:definition_varphi_image_of_preimage}.
	Then $\bar x$ is an essential local minimizer of second order 
	for \eqref{eq:constraint_problem}.
\end{theorem}
\begin{proof}
	For the proof, we are going to apply \cref{thm:SOSC_pre-image}.
	Due to \cite[Theorem~4.1]{Be19}, which is
	applicable since $\Phi$ is inner calm* at $F(\bar x)$ w.r.t.\ $C$, we find
	\[
		\mathcal C(\bar x)
		\subset
		\left\{
			u\in\R^n\,\middle|\,
				\begin{aligned}
					&\exists z\in\Phi(F(\bar x)),\,\exists w\in\R^\ell\colon\\
					&\quad\nabla F(\bar x)u=\nabla H(z)w,\,\nabla G(z)w\in T_D(G(z)),\,\skalp{\nabla f_0(\bar x),u}\leq 0
				\end{aligned}
		\right\}.
	\]
	Keeping \cref{lem:second_subderivative_indicator_image_of_pre-image} in mind, 
	the assumptions of the theorem
	guarantee that for each $u\in\mathcal C(\bar x)\setminus\{0\}$,
	there are $\alpha\geq 0$ and a multiplier $\lambda\in\Lambda^\alpha(\bar x)$ such that
	\[
		\nabla^2_{xx}L^\alpha(\bar x,\lambda)(u,u)
		+
		\dr^2\delta_C(F(\bar x);\lambda)(\nabla F(\bar x)u)>0,
	\]
	i.e., \cref{thm:SOSC_pre-image} shows that $\bar x$ is an essential local
	minimizer of second order.
	Note that the above inequality clearly holds if the infimum in \cref{lem:second_subderivative_indicator_image_of_pre-image}
	is attained, but it also holds if it is not, since in that case $\dr^2\delta_C(F(\bar x);\lambda)(\nabla F(\bar x)u) = \infty$.
\end{proof}

	Note that incorporating the directional proximal normal cone $\widehat N^p_D(G(z);\nabla G(z)w)$ 
	into the definition of the multiplier set $\Sigma(z,w,\lambda)$
	in \cref{thm:SOSC_image_of_preimage} is not restrictive.
	For fixed $z\in\Phi(F(\bar x))$, $w\in\R^\ell$, and $\eta\in\R^p$ such that $\nabla H(z)w=\nabla F(\bar x)u$,
	$\lambda\in\Lambda^\alpha(\bar x)$, and $\nabla G(z)^\top\eta=\nabla H(z)^\top\lambda$,
	the implicitly postulated lower estimate $\dr^2\delta_D(G(z);\eta)(\nabla G(z)w)>-\infty$
	already implies $\eta\in\Np_D(G(z);\nabla G(z)w)$, see \cref{prop:basic_properties_second_subderivative_indicator}. 
	By definition of the proximal pre-normal cone, this gives $\skalp{\eta,\nabla G(z)w}\leq 0$. 
	On the other hand, by choice of $u$, we also have
	\begin{align*}
		\skalp{\eta,\nabla G(z)w}
		=
		\skalp{\nabla H(z)^\top\lambda,w}
		=
		\skalp{\lambda,\nabla F(\bar x)u}
		=
		\skalp{-\nabla f_0(\bar x),u}
		\geq 
		0
	\end{align*}
	giving $\eta\in\{\nabla G(z)w\}^\perp$, i.e., $\eta\in\widehat N^p_D(G(z);\nabla G(z)w)$.

\section{Second-order sufficient conditions in composite optimization}\label{sec:SOSC_composite}

Let us investigate the composite minimization problem
\begin{equation}\label{eq:composite_program}\tag{CP}
	\min\{f_0(x)+g(F(x))\,|\,x\in\R^n\}
\end{equation}
where $f_0\colon\R^n\to\R$ and $F\colon\R^n\to\R^m$ are twice continuously
differentiable and $g\colon\R^m\to\overline\R$ is proper and lower semicontinuous.
In comparison with \eqref{eq:constraint_problem}, this is a more
general optimization problem, since \eqref{eq:constraint_problem} results
from \eqref{eq:composite_program} by setting $g:=\delta_C$.
We make use of the function $L\colon\R^n\times\R^m\to\R$ given via
\[
	L(x,\lambda)
	:=
	f_0(x)
	+
	\skalp{\lambda,F(x)}
	\qquad
	\forall x\in\R^n,\,\forall\lambda\in\R^m
\]
and note that this notation is consistent with \eqref{eq:Lagrangian} since
we have $L=L^1$. 
For each $\bar x\in\R^n$ such that $|g(F(\bar x))|<\infty$, it is, thus,
also reasonable to work with the multiplier set $\Lambda(\bar x):=\Lambda^1(\bar x)$,
see \eqref{eq:multiplier_set}, and we define the
so-called critical cone given by means of
\begin{align*}
	\mathcal C(\bar x)
	:=
	\{u\in\R^n\,|\,\skalp{\nabla f_0(\bar x),u}+\dr g(F(\bar x))(\nabla F(\bar x)u)\leq 0\}.
\end{align*}
Again, this is compatible with the definition of the critical cone
from \eqref{eq:critical_cone} in the setting of constrained optimization, 
see \cref{ex:proximal_subdifferential_of_indicator}.
Observe that
\begin{equation}\label{eq:lower_estimate_first_subderivative}
	\dr (f_0+g\circ F)(\bar x)(u)
	\geq
	\skalp{\nabla f_0(\bar x),u}+\dr g(F(\bar x))(\nabla F(\bar x)u)\qquad\forall u\in\R^n
\end{equation}
by \cite[Theorem~10.6, Corollary~10.9]{RoWe98}, 
see \cref{The:Composition} as well.

Based on \cref{prop:SOC_unconstrained} and our chain rule from \cref{The:Composition}, 
we are in position to state second-order sufficient optimality conditions 
for \eqref{eq:composite_program} very easily.
\begin{theorem}\label{thm:SOSC_composite}
	Let $\bar x\in\R^n$ be a point such that $|g(F(\bar x))|<\infty$.
	Furthermore, assume that, for each $u\in\mathcal C(\bar x)\setminus\{0\}$, there exists
	a multiplier $\lambda\in\Lambda(\bar x)$ such that
	\begin{equation}\label{eq:SOSC_composite}
		\nabla^2_{xx}L(\bar x,\lambda)(u,u)
		+
		\dr^2g(F(\bar x);\lambda)(\nabla F(\bar x)u)
		>
		0.
	\end{equation}
	Then the second-order growth condition holds for $f_0+g\circ F$ at $\bar x$. 
	Particularly, $\bar x$ is a strict local minimizer of \eqref{eq:composite_program}.
\end{theorem} 
\begin{proof}
	We will show that $\dr^2(f_0+g\circ F)(\bar x;0)(u) > 0$ for each $u \in \R^n\setminus\{0\}$
	in order to distill the result from \cref{prop:SOC_unconstrained}.
	Thus, fix $u\in\R^n\setminus\{0\}$.
	If $u \notin \mathcal C(\bar x)$, we obtain
	\[
		\dr (f_0+g\circ F)(\bar x)(u)
		\geq
		\skalp{\nabla f_0(\bar x),u}+\dr g(F(\bar x))(\nabla F(\bar x)u)
		>
		0
	\]
	from \eqref{eq:lower_estimate_first_subderivative},
	and \eqref{eq:trivial_infinity_of_second_subderivative} yields $\dr^2(f_0+g\circ F)(\bar x;0)(u) = \infty$.
	For $u \in \mathcal C(\bar x)$, the assumptions of the theorem guarantee the existence
	of a multiplier $\lambda\in\Lambda(\bar x)$ satisfying
	\begin{equation}\label{eq:lower_estimate_SOSC_composite}
		\begin{aligned}
		\dr^2(f_0+g\circ F)(\bar x;0)(u)
		&=
		\nabla^2f_0(\bar x)(u,u)
		+
		\dr^2(g\circ F)(\bar x,-\nabla f_0(\bar x))(u)
		\\
		&\geq
			\nabla^2_{xx}L(\bar x,\lambda)(u,u)+\dr^2g(F(\bar x);\lambda)(\nabla F(\bar x)u)
		>
		0
		\end{aligned}
	\end{equation}
	by \cref{prop:sum_rules}, \cref{The:Composition}, and the definitions of 
	$L$ and $\Lambda(\bar x)$.
\end{proof}

As already seen in \cref{sec:SOSC_geometric}, 
there is additional information hidden in \eqref{eq:SOSC_composite}.
Similarly as in \cref{rem:additional_information_multiplier} and exploiting
\eqref{eq:trivial_estimates_proximal_pre-subdifferential}, one can show that, for
each $u\in\mathcal C(\bar x)\setminus\{0\}$, the multiplier $\lambda\in\Lambda(\bar x)$
which satisfies \eqref{eq:SOSC_composite} is an element of 
$\widehat{\partial}^pg(F(\bar x);\nabla F(\bar x)u)$.

Let us mention that under some additional assumptions and in a more specific setting, 
second-order sufficient optimality conditions for composite optimization problems 
have been shown in 
\cite[Proposition~2.5]{Sarabi2022} and \cite[Exercise~13.26]{RoWe98}.

When comparing \cref{thm:SOSC_pre-image} and \cref{thm:SOSC_composite}, the natural
question arises whether the assertion of \cref{thm:SOSC_composite} remains true when
using the more general Lagrangian $L^\alpha$ and the more general multiplier set
$\Lambda^\alpha(\bar x)$ for some $\alpha\geq 0$, see
\eqref{eq:Lagrangian} and \eqref{eq:multiplier_set}.
The following examples illustrates that this is indeed not the case.

\begin{example}\label{ex:general_Lagrangian_composite}
	We consider \eqref{eq:composite_program} with the data functions given by
	\[
		f_0(x):=-\tfrac12 x^2,\qquad
		F(x):=\tfrac12 x^2,\qquad
		g(x):=x
		\qquad
		\forall x\in\R
	\]
	and fix the point $\bar x:=0$ which, obviously, is a (global) minimizer of
	$f_0+g\circ F$ where the second-order growth condition fails.
	However, for each $u\in\mathcal C(\bar x)\setminus\{0\}=\R\setminus\{0\}$,
	we find $1\in\Lambda^0(\bar x)$,
	$\nabla^2_{xx}L^0(\bar x,1)(u,u)=u^2>0$
	and
	$\dr^2g(F(\bar x);1)(\nabla F(\bar x)u)=\dr^2g(0;1)(0)=0$
	which gives the estimate
	$\nabla^2_{xx}L^0(\bar x,1)(u,u)+\dr^2g(F(\bar x);1)(\nabla F(\bar x)u)>0$.
	Thus, a potential generalization of the second-order condition from
	\cref{thm:SOSC_pre-image} holds, but the second-order growth condition fails.
\end{example}

The following example visualizes the result of \cref{thm:SOSC_composite} in terms
of so-called sparse optimization.

\begin{example}\label{ex:sparse_programs}
	We consider the problem \eqref{eq:composite_program} with $g(y):=\norm{y}_0$, $y\in\R^m$,
	where $\norm{\cdot}_0\colon\R^m\to\R$ is the so-called $\ell_0$-quasi-norm which counts
	the non-zero entries of the argument vector. Thus, in the associated program \eqref{eq:composite_program},
	those points $\bar x\in\R^n$ are preferred that come along with many zero entries in $F(\bar x)$.
	This is of particular interest whenever the potential constraint
	system $F(x)=0$ of equations does not possess a solution.
	
	Introducing $\phi\colon\R\to\R$ by means of
	\[
		\phi(t)
		:=
		\begin{cases}
			0	&	t=0,\\
			1	&	t\neq 0
		\end{cases}
		\qquad
		\forall t\in\R,
	\]
	we have $\norm{y}_0=\sum_{i=1}^m\phi(y_i)$ for each $y\in\R^m$, i.e., the nonsmoothness of $\norm{\cdot}_0$
	is separable and we can exploit the sum rule from \cref{prop:sum_rules}\,\ref{item:separable_sum}
	in order to compute the (first and) second subderivative of $\norm{\cdot}_0$.
	Some easy calculations show
	\[
		\dr\phi(t)(r)
		=
		\begin{cases}
			0		&	t\neq 0\text{ or }t=r=0,\\
			\infty	&	t=0,\,r\neq 0,
		\end{cases}
		\qquad
		\forall t,r\in\R
	\]
	and
	\[
		\dr^2\phi(t;t^*)(r)
		=
		\begin{cases}
			\infty	&	t^*r<0\text{ or }t=0,\,r\neq 0,\\
			-\infty	&	t^*r>0\text{ and }t\neq 0,\\
			0		&	\text{otherwise}
		\end{cases}
		\qquad
		\forall t,t^*,r\in\R.
	\]
	This can be used to show that
	\[
		\mathcal C(\bar x)
		=
		\{u\in\R^n\,|\,\skalp{\nabla f_0(\bar x),u}\leq 0,\,\skalp{\nabla F_i(\bar x),u}=0\,\forall i\in I_0(\bar x)\}
	\]
	where, for some point $\bar x\in\R^n$,
	we use $I_0(\bar x):=\{i\in\{1,\ldots,m\}\,|\,F_i(\bar x)=0\}$ and $F_1,\ldots,F_m\colon\R^n\to\R$ are the component
	functions of $F$.
	
	Thus, whenever for each $u\in\mathcal C(\bar x)\setminus\{0\}$, there is a $\lambda\in\Lambda(\bar x)$ such that
	$\lambda_i\skalp{\nabla F_i(\bar x),u}=0$ for all $i\notin I_0(\bar x)$ and $\nabla^2_{xx}L(\bar x,\lambda)(u,u)>0$, then the second-order
	growth condition holds for $f_0+\norm{\cdot}_0\circ F$ at $\bar x$.
	Note that the above is obviously less restrictive than a standard second-order sufficient optimality condition 
	for the equality-constrained optimization problem
	\[
		\min\{f_0(x)\,|\,F_i(x)=0\,\forall i\in I_0(\bar x)\}
	\]
	since the multiplier $\lambda$ does not necessarily need to vanish on 
	$\{1,\ldots,m\}\setminus I_0(\bar x)$.
\end{example}

\section{Concluding remarks}\label{sec:conclusions}

In this paper, we derived calculus rules for the second subderivative of
lower semicontinuous functions, comprising a composition rule, a marginal
function rule, an image rule, and a pre-image rule. 
Our findings throw some new light on the results from \cite[Section~13]{RoWe98}.
Moreover, we worked out the precise role of the comparatively new inner calmness* property
from \cite{Be19} in the context of the marginal function and image rule.
We introduced the directional proximal subdifferential of
a given lower semicontinuous function which captures the finiteness of
the second subderivative.
Based on the derived results, we were in position to easily obtain second-order
sufficient optimality conditions in constrained and composite optimization
which are given in terms of initial problem data. Exemplary, this has been
illustrated in terms of disjunctive and nonlinear second-order cone optimization.
Let us point out that our findings are applicable to inherently difficult problem classes
such as optimization problems with (quasi-) variational inequality constraints
(like abstract complementarity constraints induced by non-polyhedral convex
cones) or bilevel optimization problems as well.
However, the calculations which are necessary to estimate the appearing
curvature term from below still are a slightly laborious task which
is why we abstained from presenting them here but leave them as 
a promising topic of future research.
Keeping in mind that second-order sufficient optimality conditions
guarantee local fast convergence of diverse numerical solution 
methods like Newton-type or multiplier-penalty-algorithms, it should
be studied whether our new second-order conditions can be employed
beneficially in this area.
Some first steps in this direction have been done recently in
\cite{HangMordukhovichSarabi2022,HangSarabi2021,Sarabi2022}, where the authors investigate
augmented Lagrangian and sequential-quadratic-programming methods 
for nonlinear second-order cone programs and composite optimization
problems with piecewise linear-quadratic nonsmooth terms.
Finally, let us recall that our approach may not be suitable in order to obtain applicable
second-order necessary optimality conditions. 

\paragraph*{Acknowledgments}

We would like to thank the referees for their valuable comments which helped to improve the quality of this paper.

\paragraph*{Funding}

The research of Mat\'u\v{s} Benko was supported by the Austrian Science Fund (FWF) under grant P32832-N
as well as by the infrastructure of the Institute of Computational Mathematics, Johannes Kepler University Linz, Austria.

\end{document}